\documentclass[sn-mathphys,Numbered]{sn-jnl}% Math and Physical Sciences Reference Style
\usepackage{graphicx}%
\usepackage{multirow}%
\usepackage{amsmath,amssymb,amsfonts}%
\usepackage{amsthm}%
\usepackage{mathrsfs}%
\usepackage[title]{appendix}%
\usepackage[svgnames]{xcolor}%
\usepackage{textcomp}%
\usepackage{manyfoot}%
\usepackage{booktabs}%
\usepackage{algorithm}%
\usepackage{algorithmicx}%
\usepackage{algpseudocode}%
\usepackage{listings}%
\usepackage{subfig}
\usepackage{tikz}
\usepackage[export]{adjustbox}
\usepackage{afterpage}
\usepackage{caption}
\usepackage[skip=0.5ex]{subcaption}

\theoremstyle{thmstyleone}%
\theoremstyle{thmstyletwo}%

\theoremstyle{thmstylethree}%

\begin{document}

\title[High-order exponential integration for seismic wave modeling]{High-order exponential integration for seismic wave modeling}

%%=============================================================%%
%% Prefix	-> \pfx{Dr}
%% GivenName	-> \fnm{Joergen W.}
%% Particle	-> \spfx{van der} -> surname prefix
%% FamilyName	-> \sur{Ploeg}
%% Suffix	-> \sfx{IV}
%% NatureName	-> \tanm{Poet Laureate} -> Title after name
%% Degrees	-> \dgr{MSc, PhD}
%% \author*[1,2]{\pfx{Dr} \fnm{Joergen W.} \spfx{van der} \sur{Ploeg} \sfx{IV} \tanm{Poet Laureate} 
%%                 \dgr{MSc, PhD}}\email{iauthor@gmail.com}
%%=============================================================%%

\author*[1]{\fnm{Fernando} \sur{V. Ravelo}}\email{fernanvr@ime.usp.br\footnote{
Fernando V. Ravelo (ORCID: 0000-0003-1867-2123)\\
Martin Schreiber (ORCID: 0000-0002-4430-6779)\\
Pedro S. Peixoto (ORCID: 0000-0003-2358-3221)}}

\author[2,3]{\fnm{Martin} \sur{Schreiber}}\email{martin.schreiber@univ-grenoble-alpes.fr}
\equalcont{These authors contributed equally to this work.}

\author[1]{\fnm{Pedro} \sur{S. Peixoto}}\email{ppeixoto@usp.br}
\equalcont{These authors contributed equally to this work.}

\affil*[1]{\orgdiv{Instituto de Matem\'atica e Estat\'istica}, \orgname{Universidade de S\~ao Paulo}, \orgaddress{\street{R. do Matão, 1010}, \city{S\~ao Paulo}, \postcode{05508-090}, \state{S\~ao Paulo}, \country{Brazil}}}

\affil[2]{\orgdiv{Laboratoire Jean Kuntzmann / Inria}, \orgname{Université Grenoble Alpes}, \orgaddress{\street{place du Torrent, 150}, \city{Saint-Martin-d'Hères}, \postcode{38400}, \state{Isère}, \country{France}}}

\affil[3]{ \orgname{Technical University of Munich}, \orgaddress{\street{Arcisstraße, 21}, \city{München}, \postcode{80333}, \state{München}, \country{Germany}}}

%%==================================%%
%% sample for unstructured abstract %%
%%==================================%%

\abstract{Seismic imaging is a major challenge in geophysics with broad applications. It involves solving wave propagation equations with absorbing boundary conditions (ABC) multiple times. This drives the need for accurate and efficient numerical methods. This study examines a collection of exponential integration methods, known for their good numerical properties on wave representation, to investigate their efficacy in solving the wave equation with ABC. The purpose of this research is to assess the performance of these methods.
We compare a recently proposed Exponential Integration based on Faber polynomials with well-established Krylov exponential methods alongside a high-order Runge-Kutta scheme and low-order classical methods.
Through our analysis, we found that the exponential integrator based on the Krylov subspace exhibits the best convergence results among the high-order methods.
We also discovered that high-order methods can achieve computational efficiency similar to lower-order methods while allowing for considerably larger time steps.
Most importantly, the possibility of undertaking large time steps could be used for important memory savings in full waveform inversion imaging problems.}


%%\pacs[JEL Classification]{D8, H51}


\maketitle

\section{Introduction}\label{sec_intro}

The resolution of wave propagation equations is a widely researched topic due to its broad range of applications in various fields. One particularly prominent application is seismic imaging, where material parameters of underground regions are estimated based on seismic data. This technique is extensively utilized in the industry for the exploration and extraction of fossil fuels \citep{ikelle2018introduction}.

The numerical approximation of propagating wave equations is a critical stage in this procedure. Consequently, the complexity of the problem impels the development of novel techniques competitive to the efficiency and accuracy of existing schemes \citep{lee2023consistent,alkhadhr2021modeling,kwon2020analysis}.

The propagation of elastic waves can be described as a linear hyperbolic system of PDEs. Nonetheless, the addition of absorbing boundary conditions to replicate an infinite domain modifies the eigenvalues, and they are no longer purely imaginary. In this context, low-order classical explicit schemes such as the Leap-Frog approximation, fourth-order Runge-Kutta, and similar methods have proven effective. Nevertheless, despite their computational speed, these algorithms require very small time steps to approximate the solution accurately. Consequently, this leads to high memory requirements, which can be a significant challenge in solving inverse problems, which is another crucial step in seismic imaging.

In recent decades, a class of numerical algorithms known as exponential integrators have emerged and demonstrated successful applications in various fields. These algorithms have been effectively utilized in areas such as photonics \citep{pototschnig2009time}, the development of numerical methods for weather prediction \citep{peixoto2019semi}, and the modeling of diverse physical phenomena \citep{loffeld2013comparative}, often surpassing the performance of classical schemes. Exponential integrators are typically employed to preserve favorable dispersion properties while allowing for larger time steps \cite{schreiber2019exponential}.

Exponential integrators can be categorized into two types: one primarily concerned with approximating the exponential (or related $\varphi$-functions) of a large matrix resulting from the spatial discretization of the linear term of a system of PDE, and the other focused on different schemes to approximate the non-linear term \citep{hochbruck2010exponential,mossaiby2015}. In the context of wave propagation equations with absorbing boundary conditions, these equations are primarily governed by the linear term, and a source function replaces the non-linear term with a well-defined analytic representation. This leads to a transformation of the problem, as demonstrated by \cite{al2011computing}, which is a generalization of the work of \cite{sidje1998expokit}, where the problem transforms into calculating the exponential of a slightly enlarged matrix.

The approximation of a matrix exponential has received significant attention \citep{moler2003nineteen, alonso2023euler,acebron2019monte}. Numerous exponential integrators have been developed to address this matrix function calculation \citep{hochbruck2010exponential,al2011computing, lu2003computing}. One notable exponential integrator is based on the Krylov subspace, with several schemes utilizing this approach and demonstrating good performance \citep{sidje1998expokit,gaudreault2021kiops,niesen2009krylov}. Another method relies on rational approximations \citep{al2010new}, which are generally implicit and less suitable for large operators. However, they can be combined with the Krylov method to reduce matrix dimensions \citep{al2011computing}. Another approach utilizes Chebyshev polynomials, an explicit method that can be formulated as a three-term recurrence relation \citep{kole2003solving, bergamaschi2000efficient}. Additionally, there are other methodologies, such as Leja points interpolation \citep{bergamaschi2004relpm, deka2023lexint}, optimized Taylor approximations \citep{bader2019computing}, and contour integrals \citep{schmelzer2006evaluating}, among others.

However, when applied to solve hyperbolic systems, such as the wave equations in heterogeneous media, their performance is poorly understood. To the best of our knowledge, only a limited number of literature publications have focused on methods of practical relevance for this specific problem \citep{zhang2014nearly,kole2003solving,fernando2022}. 

In \cite{zhang2014nearly}, an implicit exponential integrator method is developed, and a comparison with other methods is presented, demonstrating superior results in terms of accuracy and dispersion. However, a notable drawback of the implicit method is its high computational cost for each time step, making it primarily suitable for very stiff problems.

\cite{kole2003solving} proposes an explicit exponential integrator based on Chebyshev polynomial approximations, which achieves high solution accuracy and permits large time steps. Nevertheless, the applicability of Chebyshev polynomials for approximating the solution is limited to cases where the system matrix is symmetric or antisymmetric, preventing the modeling of absorbing boundary conditions. As a result, its usage in seismic applications is constrained.

In previous work \citep{fernando2022}, we explored a generalization of the exponential integrator using Faber polynomials, a variant of Chebyshev polynomials. This approach enabled us to solve the wave equations with absorbing boundary conditions. We found that employing higher approximation degrees in the Faber polynomial-based method allows for increased time step sizes without incurring additional computational costs. Furthermore, the augmented time step approximations exhibit favorable accuracy and dispersion properties.

A notable gap in the existing literature is the absence of experiments comparing high-order methods with classical low-order schemes for solving wave equations with absorbing boundary conditions. Our work fills in this gap by comparing exponential integrators based on Faber polynomials, Krylov subspace projection, and High-order Runge-Kutta with various classical methods. Specifically, we consider classical low-order methods such as Leap-frog, fourth-order and four-stage Runge-Kutta (RK4-4), second-order and three-stage Runge-Kutta (RK3-2), and seventh-order and nine-stage Runge-Kutta (RK9-7). Detailed descriptions of these methods can be found in Section \ref{sec_methods}. The comparison between these algorithms focuses on several key characteristics, including numerical dispersion, dissipation, convergence, and computational cost, which are thoroughly discussed in Sections \ref{sec_disp_diss} and \ref{sec_convergence}. By investigating these aspects, we aim to comprehensively evaluate the different methods and their suitability for solving wave equations with absorbing boundary conditions. Finally, in Section \ref{sec_discussion}, we summarize the main findings and draw conclusive remarks based on our research.

\section{The wave equation}

The execution of finite difference methods when solving a system of partial differential equations depends on the continuum formulation and the approximation of the spatial derivatives \citep{thomas2013numerical}. These factors directly impact the discrete operator used in the computations. This section lays the groundwork for the entire analysis presented in the manuscript. We discuss the fundamental elements defining the discrete spatial operator present in seismic imaging applications. These elements include formulating wave propagation equations with absorbing boundary conditions (ABC), spatial discretization using derivative approximations, and free surface treatment.

We employ Perfectly Matching Layers (PML) as the absorbing boundary condition \citep{assi2017compact} to simulate an infinite domain. Despite the significant computational cost associated with PML absorbing boundaries, they remain widely used in numerous numerical studies within the field of seismic imaging \citep{tago2012modelling,jing2019highly,chern2019reflectionless}. For computational efficiency, we implement the PML for the two-dimensional acoustic wave propagation equations. While we can extend our analysis to propagating waves in three dimensions, the complexity of the equations substantially increases, resulting in a significant rise in computational requirements. Thus, for our purposes, we define the system of equations within a rectangular domain $\Omega=[0,a]\times[0,-b]$ for $t>0$, as follows:

\begin{equation}
\frac{\partial}{\partial t}\begin{pmatrix}
    u\\v\\w_x\\w_y
    \end{pmatrix}
        =\begin{pmatrix}
        0&1&0&0\\
        -\beta_x\beta_y+c^2\left(\frac{\partial^2}{\partial x^2}+\frac{\partial^2}{\partial y^2}\right) & -(\beta_x+\beta_y)&c^2\frac{\partial}{\partial x}&c^2\frac{\partial}{\partial y}\\
        (\beta_y-\beta_x)\frac{\partial}{\partial x}& 0& -\beta_x&0\\
        (\beta_x-\beta_y)\frac{\partial}{\partial y}& 0& 0&-\beta_y
        \end{pmatrix}
        \begin{pmatrix}
    u\\v\\w_x\\w_y
    \end{pmatrix}
 	+\begin{pmatrix}
    0\\f\\0\\0
    \end{pmatrix},\label{eq_wave_equation}
\end{equation}
where, $u=u(t,x,y)$ is the displacement, $c=c(x,y)$ is the given velocity distribution in the medium, $v=v(t,x,y)$ is the wave velocity, and $f=f(x,y,t)$ is the source term. The $w$-functions, $(w_x,w_y)=(w_x(t,x,y),w_y(t,x,y))$, are the auxiliary variables of the PML approach and the $\beta$-functions are known and control the damping factor in the absorbing layer.

\begin{equation}
\beta_z(z)=\left\{\begin{array}{ll}
0,&\text{ if } d(z,\partial\Omega)>\delta\\
\beta_0\left(1-\frac{d(z,\Omega)}{\delta}\right)^2,& \text{ if }d(z,\partial\Omega)\leq\delta
\end{array}\right.,\quad z\in\{x,\;y\}
\end{equation}
where $d(z,\partial \Omega)$ is the distance from $z$ to the boundary of $\Omega$, $\delta$ is the thickness of the PML domain, and $\beta_0$ is the magnitude of the absorption factor. Thus, the domain $\Omega$ comprises a physical domain, where the wave propagates normally, and an outer layer of thickness $\delta$ (the domain of the PML), where the waves dampen.

Due to the attenuation of displacement within the PML domain, we opt for a Dirichlet boundary condition (null displacement) along three sides of the rectangular domain $\Omega$. However, this boundary condition does not apply to the top side, as a free-surface boundary condition is more suitable for seismic-imaging simulations. Therefore, on the upper side of $\Omega$, we exclude the PML domain ($\beta_y(y)=0,$ for all $y\in[0,\delta]$), and determine the solution approximation at the upper boundary based on the chosen spatial discretization.

\subsection{Spatial discretization}

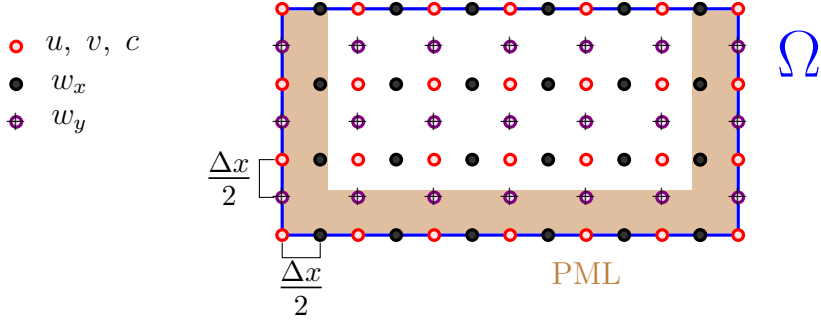
\begin{figure}[bth]
\centering
\begin{tikzpicture}
\draw[blue, very thick] (0,0) -- (0,3) -- (6,3) -- (6,0) -- (0,0);
\node at (6.8,2.4) {\Huge\textcolor{blue}{$\Omega$}};

\fill[fill=brown!50] (0.03,0.6) -- (5.97,0.6) -- (5.97,0.03) -- (0.03,0.03) -- (0.03,0.6);
\fill[fill=brown!50] (0.03,2.97) -- (0.6,2.97) -- (0.6,0.6) -- (0.03,0.6) -- (0.03,2.97);
\fill[fill=brown!50] (5.97,2.97) -- (5.4,2.97) -- (5.4,0.6) -- (5.97,0.6) -- (5.97,2.97);
\node at (4,-0.5) {\large\textcolor{brown}{PML}};

\filldraw[color=red, fill=red!10, very thick](0,0) circle (0.07);
\filldraw[color=red, fill=red!10, very thick](1,0) circle (0.07);
\filldraw[color=red, fill=red!10, very thick](2,0) circle (0.07);
\filldraw[color=red, fill=red!10, very thick](3,0) circle (0.07);
\filldraw[color=red, fill=red!10, very thick](4,0) circle (0.07);
\filldraw[color=red, fill=red!10, very thick](5,0) circle (0.07);
\filldraw[color=red, fill=red!10, very thick](6,0) circle (0.07);
\filldraw[color=black, fill=black!80, very thick](0.5,0) circle (0.07);
\filldraw[color=black, fill=black!80, very thick](1.5,0) circle (0.07);
\filldraw[color=black, fill=black!80, very thick](2.5,0) circle (0.07);
\filldraw[color=black, fill=black!80, very thick](3.5,0) circle (0.07);
\filldraw[color=black, fill=black!80, very thick](4.5,0) circle (0.07);
\filldraw[color=black, fill=black!80, very thick](5.5,0) circle (0.07);
\filldraw[color=black, fill=black!80, very thick](5.5,0) circle (0.07);

\filldraw[color=red, fill=red!10, very thick](0,1) circle (0.07);
\filldraw[color=red, fill=red!10, very thick](1,1) circle (0.07);
\filldraw[color=red, fill=red!10, very thick](2,1) circle (0.07);
\filldraw[color=red, fill=red!10, very thick](3,1) circle (0.07);
\filldraw[color=red, fill=red!10, very thick](4,1) circle (0.07);
\filldraw[color=red, fill=red!10, very thick](5,1) circle (0.07);
\filldraw[color=red, fill=red!10, very thick](6,1) circle (0.07);
\filldraw[color=black, fill=black!80, very thick](0.5,1) circle (0.07);
\filldraw[color=black, fill=black!80, very thick](1.5,1) circle (0.07);
\filldraw[color=black, fill=black!80, very thick](2.5,1) circle (0.07);
\filldraw[color=black, fill=black!80, very thick](3.5,1) circle (0.07);
\filldraw[color=black, fill=black!80, very thick](4.5,1) circle (0.07);
\filldraw[color=black, fill=black!80, very thick](5.5,1) circle (0.07);
\filldraw[color=black, fill=black!80, very thick](5.5,1) circle (0.07);

\filldraw[color=red, fill=red!10, very thick](0,2) circle (0.07);
\filldraw[color=red, fill=red!10, very thick](1,2) circle (0.07);
\filldraw[color=red, fill=red!10, very thick](2,2) circle (0.07);
\filldraw[color=red, fill=red!10, very thick](3,2) circle (0.07);
\filldraw[color=red, fill=red!10, very thick](4,2) circle (0.07);
\filldraw[color=red, fill=red!10, very thick](5,2) circle (0.07);
\filldraw[color=red, fill=red!10, very thick](6,2) circle (0.07);
\filldraw[color=black, fill=black!80, very thick](0.5,2) circle (0.07);
\filldraw[color=black, fill=black!80, very thick](1.5,2) circle (0.07);
\filldraw[color=black, fill=black!80, very thick](2.5,2) circle (0.07);
\filldraw[color=black, fill=black!80, very thick](3.5,2) circle (0.07);
\filldraw[color=black, fill=black!80, very thick](4.5,2) circle (0.07);
\filldraw[color=black, fill=black!80, very thick](5.5,2) circle (0.07);
\filldraw[color=black, fill=black!80, very thick](5.5,2) circle (0.07);

\filldraw[color=red, fill=red!10, very thick](0,3) circle (0.07);
\filldraw[color=red, fill=red!10, very thick](1,3) circle (0.07);
\filldraw[color=red, fill=red!10, very thick](2,3) circle (0.07);
\filldraw[color=red, fill=red!10, very thick](3,3) circle (0.07);
\filldraw[color=red, fill=red!10, very thick](4,3) circle (0.07);
\filldraw[color=red, fill=red!10, very thick](5,3) circle (0.07);
\filldraw[color=red, fill=red!10, very thick](6,3) circle (0.07);
\filldraw[color=black, fill=black!80, very thick](0.5,3) circle (0.07);
\filldraw[color=black, fill=black!80, very thick](1.5,3) circle (0.07);
\filldraw[color=black, fill=black!80, very thick](2.5,3) circle (0.07);
\filldraw[color=black, fill=black!80, very thick](3.5,3) circle (0.07);
\filldraw[color=black, fill=black!80, very thick](4.5,3) circle (0.07);
\filldraw[color=black, fill=black!80, very thick](5.5,3) circle (0.07);
\filldraw[color=black, fill=black!80, very thick](5.5,3) circle (0.07);

\filldraw[color=Purple, fill=Purple!10, very thick](0,0.5) circle (0.07);
\filldraw[color=Purple, fill=Purple!10, very thick](1,0.5) circle (0.07);
\filldraw[color=Purple, fill=Purple!10, very thick](2,0.5) circle (0.07);
\filldraw[color=Purple, fill=Purple!10, very thick](3,0.5) circle (0.07);
\filldraw[color=Purple, fill=Purple!10, very thick](4,0.5) circle (0.07);
\filldraw[color=Purple, fill=Purple!10, very thick](5,0.5) circle (0.07);
\filldraw[color=Purple, fill=Purple!10, very thick](6,0.5) circle (0.07);
\filldraw[color=Purple, fill=Purple!10, very thick](0,1.5) circle (0.07);
\filldraw[color=Purple, fill=Purple!10, very thick](1,1.5) circle (0.07);
\filldraw[color=Purple, fill=Purple!10, very thick](2,1.5) circle (0.07);
\filldraw[color=Purple, fill=Purple!10, very thick](3,1.5) circle (0.07);
\filldraw[color=Purple, fill=Purple!10, very thick](4,1.5) circle (0.07);
\filldraw[color=Purple, fill=Purple!10, very thick](5,1.5) circle (0.07);
\filldraw[color=Purple, fill=Purple!10, very thick](6,1.5) circle (0.07);
\filldraw[color=Purple, fill=Purple!10, very thick](0,2.5) circle (0.07);
\filldraw[color=Purple, fill=Purple!10, very thick](1,2.5) circle (0.07);
\filldraw[color=Purple, fill=Purple!10, very thick](2,2.5) circle (0.07);
\filldraw[color=Purple, fill=Purple!10, very thick](3,2.5) circle (0.07);
\filldraw[color=Purple, fill=Purple!10, very thick](4,2.5) circle (0.07);
\filldraw[color=Purple, fill=Purple!10, very thick](5,2.5) circle (0.07);
\filldraw[color=Purple, fill=Purple!10, very thick](6,2.5) circle (0.07);
\node at (0-0.01,0.51) {+};
\node at (1-0.01,0.51) {+};
\node at (2-0.01,0.51) {+};
\node at (3-0.01,0.51) {+};
\node at (4-0.01,0.51) {+};
\node at (5-0.01,0.51) {+};
\node at (6-0.01,0.51) {+};
\node at (0-0.01,1.51) {+};
\node at (1-0.01,1.51) {+};
\node at (2-0.01,1.51) {+};
\node at (3-0.01,1.51) {+};
\node at (4-0.01,1.51) {+};
\node at (5-0.01,1.51) {+};
\node at (6-0.01,1.51) {+};
\node at (0-0.01,2.51) {+};
\node at (1-0.01,2.51) {+};
\node at (2-0.01,2.51) {+};
\node at (3-0.01,2.51) {+};
\node at (4-0.01,2.51) {+};
\node at (5-0.01,2.51) {+};
\node at (6-0.01,2.51) {+};

\filldraw[color=red, fill=red!10, very thick](-3.5,2.5) circle (0.07);
\filldraw[color=black, fill=black!80, very thick](-3.5,2) circle (0.07);
\filldraw[color=Purple, fill=Purple!10, very thick](-3.5,1.5) circle (0.07);
\node at (-3.5-0.01,1.51) {+};

\node at (-2.5,2.5) {\large $u,\;v,\;c$};
\node at (-2.8,2) {\large $w_x$};
\node at (-2.8,1.5) {\large $w_y$};
\draw[black] (0,-0.1) -- (0,-0.3) -- (0.5,-0.3) -- (0.5,-0.1);
\node at (0.25,-0.7) {\Large$\frac{\Delta x}{2}$};
\draw[black] (-0.1,0.5) -- (-0.3,0.5) -- (-0.3,1) -- (-0.1,1);
\node at (-0.7,0.75) {\Large$\frac{\Delta x}{2}$};
\end{tikzpicture}
\caption{Uniform staggered grid in 2D with the relative positions of the acoustic wave equations' variables and parameters. $u,\;v$ and $c$ are collocated. The shaded region represents the PML domain.}\label{fig_staggered_2D}
\end{figure}

The spatial discretization consists of a uniform staggered grid ($\Delta x=\Delta y$) of 8th-order, guaranteeing the representation of waves up to frequencies of $\frac{2}{\Delta x}$. The positions of the discrete points are depicted in Fig.\,\ref{fig_staggered_2D}.

For the inner discrete points, the 8th-order approximation of the derivatives is given by

\begin{align}
    \frac{\partial u_{i+\frac{1}{2}}}{\partial x}&\approx \frac{1225}{1024\Delta x}\left(u_{i+1}-u_{i}-\frac{u_{i+2}-u_{i-1}}{15}+\frac{u_{i+3}-u_{i-2}}{125}-\frac{u_{i+4}-u_{i-3}}{1715}\right) \label{eq_spatial_8th_1}\\
    \frac{\partial^2 u_i}{\partial x^2}&\approx -\frac{205}{72}u_i+\frac{8}{5}\left(u_{i+1}-u_{i-1}\right)-\frac{1}{5}\left(u_{i+2}-u_{i-2}\right)+\frac{8}{315}\left(u_{i+3}-u_{i-3}\right)\nonumber\\
    &-\frac{1}{560}\left(u_{i+4}-u_{i-4}\right) \label{eq_spatial_8th_2}
\end{align}
with analogous expressions for the $y$-coordinate in the 2D discretization.

The approximation of derivatives near the sides and bottom boundaries, where Dirichlet boundary conditions within a PML domain are applied, is performed using the formulas \eqref{eq_spatial_8th_1} and \eqref{eq_spatial_8th_2}. In these cases, the function values required outside the domain $\Omega$ are set to zero. However, this does not impact the accuracy of the numerical approximations because, within the PML domain, the wave amplitudes decrease to zero.

A different strategy is necessary to approximate derivatives at points near the upper boundary. Since there is no PML domain, and the boundary condition corresponds to a free surface.

\subsection{Free surface}

From the free-surface condition $\nabla u \cdot (0,1) = 0$, we deduce the Neumann boundary condition $\frac{\partial u}{\partial y} = 0$. Additionally, by substituting $\beta_y = 0$ at the free surface in the third equation of \ref{eq_wave_equation}, we obtain $w_y = 0$ at the free surface. Utilizing these values, we can approximate the required spatial derivatives of the functions in \eqref{eq_wave_equation} concerning the variable $y$.

There are two main approaches for approximating the spatial derivatives concerning $y$.
The first approach introduces artificial points outside $\Omega$, assigning function values at these points to satisfy the conditions at the free surface.
The second procedure involves approximating the derivatives at the free surface and its nearest points using only the function values within the domain $\Omega$, without artificially extending the functions. According to \cite{kristek2002efficient}, the latter alternative brings greater accuracy to the numerical solution and is the approach employed throughout this work.

Next, assuming that the grid points lying on the free surface correspond to the evaluation of the displacement $u$ (i.e., the free surface is at $y=0$), we need 8th-order approximations for
\begin{enumerate}
    \item the second derivative $\frac{\partial^2 u}{\partial y^2}$ at the points with $y=\{0,-\Delta x, -2\Delta x,-3\Delta x\}$.
    \item the first derivative $\frac{\partial u}{\partial y}$ at the points with $y=\{-\frac{1}{2}\Delta x,-\frac{3}{2}\Delta x,-\frac{5}{2}\Delta x\}$.
    \item the first derivative $\frac{\partial w_y}{\partial y}$ at the points with $y=\{0,-\Delta x,-2\Delta x,-3\Delta x\}$.
\end{enumerate}

The referred approximations for the derivatives $\frac{\partial^2 u}{\partial y^2}$ and $\frac{\partial u}{\partial y}$, using Taylor expansions, can be found in the Appendix \ref{sec_appendix_approximation_free_surface}. As for the derivative $\frac{\partial w_y}{\partial y}$, we apply the algorithm outlined in \cite{fornberg1988generation}. This algorithm computes the derivative with any approximation order and utilizes an arbitrary points distribution where the values of the derived function are known.

\section{Time integration methods}\label{sec_methods}

After characterizing the spatial discretization and the approximation of the spatial derivatives, we obtain the following linear system of equations:

\begin{equation}\label{eq_1st_ord}
\frac{d}{dt}\boldsymbol{U}(t)=\boldsymbol{H}\boldsymbol{U}(t)+\boldsymbol{f}(t),\quad \boldsymbol{U}(t_0)=\boldsymbol{U}_0.
\end{equation}

Here, $\boldsymbol{U}(t)$ is a vector comprising the discretized functions $u$, $v$, $w_x$, and $w_y$, while the matrix $\boldsymbol{H}$ represents the discretized spatial operator of the system \eqref{eq_wave_equation}. The vector $\boldsymbol{f}$ consists of the source function evaluated at each grid point.

Most of the numerical methods described in this section solve the first-order system of ordinary differential equations \eqref{eq_1st_ord}. Our primary focus lies in approximating the time dimension, leading to the classification of methods as either low or high order concerning time. The following subsections present the numerical schemes employed in the former classifications.

\subsection{Low order methods}

We consider four low-order methods that offer attractive features for approximating the solution of wave equations. Three of these methods are based on the Runge-Kutta (RK) approach, while the fourth is the Leap-frog scheme.

\begin{itemize}
    \item \textbf{2nd order Runge-Kutta (RK3-2):} The RK3-2 method is a second-order RK scheme with three stages. It is a modification of the classical RK2-2 method designed to increase its stability region \citep{crouseilles2020exponential}, enabling its application to hyperbolic problems. The scheme can be expressed as follows:

    \begin{align*}
        \boldsymbol{k_1}&=H\boldsymbol{u}^n+\boldsymbol{f}(t_n),\\
        \boldsymbol{k_2}&=H(\boldsymbol{u}^n+(\Delta t/2)\boldsymbol{k_1})+\boldsymbol{f}(t_n+\Delta t/2),\\
        \boldsymbol{k_3}&=H(\boldsymbol{u}^n+(\Delta t/2)\boldsymbol{k_2})+\boldsymbol{f}(t_n+\Delta t/2),\\
        \boldsymbol{u}^{n+1}&=\boldsymbol{u}^n+\Delta t \boldsymbol{k_3}.
    \end{align*}
    
    \item \textbf{4th order Runge-Kutta (RK4-4):}
    The classical RK4-4 scheme balances stability region and computational requirements \citep{burden2015numerical}.
    
    \item \textbf{7th order Runge-Kutta of nine stages (RK9-7):} This scheme has been specifically constructed for hyperbolic equations and exhibits favorable dispersion properties \citep{calvo1996explicit}.
    
    \item \textbf{Two step method (Leap-frog):}
    The Leap-frog method is highly efficient for solving wave equations. It utilizes two time steps to approximate the second-order time derivative. The equations solved by the Leap-frog method are
    \begin{align*}
        \frac{\partial^2 u}{\partial t^2}&=-\beta_x\beta_y u-(\beta_x+\beta_y)\frac{\partial u}{\partial t}+c^2\left(\frac{\partial^2 u}{\partial x^2}+\frac{\partial^2 u}{\partial y^2}+\frac{\partial\omega_x}{\partial x}+\frac{\partial\omega_y}{\partial y}\right)+f\\
        \frac{\partial w_x}{\partial t}&=-\beta_xw_x+(\beta_y-\beta_x)\frac{\partial u}{\partial x}\\
        \frac{\partial w_y}{\partial t}&=-\beta_yw_y+(\beta_x-\beta_y)\frac{\partial u}{\partial y}
    \end{align*}
    with the discrete approximations in time
    \begin{align*}
        \frac{\partial^2 u^n_i}{\partial t^2}&\approx\frac{u^{n+1}_i-2u^n_i+u^{n-1}_i}{\Delta t^2},\\
        \frac{\partial w_{z_{i+1/2}}^n}{\partial t}&\approx\frac{w_{z_{i+1/2}}^{n+1}-w_{z_{i+1/2}}^{n-1}}{2\Delta t},\mbox{     with  } z\in\{x,y\}.
    \end{align*}
   
\end{itemize}

\subsection{High order methods}

The methods presented in this section are of arbitrary order and utilize exponential integrators based on Faber polynomials, Krylov subspaces, and a high-order Runge-Kutta method.

According to \cite{hochbruck2010exponential}, an exponential integrator approximates the semi-analytic solution of \eqref{eq_1st_ord} using the formula of constant variation
\begin{equation*}\label{eq_const_var}
     \boldsymbol{U}(t)=e^{(t-t_0)\boldsymbol{H}}\boldsymbol{U}_0+\int\limits_{t_0}^te^{(t-\tau)\boldsymbol{H}}\boldsymbol{f}(\tau)d\tau.
\end{equation*}

Expanding the function $\boldsymbol{f}$ in a Taylor series, the solution of \eqref{eq_const_var} can be expressed as the matrix exponential \citep{al2011computing}
\begin{equation}\label{eq_matrix_ampli}
    \boldsymbol{u}(t)=\begin{bmatrix}\boldsymbol{I}_{n\times n}&\boldsymbol{0}\end{bmatrix}e^{(t-t_0)\tilde{\boldsymbol{H}}}\begin{bmatrix}\boldsymbol{u}_0\\\boldsymbol{e}_p\end{bmatrix},
\end{equation}
where $\boldsymbol{e}_p\in\mathbb{R}^p$ is a vector with zeros in its first $p-1$ elements and one in its last element, $\boldsymbol{I}_{n\times n}$ is the identity matrix of dimension $n$, and
\begin{equation*}
    \tilde{\boldsymbol{H}}=\begin{pmatrix}
    \boldsymbol{H}& \boldsymbol{W}\\
    \boldsymbol{0} & \boldsymbol{J}_{p-1}
    \end{pmatrix},
\end{equation*}
where the columns of the matrix $\boldsymbol{W}$ consist of the values of the function $\boldsymbol{f}$ and the approximations of the first $p-1$ derivatives of $\boldsymbol{f}$, and $\boldsymbol{J}_{p-1}$ is a square matrix of dimensions $p\times p$ with ones in the upper diagonal and zeros elsewhere.

Equation \eqref{eq_matrix_ampli} forms the basis for the exponential integrator methods implemented in this research, and the approach used to compute the matrix exponential in \eqref{eq_matrix_ampli} determines each of the following exponential integrators.

\begin{itemize}
    \item \textbf{Faber approximation (FA):}
    This method is an exponential integrator based on Faber polynomials. As presented in \cite{fernando2022}, the exponential approximation is carried on with the three-term recurrence Faber series
    \begin{align*}
\boldsymbol{F}_0(\boldsymbol{H})&=\boldsymbol{I}_{n\times n},\quad \boldsymbol{F}_1(\boldsymbol{H})=\boldsymbol{H}/\gamma-c_0\boldsymbol{I}_{n\times n},\\  \boldsymbol{F}_2(\boldsymbol{H})&=\boldsymbol{F}_1(\boldsymbol{H})\boldsymbol{F}_1(\boldsymbol{H})-2c_1\boldsymbol{I}_{n\times n},\label{eq_faber_pol_0}\\
\boldsymbol{F}_j(\boldsymbol{H})&=\boldsymbol{F}_1(\boldsymbol{H})\boldsymbol{F}_{j-1}(\boldsymbol{H})-c_1\boldsymbol{F}_{j-2}(\boldsymbol{H}),\quad j\geq 3,
\end{align*} 
where the parameters $c_0$ and $c_1$ depend on the eigenvalues distribution of the operator $\boldsymbol{H}$. Then, the solution in the next time instant is expressed as
\begin{equation}
\boldsymbol{u}^{n+1}=\sum\limits_{j=0}^m a_j\boldsymbol{F}_j(\boldsymbol{H})\boldsymbol{u}^n,
\end{equation}
    where $a_j$ are the Faber coefficients.
    
    \item \textbf{Krylov subspace projection (KRY):}
    This method is an exponential integrator utilizing operator projections within the Krylov subspace. Various proposed algorithms involve adaptive time steps and different strategies for generating the subspace basis \citep{gaudreault2021kiops}. However, to ensure an impartial comparison among all the schemes, we opt for the traditional Arnoldi algorithm to establish the vector basis and perform the projection of $H$ \citep{gallopoulos1992efficient}.

\begin{lstlisting}[mathescape=true,caption= Pseudocode of Arnoldi algorithm.,captionpos=b ]
    $\boldsymbol{u}_1=\boldsymbol{u}_0/\|\boldsymbol{u}_0\|_2$
    Do $j$ from 1 to $m$:
        $\boldsymbol{w}=\boldsymbol{H}\boldsymbol{u}_j$
        Do $k$ from 1 to $j$:
            $A_{i,j}=\boldsymbol{w}\cdot\boldsymbol{u}_k$
            $\boldsymbol{w}=\boldsymbol{w}-A_{i,j}\boldsymbol{u}_k$
        $A_{j+1,j}=\|\boldsymbol{w}\|_2$
        $\boldsymbol{u}_{j+1}=\boldsymbol{w}/A_{j+1,j}$
    Then, $e^{\boldsymbol{H}}\boldsymbol{u}_0\approx\|\boldsymbol{u}_0\|_2\boldsymbol{[}\boldsymbol{u}_1|\dots|\boldsymbol{u}_m\boldsymbol{]}e^{\boldsymbol{A}}\boldsymbol{e}_1$
\end{lstlisting}
    
    After constructing the matrix projection $\boldsymbol{A}$, we compute the reduced matrix's exponential using the Padé polynomial approximation method, as outlined in \cite{al2011computing}.

    The Arnoldi algorithm to construct an orthonormal basis is very computationally intensive, and the amount of matrix-vector operations does not represent its actual computational cost. Regarding this subject, the use of non-orthonormal bases has been proposed to greatly reduce this cost \citet{gaudreault2021kiops}. However, this implementation produces more approximation errors than the classical Arnoldi method, at least in its analytic approximation. As we aim to use the classical Krylov method, we employ the Arnoldi algorithm without considering the cost of constructing the Krylov subspace, given the potential to significantly reduce the computational cost.

    \item \textbf{High-order Runge-Kutta (HORK):}
    Runge-Kutta methods are extensively used for solving differential equations \cite{butcher1996history}, and also in combination with exponential integrator schemes \citep{lawson1967generalized,crouseilles2020exponential}.\\
    
    These methods are naturally extended to high-order schemes. They can be explicit and are easy to implement. For this research, we adopt the Runge-Kutta algorithm of arbitrary order proposed by \cite{gottlieb2003strong}, defined by the relation
    \begin{align*}
        \boldsymbol{k}_0&=\boldsymbol{u}^n\\
        \boldsymbol{k}_i&=\left(\boldsymbol{I}_{n\times n}+\Delta t \boldsymbol{H}\right)\boldsymbol{k}_{i-1},\quad i=1\dots m-1\\
        \boldsymbol{k}_m&=\sum\limits_{i=0}^{m-2}\lambda_i\boldsymbol{k}_i+\lambda_{m-1}\left(\boldsymbol{I}_{n\times n}+\Delta t \boldsymbol{H}\right)\boldsymbol{k}_{m-1}\\
        \boldsymbol{u}^{n+1}&=\boldsymbol{k}_m,        
    \end{align*}
    where $\lambda_i$ are the coefficients of the Runge-Kutta and have a straightforward computation.
    According to \cite{gottlieb2003strong}, the Runge-Kutta method exhibits strong stability-preserving properties if the coefficients $\lambda_i$ are non-negative.
\end{itemize}

\subsection{Computational cost and memory usage}\label{sec_cost_memory}

In addition to the accuracy of the numerical solution when discussing the different approaches, we are also interested in their resource consumption. Specifically, we focus on the computational operations required by each algorithm and their utilization of computational memory.

Determining the exact number of computations performed by these methods is a complex task, further complicated by the fact that sparse matrix-vector multiplications are known to be bandwidth-limited in terms of performance \cite{huber2020cache,alappat2022level}.
Therefore, we adopt a simplified model that focuses exclusively on counting the loading and storing of elements. We consider only the matrix-vector operations, as the other vector operations introduce, at most, small variations in the number of operations. Consequently, the cost of each method by time step will be its number of stages or matrix-vector operations (MVOs). Therefore, the overall number of MVOs of a method for computing the solution up to a fixed time $T$ and using a time step size $\Delta t$  can be expressed as:

\begin{equation*}
    \text{N}_\text{op}=\#\text{MVOs}\frac{T}{\Delta t}=\frac{\#\text{MOVs}}{\Delta t} T,
\end{equation*}
where the value of $T$ can be disregarded when comparing the methods since it remains constant within a numerical experiment.

Memory consumption becomes a critical factor when solving the three-dimensional wave equation for seismic imaging applications.
The primary concern is for the inverse problem, where the solution for each time step must be stored to be accessed later. Therefore, the number of time steps required for each method
\begin{equation*}
    N_\text{mem}=\frac{T}{\Delta t},
\end{equation*}
is also an important variable that we will take into account afterward.

\section{Analysis on homogeneous media}\label{sec_disp_diss}

A common challenge arises when utilizing finite difference methods to solve wave equations due to numerical dispersion and dissipation. Numerical dispersion occurs when phase velocities depend on the frequency, leading to distortions in wave signals. On the other hand, numerical dissipation is associated with wave amplitude and is responsible for the emergence of high-frequency waves with small amplitudes in finite difference methods (Section 5.1 of \cite{strikwerda2004finite}).\\

Since the continuous wave equation is non-dispersive and non-dissipative, it is essential to ensure that the numerical methods used to solve it do not introduce excessive dispersion and dissipation. In seismic imaging problems, these issues can lead to significant inaccuracies in estimating the velocity field. Therefore, special attention must be given to identifying and mitigating these errors.\\

In this section, we conduct a comparative analysis of the methods introduced in Section \ref{sec_methods} within the context of a homogeneous velocity field and a single wave signal. We focus on evaluating their dispersion and dissipation errors and examining how these errors depend on the choice of time-step size.

\subsection{Numerical dispersion and dissipation by Fourier transform}

Our analysis investigates numerical dispersion and dissipation by quantifying variations in phase velocities of numerical approximations concerning a reference solution. To achieve this, we conduct a comparison in the frequency domain and estimate velocity changes for each frequency. For this analysis, a Fourier transform is applied to the solution, consisting of a single signal of a Ricker wavelet \citep{ricker1994}. Consequently, we consider a homogeneous medium with a source point and a receiver (a spatial position where the signal is recorded over time).\\

Let $\mathcal{F}{\small \text{appr}}(\omega)$ and $\mathcal{F}{\small \text{ref}}(\omega)$ denote the Fourier transforms of the approximated and reference signals, respectively, with $\omega$ representing the frequency. Thus, we establish the relationship as follows:

\begin{equation*}
\mathcal{F}{\small\text{ref}}(\omega) = e^{k(\omega) + il(\omega)}\mathcal{F}{\small\text{appr}}(\omega),
\end{equation*}
where the real functions $k(\omega)$ and $l(\omega)$ account for the numerical dissipation and dispersion errors, respectively, present in the approximated solution.\\

It is important to note that minimizing dissipation and dispersion errors hinges upon the extent to which the functions $k(\omega)$ and $l(\omega)$ approaches zero.  As the numerical solution is computed at a finite number of time instants, $\omega$ is also limited to a finite range. Then, we calculate the mean of the absolute values of $k(\omega)$ and $l(\omega)$, which can be considered an approximation of the integral of their absolute values. Hereafter, we refer to these metrics as the dissipation and dispersion error.
Furthermore, to mitigate potential numerical errors arising from divisions by small quantities during the computation of dispersion and dissipation errors, we exclusively consider frequencies where the amplitudes in the reference or approximated solutions surpass 1\% of the peak amplitude of the reference solution.

In the next section, we will outline the numerical features of the Ricker signal experiment. Following that, in the subsequent two sections, we will apply the criteria discussed here to assess the numerical dissipation and dispersion errors.

\subsection{Single signal experiment}

The numerical solutions for wave propagation equations are computed in the homogeneous medium $\Omega=[0\text{km},6\text{km}]\times[0\text{km},5\text{km}]$, with a velocity $c=3\text{km/s}$. A Ricker source is placed at position $(3\text{km},4.99\text{km})$ (with a delay of $t_0=0.18\text{s}$), and a receiver is positioned at $(3\text{km},2.5\text{km})$. The time integration is carried out until $T=1.3\text{s}$ without applying any absorbing boundary conditions, as the reflections at the boundary have not yet reached the receiver by the final time. The spatial discretization size used for numerical solutions of the methods is $\Delta x=10\text{m}$, while the reference solution is computed with $\Delta x=2.5\text{m}$ and $\Delta t=0.104\text{ms}$ using the RK9-7 scheme.\\

\begin{figure}[H]
\subfloat[Wave propagation at time $t=0.6s$.]{\includegraphics[scale=0.3]{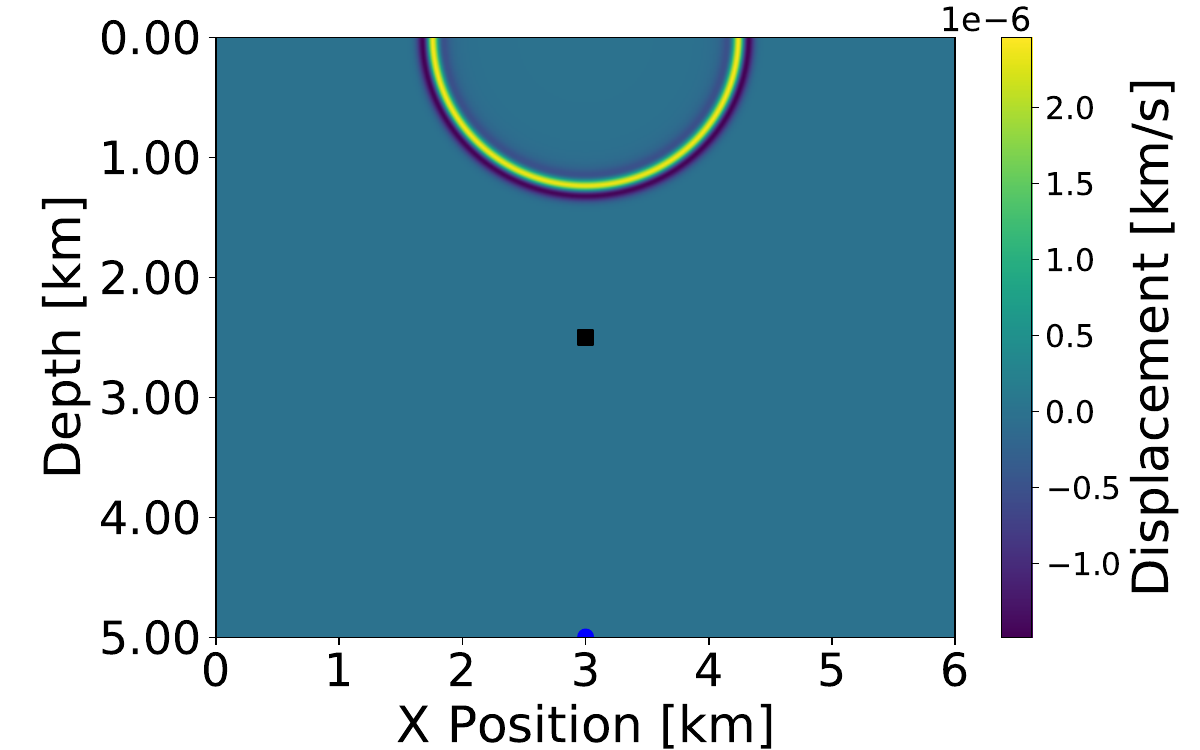}}
    \hfill
\subfloat[Wave propagation at time $T=1.3s$.]{\includegraphics[scale=0.3]{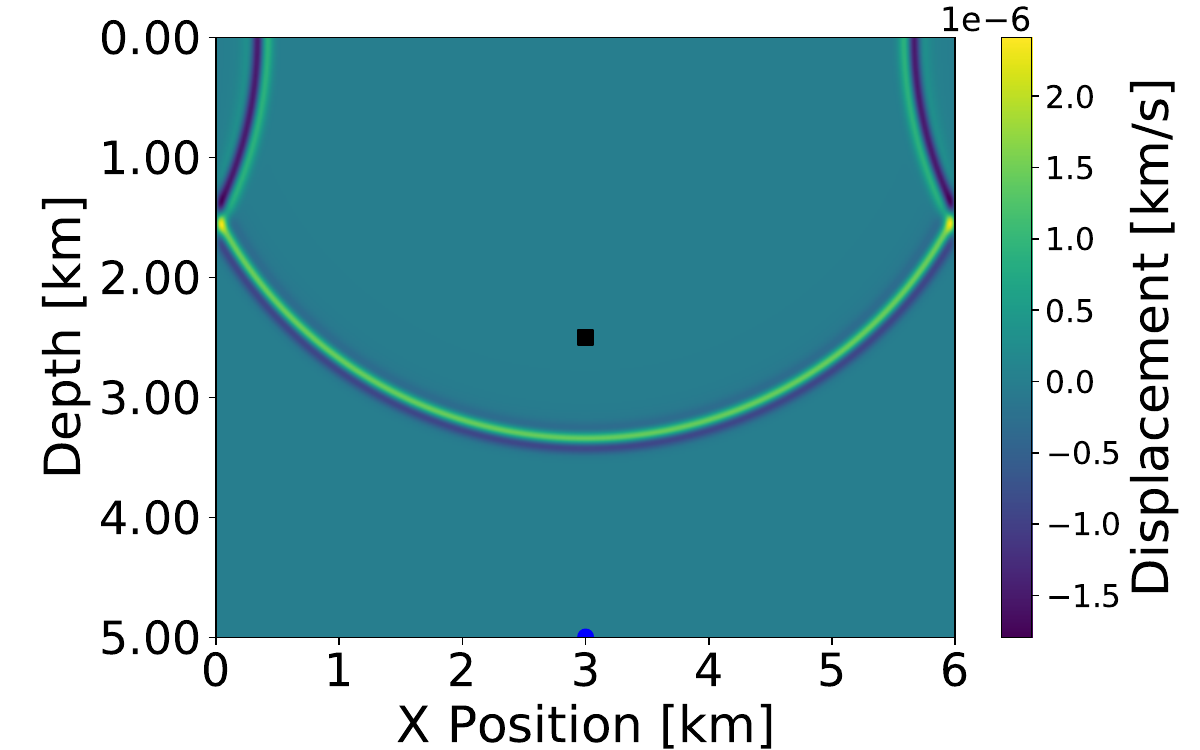}}
\caption{Snapshots of the reference solution at times $t=0.6s$ and $T=1.3s$ within the homogeneous medium $\Omega=[0\text{km},6\text{km}]\times[0\text{km},5\text{km}]$. The Ricker signal source position (blue dot) and the receiver location (black square) are highlighted. During the time interval $t\in[0.6,1.3]s$, the front wave propagates through the receiver location.
}\label{fig_wave_propagation_homo}
\end{figure}

We are mainly interested in the largest time step allowed such that the error of the methods is under a fixed threshold. However, to ensure uniform wave sampling of the numerical approximations at the receiver, we use larger time steps up to the point when the wave closely approaches the receiver ($t=0.6\text{s}$). Then, a uniform $\Delta t=0.417\text{ms}$ is employed until the final time $T=1.3\text{s}$ is reached. Figure \ref{fig_wave_propagation_homo} displays the homogeneous medium with the source and the receiver positions and the snapshots of the reference solution at times $t=0.6s$ and $T=1.3s$.

Although our primary focus lies in evaluating the time error of the methods, it is essential to acknowledge the influence of spatial discretization on numerical accuracy.  To account for this spatial effect, convergence, dispersion, and dissipation are computed for all methods with a small time-step, $\Delta t=0.417\text{ms}$ (see Figure \ref{fig_disp_diss_min_delta} in Appendix \ref{sec_appendix_dips_diss}). The minimum convergence, dispersion, and dissipation errors obtained from this computation serve as an estimation of the spatial effect. Then, we determine the maximum $\Delta t$ allowable for the methods such that the time error remains less or equal to $50\%$ of the spatial error.\\

\begin{figure}[H]
\centering
\includegraphics[scale=0.33]{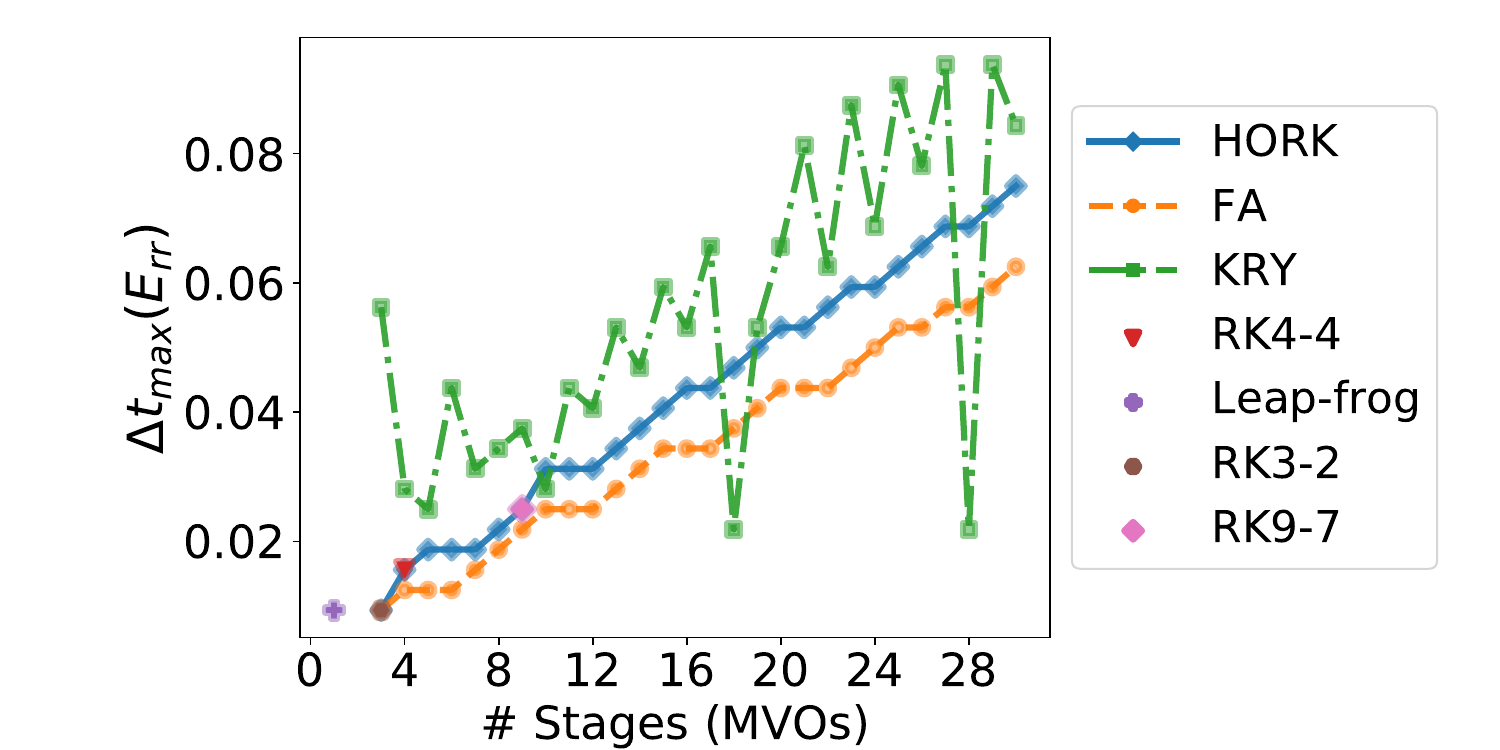}
\caption{Dependence of $\Delta t_{\text{max}}$ on the approximation degree of the numerical scheme. A higher number of stages leads to an increase in the maximum allowable time step without significantly increasing the number of computations.}\label{fig_conv_homo}
\end{figure}

For this experiment, the approximation error due to the spatial discretization is approximately $3.9\cdot 10^{-6}$ (see Figure \ref{fig_disp_diss_min_delta} in Section \ref{sec_appendix_dips_diss}). Based on this, we determine $\Delta t_{\text{max}}$ as the maximum $\Delta t$ such that the approximation error is less or equal to $E_{\text{rr}}=5.9\cdot10^{-6}$. Then, the convergence can be analyzed by investigating the signal error at a specific receiver location $(3\text{km},2.5\text{km})$. It becomes clear that an increase in the number of stages leads to an increase in $\Delta t_{\text{max}}$ (see Figure \ref{fig_conv_homo}).

Referring to Figure \ref{fig_conv_homo}, it can be observed that the Krylov method displays a highly oscillatory pattern concerning its associated $\Delta t_{\text{max}}$. Intriguingly, this pattern reaches its local peak values when the subspace dimension is an odd number. The general behavior of the methods convergence is not sensible to the cutting point of the error threshold, and for variations of $E_{\text{rr}}=5.9\cdot10^{-6}$, they remain valid. So, we expect a similar behavior when studying the dispersion and dissipation.

\subsection{Dispersion results}

The dispersion error arising from spatial discretization is estimated as $0.002$. Consequently, we permit for the time integrator methods an error threshold of $1.5\times$ higher, equating to a maximum allowable dispersion error of $0.003$. Then, we search for $\Delta t_{\text{max}}$ such that the dispersion error remains below this limit.

In addition to $\Delta t_{\text{max}}$, we introduce a computational cost measure denoted as $\text{N}^{\text{\tiny disp}}_{\text{op}}$, similar to the ideas of Section \ref{sec_cost_memory}, defined as:

\begin{equation*}
\text{N}^{\text{\tiny disp}}_{\text{op}}=\frac{\text{\# MVOs}}{\Delta t_{\text{max}}}.
\end{equation*}

\begin{figure}[H]
\centering
\includegraphics[trim=50 480 0 0,clip,scale=0.33]{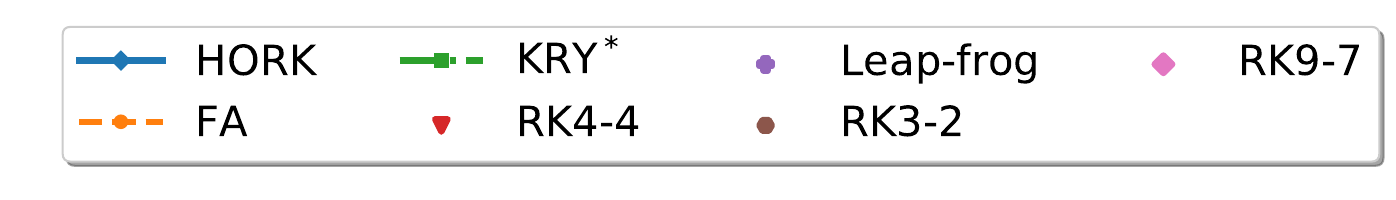}\\[-3ex]
\subfloat[Maximum time-step, $\Delta t_{\text{max}}$.]{\includegraphics[scale=0.33]{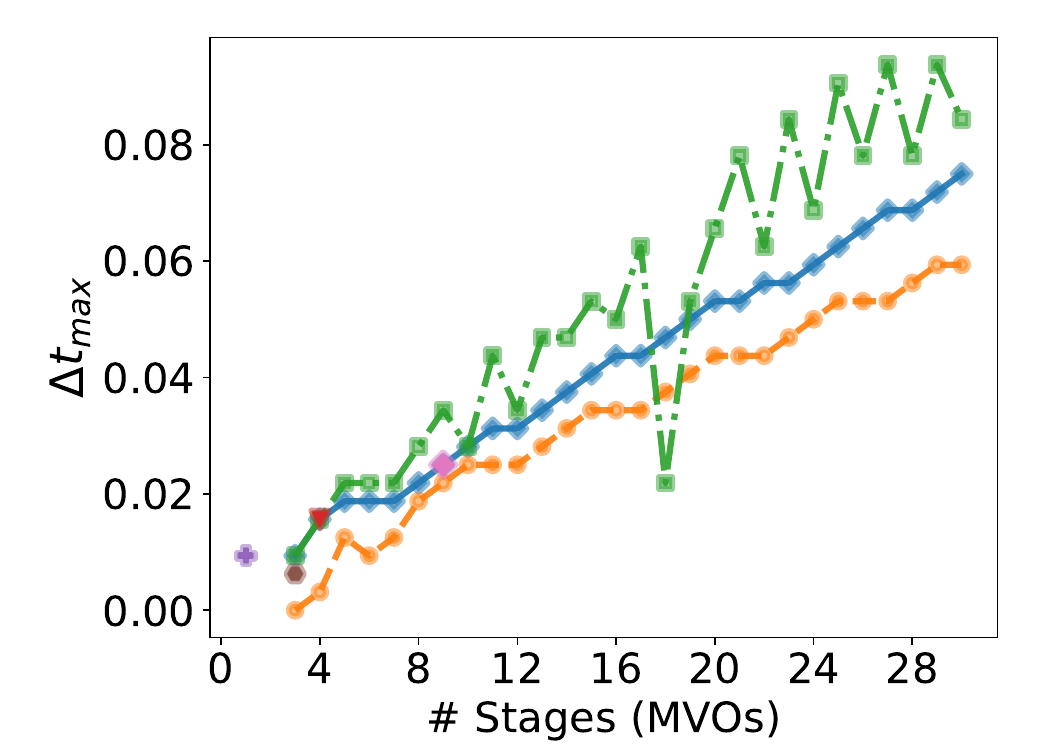}}
    \hfill
\subfloat[Number of MVOs by $\Delta t_{\text{max}}$.]{\includegraphics[scale=0.33]{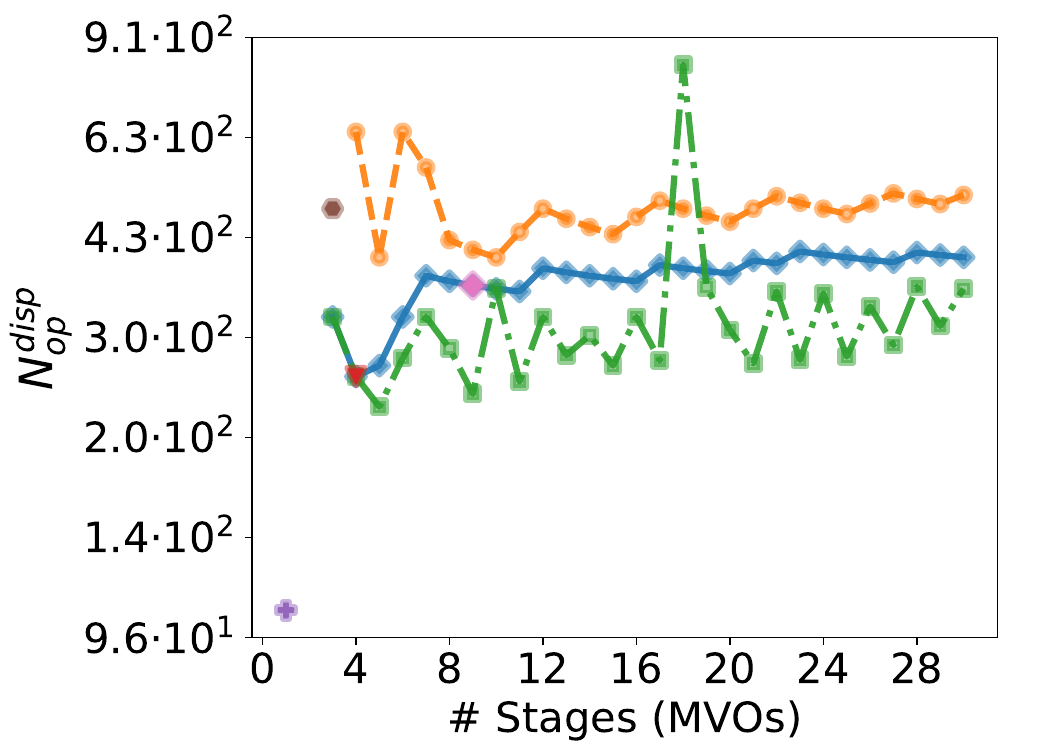}}
\caption{Variation of $\Delta t_{\text{max}}$ (left) and $\text{N}^{\text{\tiny disp}}_{\text{op}}$ (right) concerning the numerical scheme and the number of stages utilized, according to the numerical dispersion error for a Ricker source peak frequency of $f_M=15$Hz. Generally, a higher number of stages leads to an increase in the maximum allowable time step size without significantly increasing the number of computations. * Here we neglect the computational complexity of creating the Krylov subspaces.}\label{fig_disp_max_delta_eff}
\end{figure}

Based on Fig. \ref{fig_disp_max_delta_eff}, the Leap-frog algorithm is approximately two times faster than the other schemes but requires small time steps. On the other hand, the explicit exponential methods exhibit an increase in their maximum time step as the number of stages used rises, without a significant increase in the number of matrix-vector operations required. Interestingly, the peak values of the Krylov methods for the largest $\Delta t$ and the lower $\text{N}_{\text{op}}^\text{disp}$ are consistently for the odd numbers of the subspace dimension \ref{fig_disp_max_delta_eff}.\\

To ensure the robustness of our analysis, we reproduce the previous results in Appendix \ref{sec_appendix_dips_diss_disp} using various peak frequencies of the Ricker source since wave frequencies influence dispersion.

\subsection{Dissipation results}

Similar to the previous section, we estimate the minimum dispersion error, independent of the time integrator used. The minimum dissipation error is approximately $2.4\cdot10^{-7}$. Thus, we once again compute the maximum time-step, $\Delta t_{\text{max}}$, such that the dispersion error remains below $3.6\cdot10^{-7}$. Besides of $\Delta t_{\text{max}}$, we define the computational cost measure as

\begin{equation*}
\text{N}^{\text{\tiny diss}}_{\text{op}}=\frac{\text{\# MVOs}}{\Delta t_{\text{max}}},
\end{equation*}
similar to convergence and dispersion. 

\begin{figure}[H]
\centering
\includegraphics[trim=50 480 0 0,clip,scale=0.33]{figures/legend_b.pdf}\\[-3ex]
\subfloat[Maximum time-step, $\Delta t_{\text{max}}$.]{\includegraphics[scale=0.33]{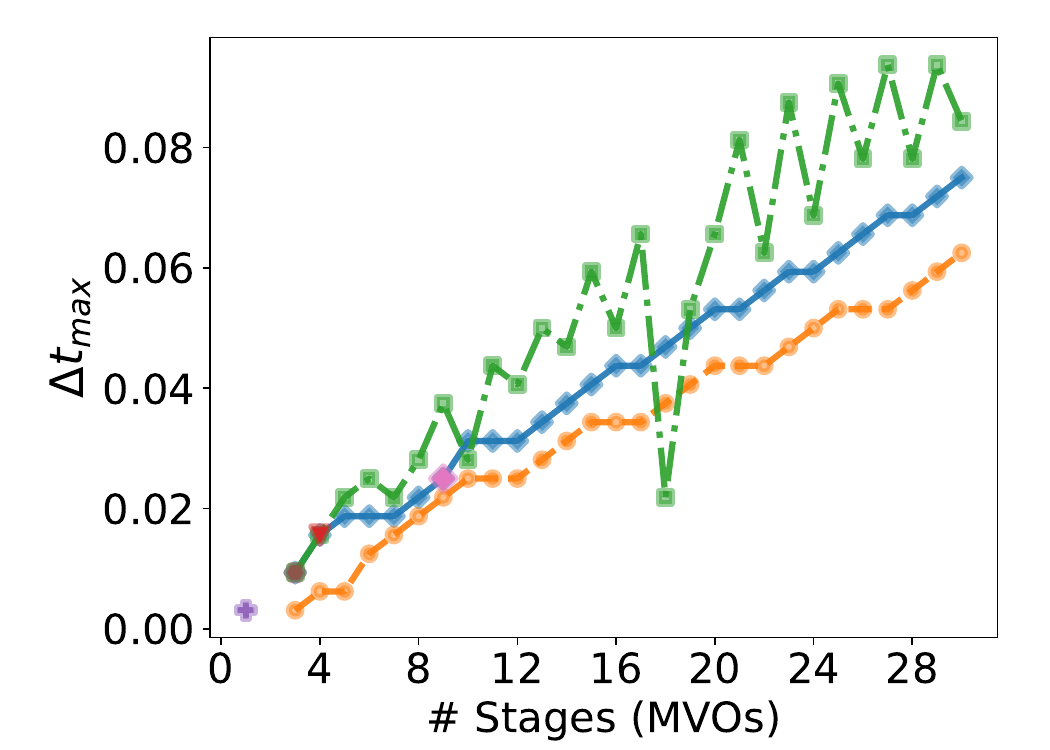}}
    \hfill
\subfloat[Number of MVOs by $\Delta t_{\text{max}}$.]{\includegraphics[scale=0.33]{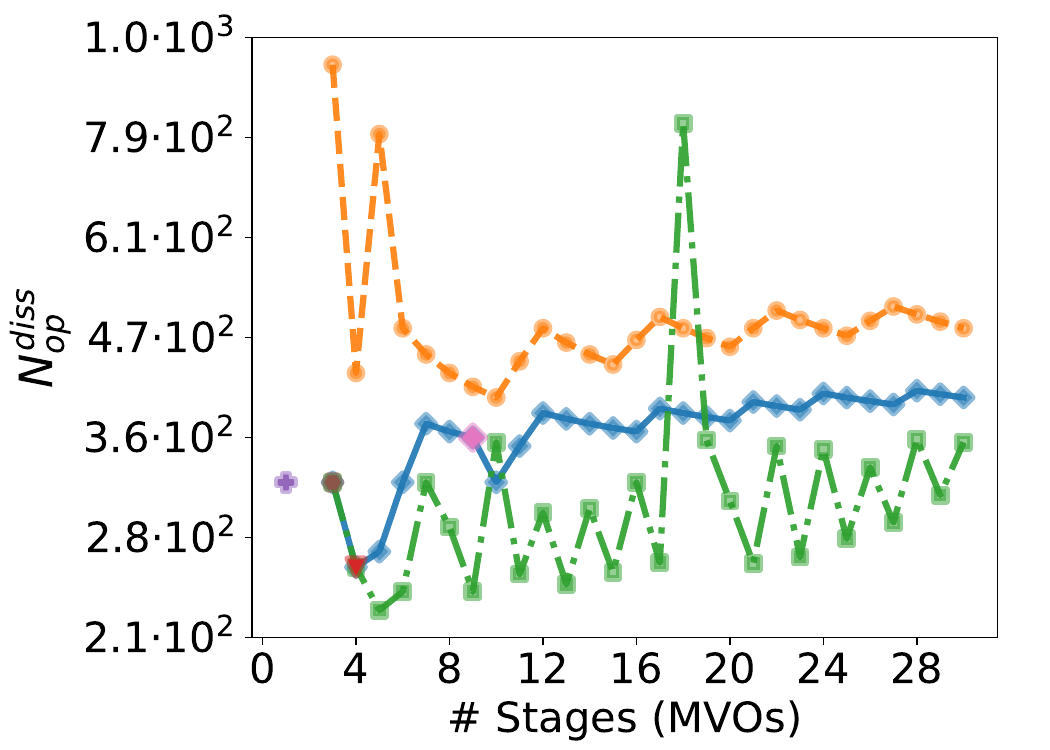}}
\caption{Variation of $\Delta t_{\text{max}}$ (left) and $\text{N}^{\text{\tiny disp}}_{\text{op}}$ (right) concerning the numerical scheme and the number of stages utilized, according to the numerical dissipation error for a Ricker source peak frequency of $f_M=15$Hz. Generally, a higher number of stages leads to an increase in the maximum allowable time step size without significantly increasing the number of computations. * Here we neglect the computational complexity of creating the Krylov subspaces.}\label{fig_diss_max_delta_eff}
\end{figure}

In Figure \ref{fig_diss_max_delta_eff}, a similar trend is observed with the dispersion error, except that the performance of the exponential integrator is better in relation to the Leap-frog when comparing the dissipation error.
Notably, the high-order methods display an increase in the time-step size with the number of stages used without significantly increasing the number of matrix-vector operations required.

As with the numerical dispersion, we reproduce the experiments for different Ricker source peak frequencies in Appendix (Section \ref{sec_appendix_dips_diss_diss}).

\section{Analysis on realistic seismic models}\label{sec_convergence}

In this section, we describe the numerical experiments we will use to compare the accuracy of the approximations of the different methods. For comparison, we generated a reference solution using the RK9-7 scheme with a finer grid ($\Delta x=5$m) and then estimated the error for each method using various time step sizes. To ensure a robust accuracy assessment, we employ two procedures. First, we compare the approximated solution across the entire physical space (excluding the PML domain) at a specific time instant. Second, we compare the seismogram data of the solution values at the upper boundary for all the simulation time. For each error evaluation, we determine the maximum time step size, $\Delta t_{\text{max}}$, that allows a scheme of a particular order to achieve a solution accuracy below a predefined threshold with the least number of MVOs. Additionally, we introduce an efficiency measure and an indicator of memory utilization derived from the number of MVOs and $\Delta t_{\text{max}}$, following the concepts outlined in Section \ref{sec_cost_memory}.

\subsection{Test cases}\label{sec_test_cases}

\begin{figure}[htb!]
\hspace*{0.3cm}\subfloat[Corner Model.]{\includegraphics[scale=0.3]{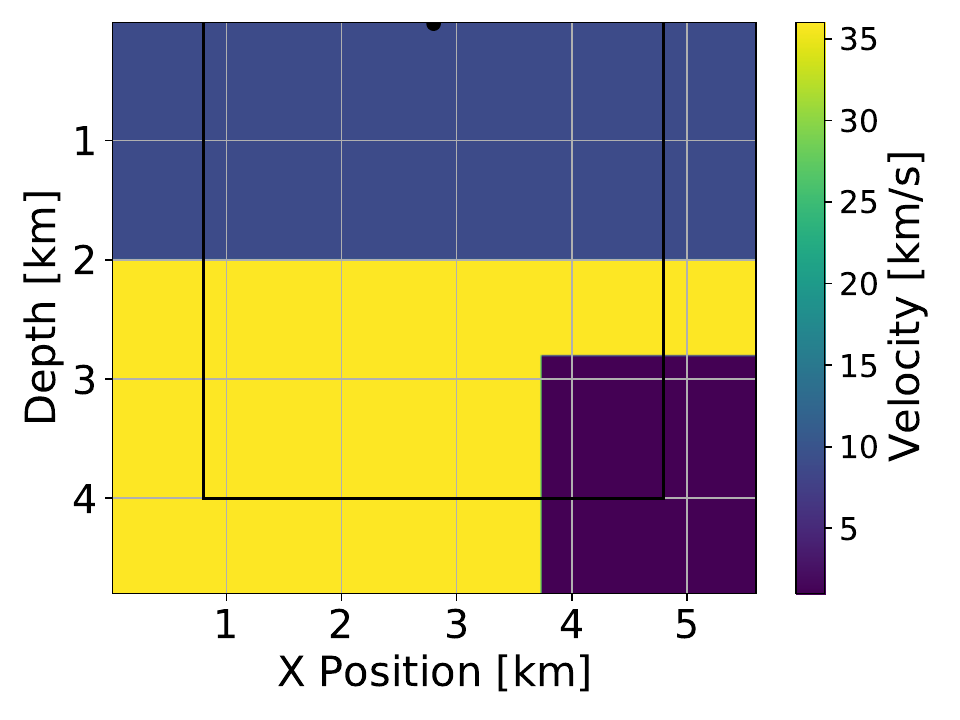}}
\hspace*{0.8cm}
\subfloat[Santos Basin.]{\includegraphics[scale=0.29]{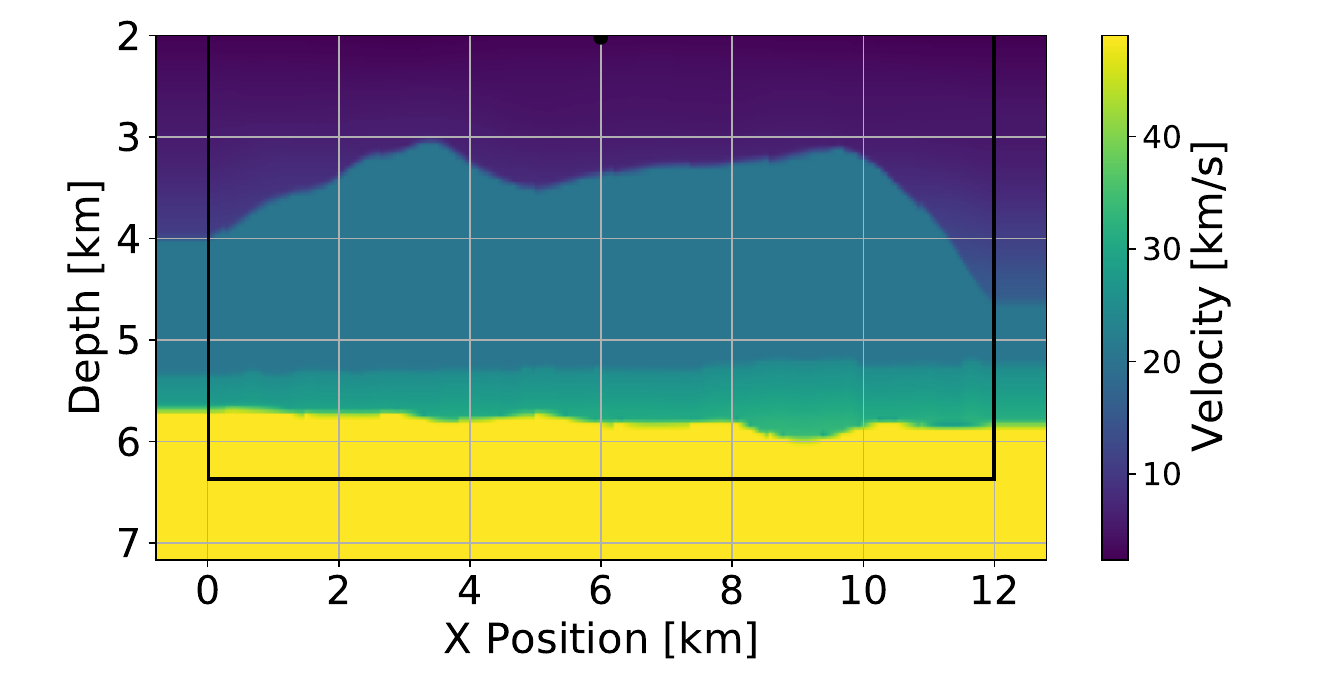}}\\
\hspace*{-0.2cm}\subfloat[Marmousi.]{\includegraphics[scale=0.29]{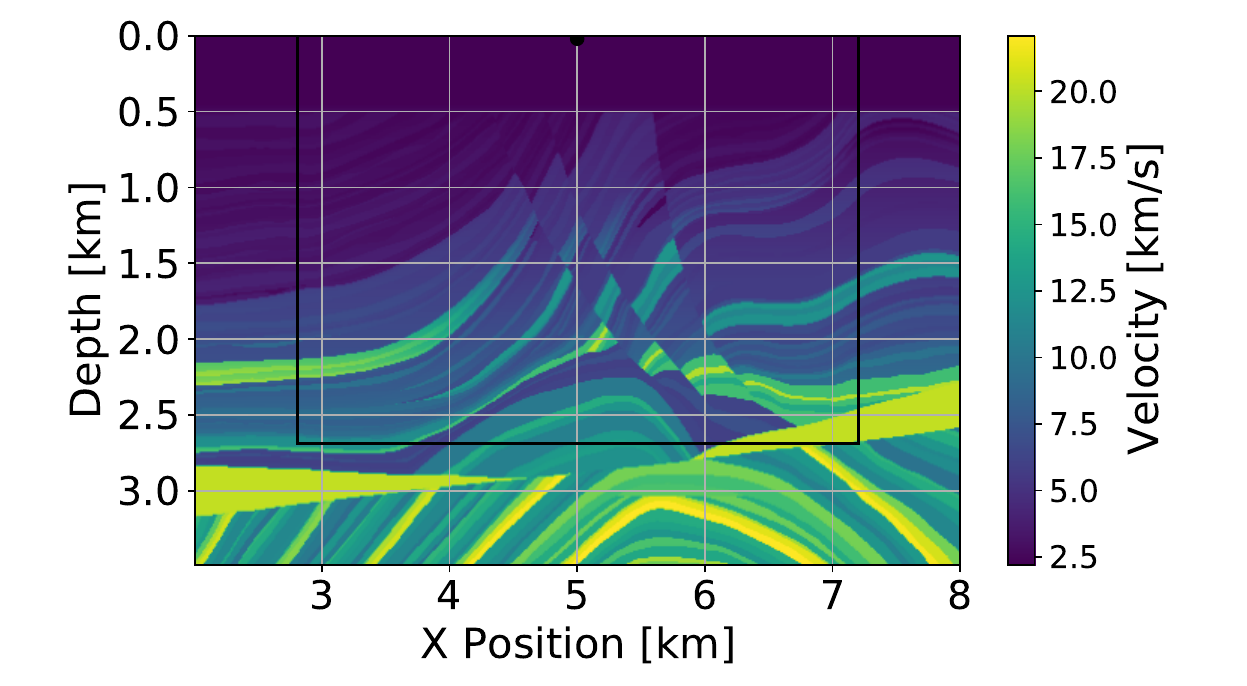}}
\hfill
\subfloat[SEG/EAGE.]{\includegraphics[scale=0.32]{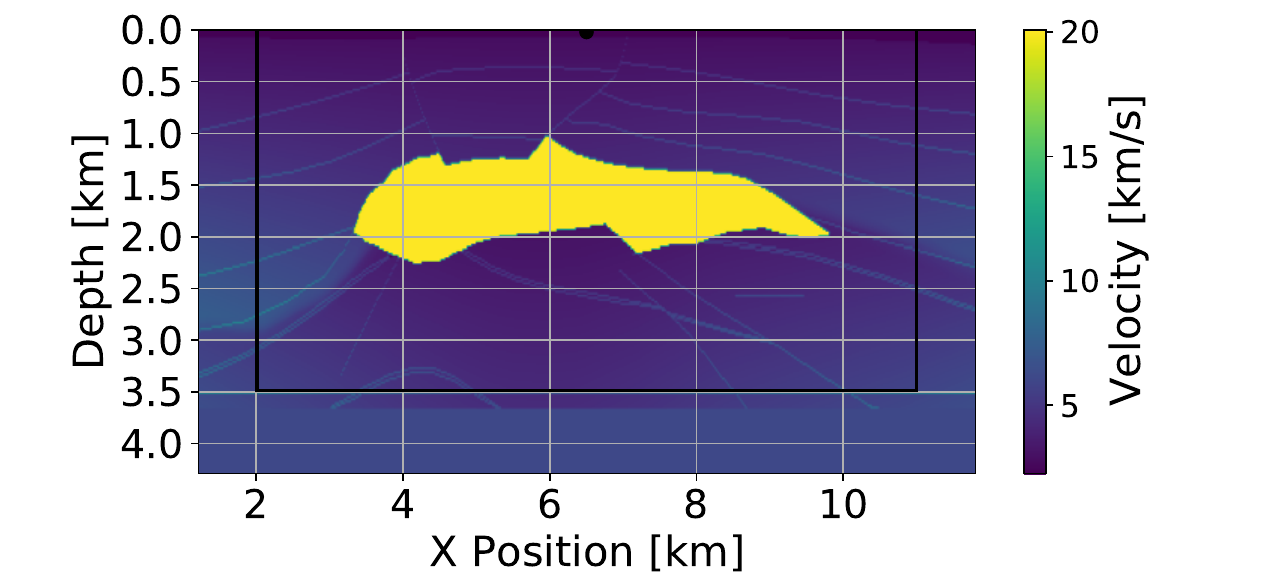}}
\caption{Velocity fields of the test cases Corner Model, Santos Basin, Marmousi, and SEG/EAGE, used to study the numerical convergence.}\label{fig_velocity_fields}
\end{figure}

We consider four numerical scenarios with different velocity fields (see Figure \ref{fig_velocity_fields}). The first is a synthetic example of a heterogeneous medium with high contrast velocities and a sharp corner. The second is a 2D slice of the velocity field of the Santos Basin\footnote{A typical velocity field of Santos Basin region, in Brazil.} oil and gas exploration region. A 2D portion of Marmousi velocity field is the third example, and the final test is the 2D SEG/EAGE synthetic model.

In all the examples, we include a source and an arrangement of receivers near the surface of the medium. The specification of this construct and other parameters of the numerical simulations are specified in Table \ref{tab_test_cases}.

\renewcommand{\arraystretch}{2}
\begin{table}[htb!]
\centering
\begin{tabular}{ |c|c|c| } 
 \hline
 \textbf{Test cases}&\textbf{Corner Model}& \textbf{Santos Basin} \\ \hline
 Domain dimensions& $\Omega=[0\text{km},4\text{km}]\times[0\text{km},4\text{km}]$ &$\Omega=[0\text{km},12\text{km}]\times[2\text{km},6.4\text{km}]$\\ \hline
 Simulation time& $T=1.1\text{s}$ &$T=1.5\text{s}$\\ \hline
 Source position& $(2\text{km},0.02\text{km})$ &$(6\text{km},2.02\text{km})$\\ \hline
 PML thickness ($\delta$) & $1.0\text{km}$&$0.8\text{km}$\\\hline
 \hline
 \textbf{Test cases}&\textbf{Marmousi}& \textbf{SEG/EAGE} \\ \hline
 Domain dimensions& $\Omega=[2\text{km},8\text{km}]\times[0\text{km},3.5\text{km}]$ &$\Omega=[2\text{km},11\text{km}]\times[0\text{km},3.5\text{km}]$\\ \hline
 Simulation time& $T=1.5\text{s}$ &$T=2\text{s}$\\ \hline
 Source position& $(5\text{km},0.02\text{km})$ &$(6.5\text{km},0.02\text{km})$\\ \hline
 PML thickness ($\delta$) & $0.8\text{km}$&$0.8\text{km}$\\\hline
\end{tabular}
\caption{Parameters of the four numerical simulations considered in this paper.}\label{tab_test_cases}
\end{table}

We save the solution at the upper boundary at each simulated time instant to construct the seismogram. We use a time span twice as long as specified in each experiment outlined in Table \ref{tab_test_cases} to allow the reflected waves to reach the surface.

\subsection{Maximum time-step}\label{sec_maximum_time_step}

We need to calculate the maximum allowable time step, denoted as $\Delta t_{\text{max}}$, for all time integration schemes. We initially consider the numerical error inherent to the spatial discretization in each numerical experiment (see Appendix \ref{sec_appendix_convergence}) since this error is independent of the time integration method. Next, we employ a tolerance level equivalent to 150\% of the spatial discretization error in each experiment. Finally. we use that tolerance to compute the value of $\Delta t_{\text{max}}$ for the numerical schemes described in Section \ref{sec_methods}. 

We consider a spatial-step size of $\Delta x=10m$ to compute the approximated solutions mentioned before. Figures \ref{fig_convergence_space} and \ref{fig_convergence_time} show the allowed $\Delta t_{\text{max}}$ by all the methods for the numerical tests Corner Model, Santos Basin, Marmousi, and SEG/EAGE.

\afterpage{
\begin{figure}[p]
\subfloat[Corner Model solution at time$\qquad$\linebreak $T=1.1s$.]{\includegraphics[scale=0.345,valign=t]{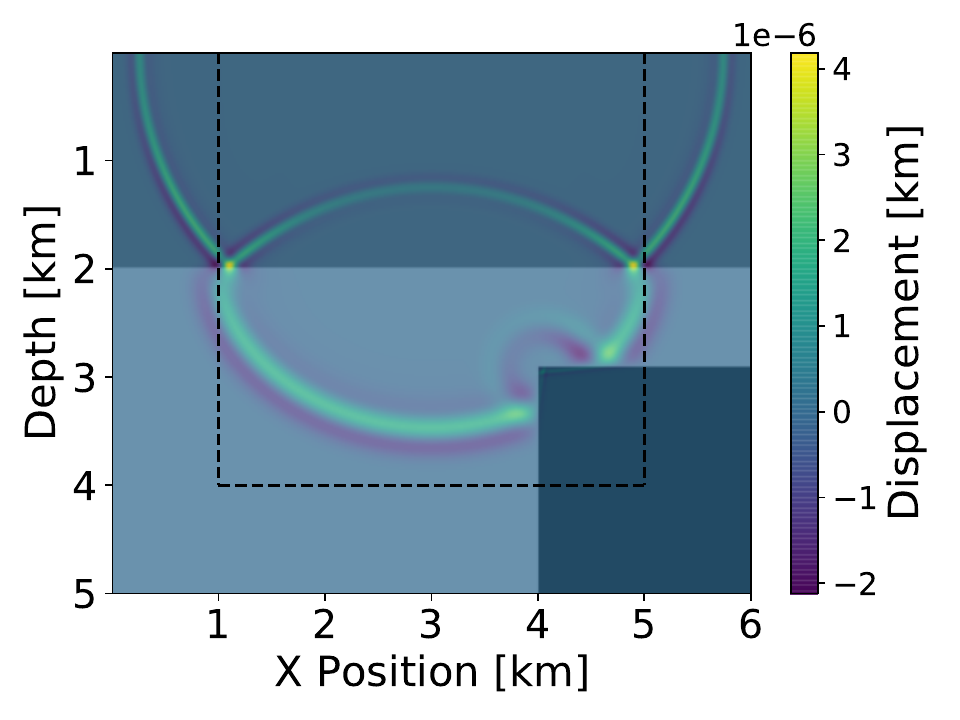}}
    \hfill
    \subfloat[Computing $\Delta t_{\text max}$ for Corner Model using an error tolerance of $3.03\cdot 10^{-7}$.]{\includegraphics[scale=0.27,valign=t]{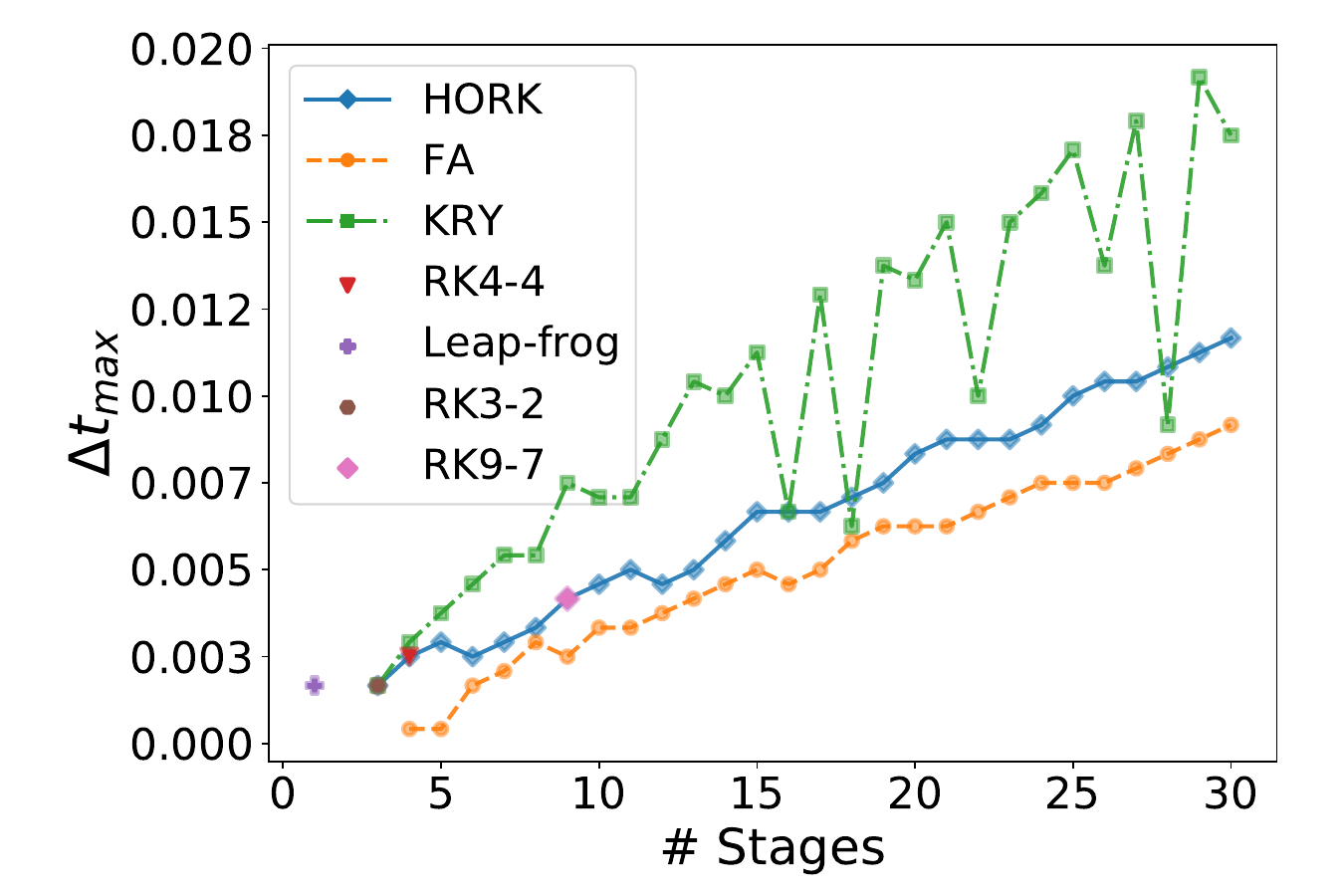}}
    \\
    \subfloat[Santos Basin solution at time$\qquad$\linebreak $T=1.5s$.]{\includegraphics[scale=0.355,valign=t]{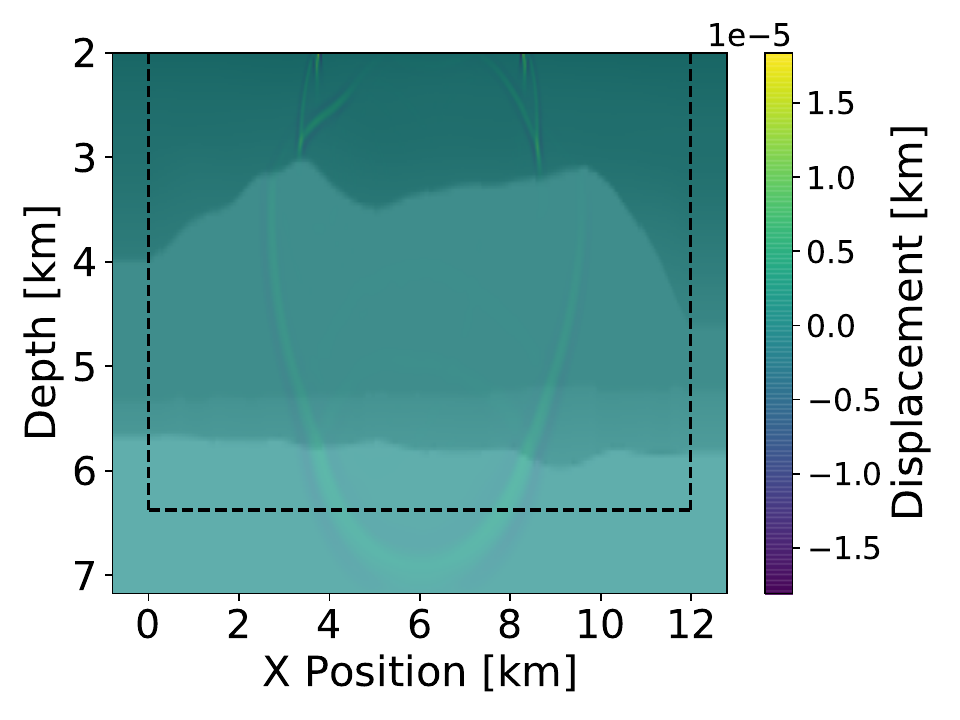}}
    \hfill
    \subfloat[Computing $\Delta t_{\text max}$ for Santos Basin using an error tolerance of $8.33\cdot 10^{-7}$.]{\includegraphics[scale=0.27,valign=t]{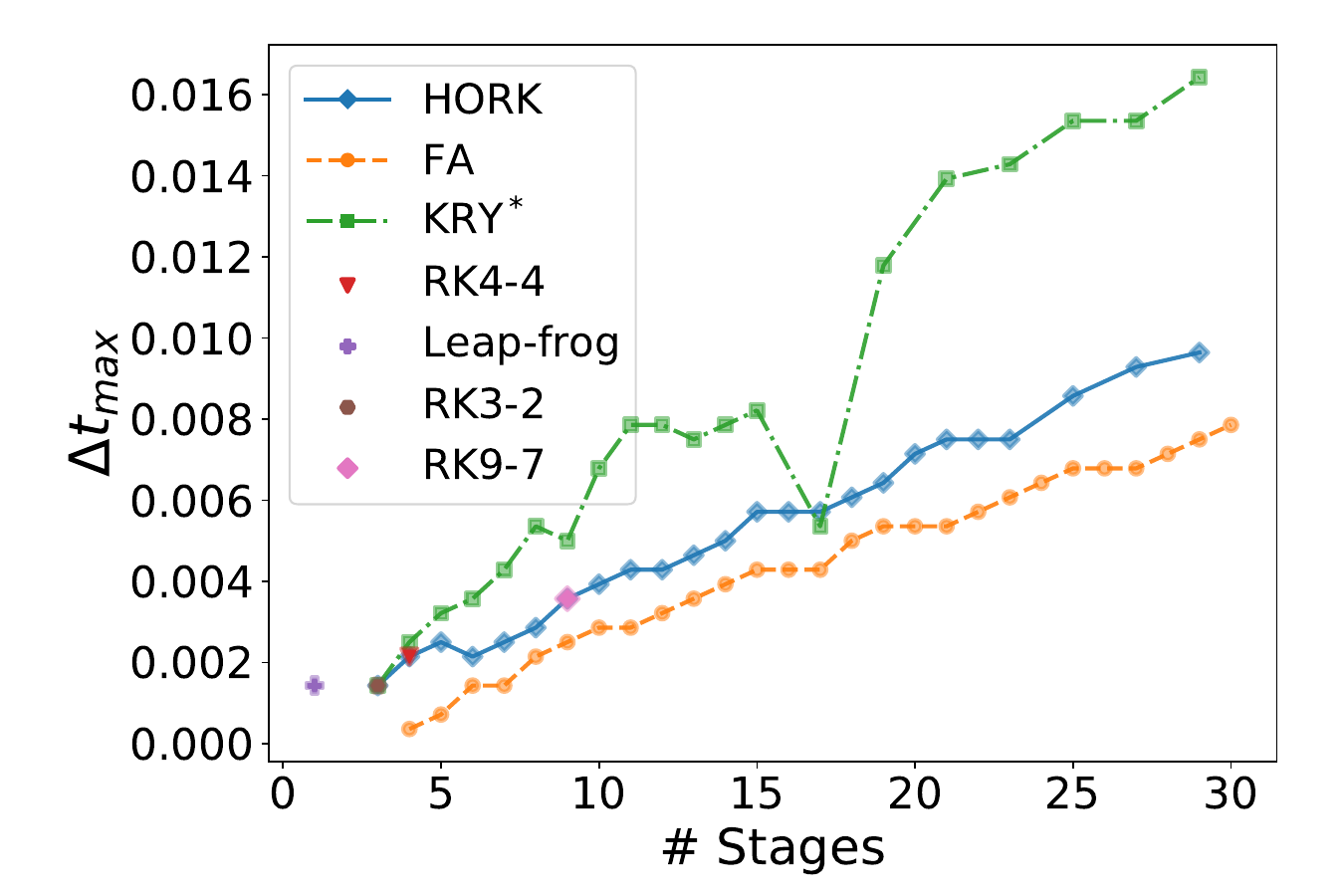}}\\
    \hspace*{-0.2cm}\subfloat[Marmousi solution at time $T=1.5s$.]{\includegraphics[scale=0.37,valign=t]{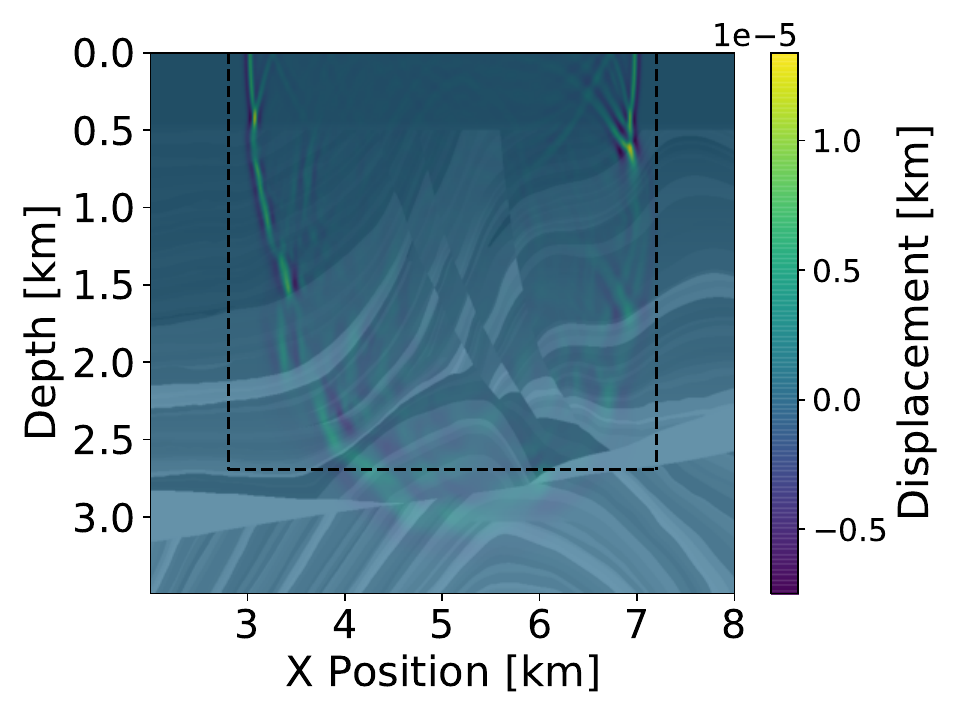}}
    \hfill
    \subfloat[Computing $\Delta t_{\text max}$ for Marmousi using an error tolerance of $9.93\cdot 10^{-7}$.]{\includegraphics[scale=0.27,valign=t]{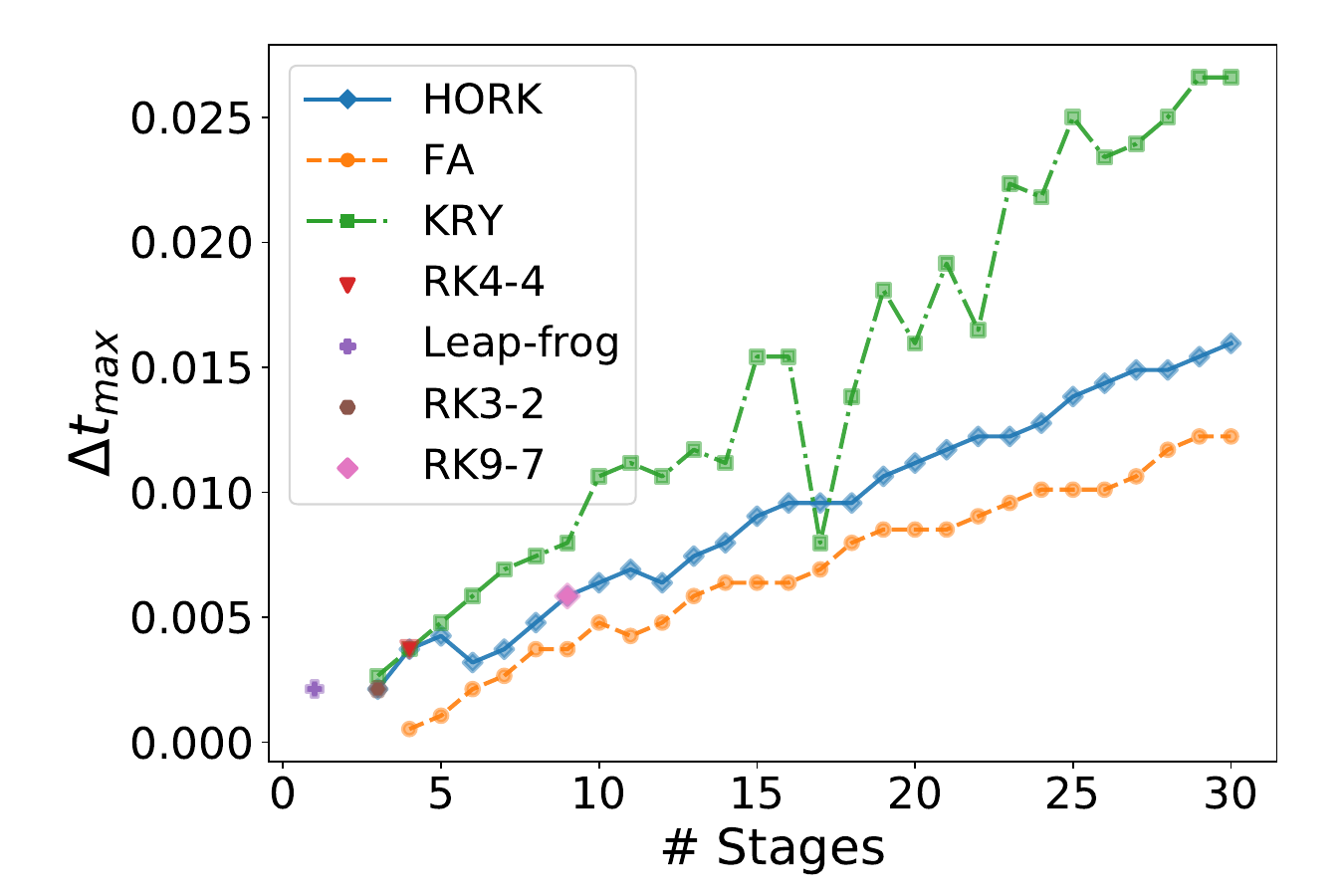}}
    \\
    \hspace*{-0.3cm}\subfloat[SEG/EAGE solution at time $T=2s$.]{\includegraphics[scale=0.37,valign=t]{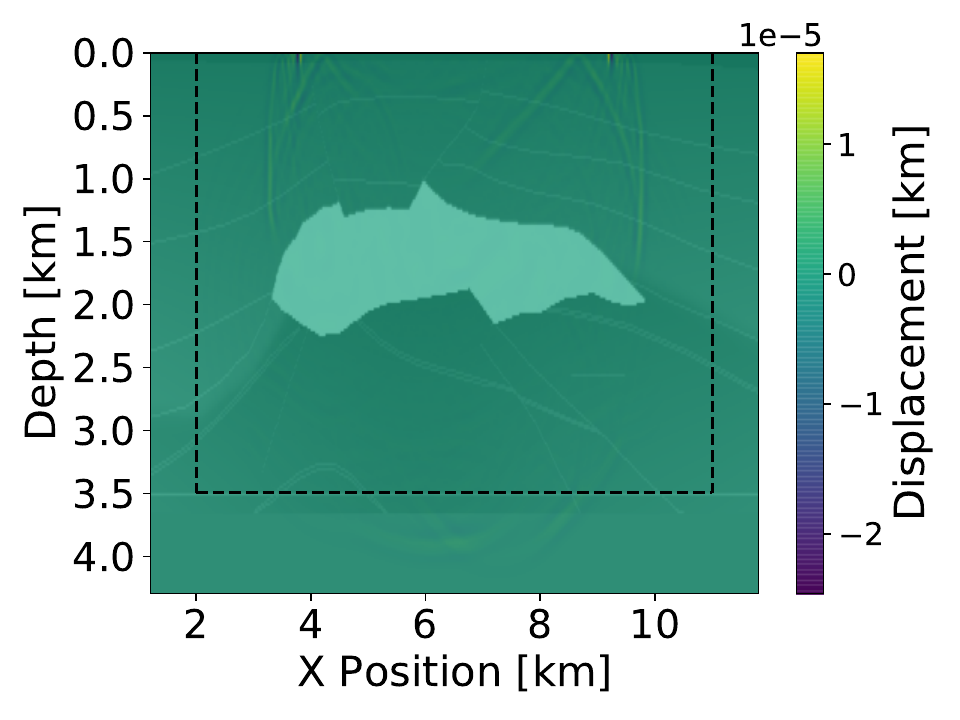}}
    \hfill
    \subfloat[Computing $\Delta t_{\text max}$ for SEG/EAGE using an error tolerance of $1.3\cdot 10^{-6}$.]{\includegraphics[scale=0.27,valign=t]{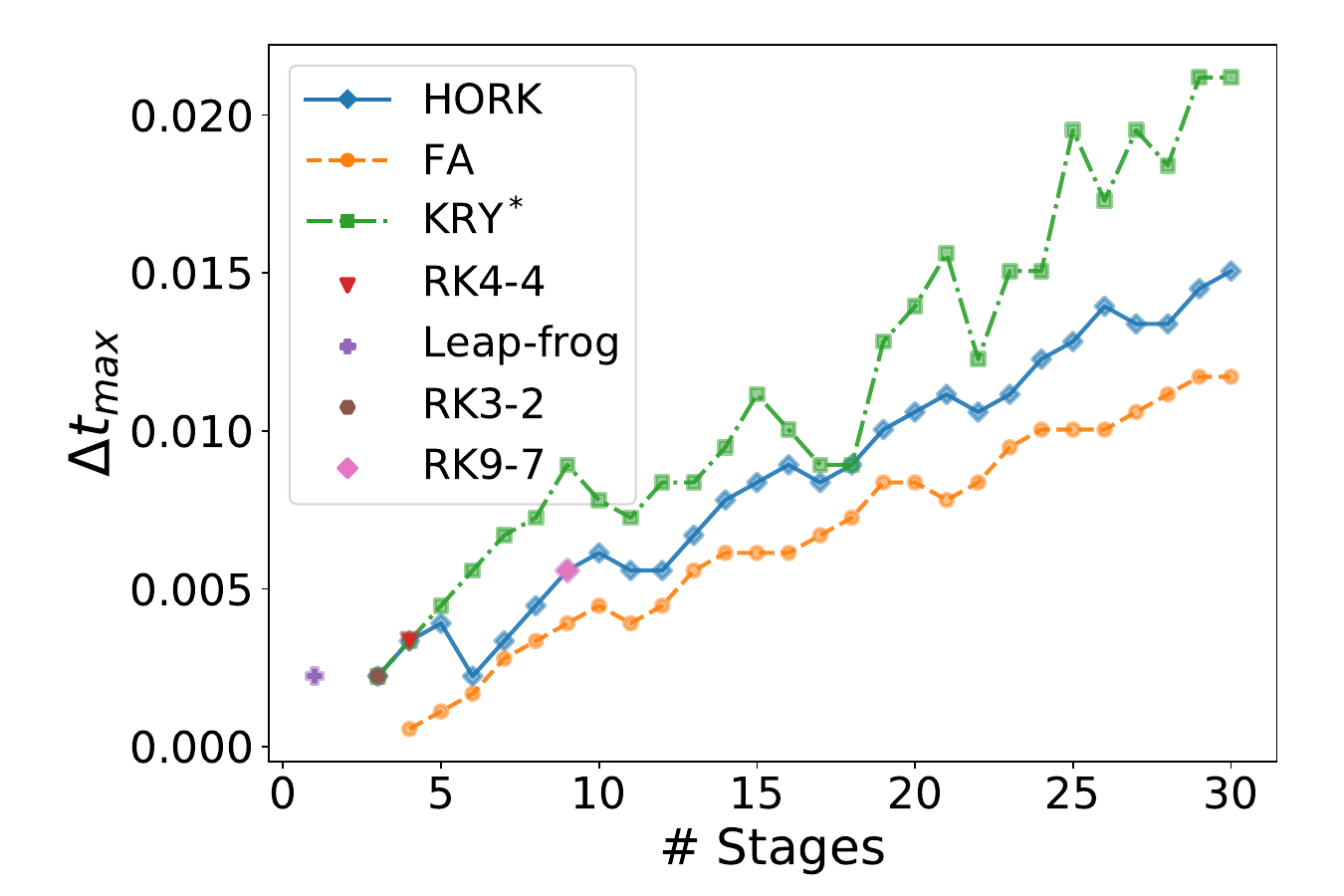}}
    \hfill
    \caption{Snapshots of the reference solution for Corner Model, Santos Basin, Marmousi, and SEG/EAGE numerical tests (left column), and the $\Delta t_{\text max}$ of each method such that the error with the reference solution is under a fixed threshold (right column). An increase in the number of stages of the method leads to a larger $\Delta t_{\text max}$.
    }\label{fig_convergence_space}
\end{figure}
\clearpage
}

Figure \ref{fig_convergence_space} presents the maximum time step, $\Delta t_\text{max}$, considering the spatial error of the solution at a time instant. Generally, when the approximation degree increases, we observe an increment in the allowed $\Delta t_\text{max}$. Moreover, the Krylov subspace approximation exhibits the largest time steps among the studied methods, followed by the other high-order methods. In contrast, low-order methods such as Leap-frog and RK3-2 require smaller time steps. 

The determination of $\Delta t_\text{max}$ based on the seismogram data is illustrated in Figure \ref{fig_convergence_time}. Similar to Figure \ref{fig_convergence_space}, an increase in the number of stages leads to a higher maximum time step. Notably, the Krylov subspace method consistently demonstrates the highest $\Delta t_\text{max}$ values, followed by other high-order methods. 

\afterpage{
\begin{figure}[p]
\subfloat[Corner Model seismogram until time $T=2.2s$.]{\includegraphics[scale=0.35,valign=t]{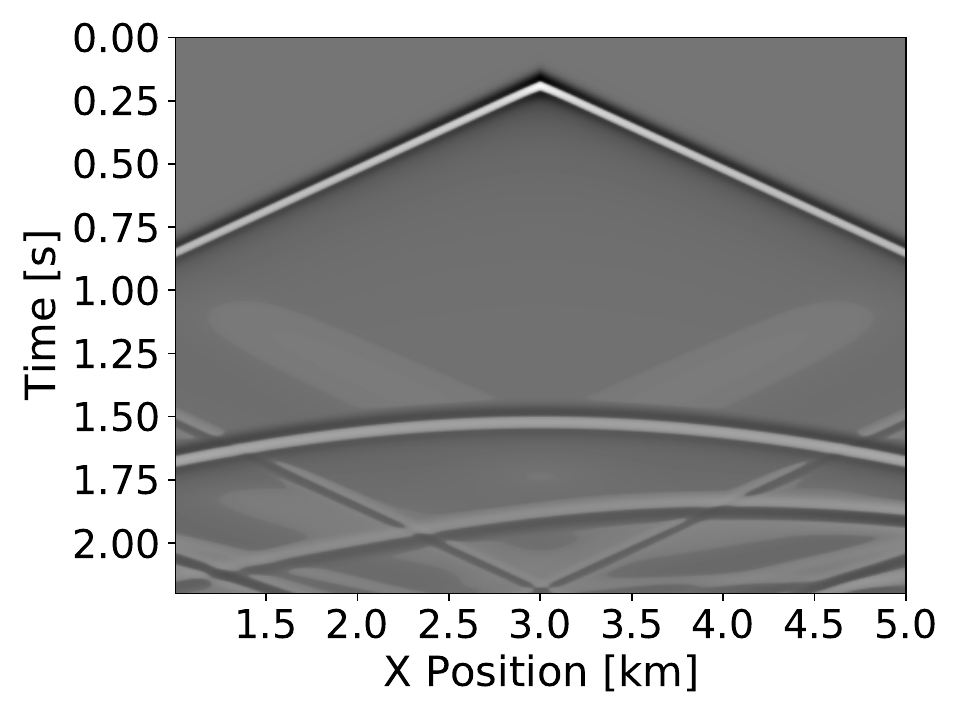}}
    \hfill
    \subfloat[Computing $\Delta t_{\text max}$ for Corner Model using an error tolerance of $2.92\cdot10^{-7}$.]{\includegraphics[scale=0.27,valign=t]{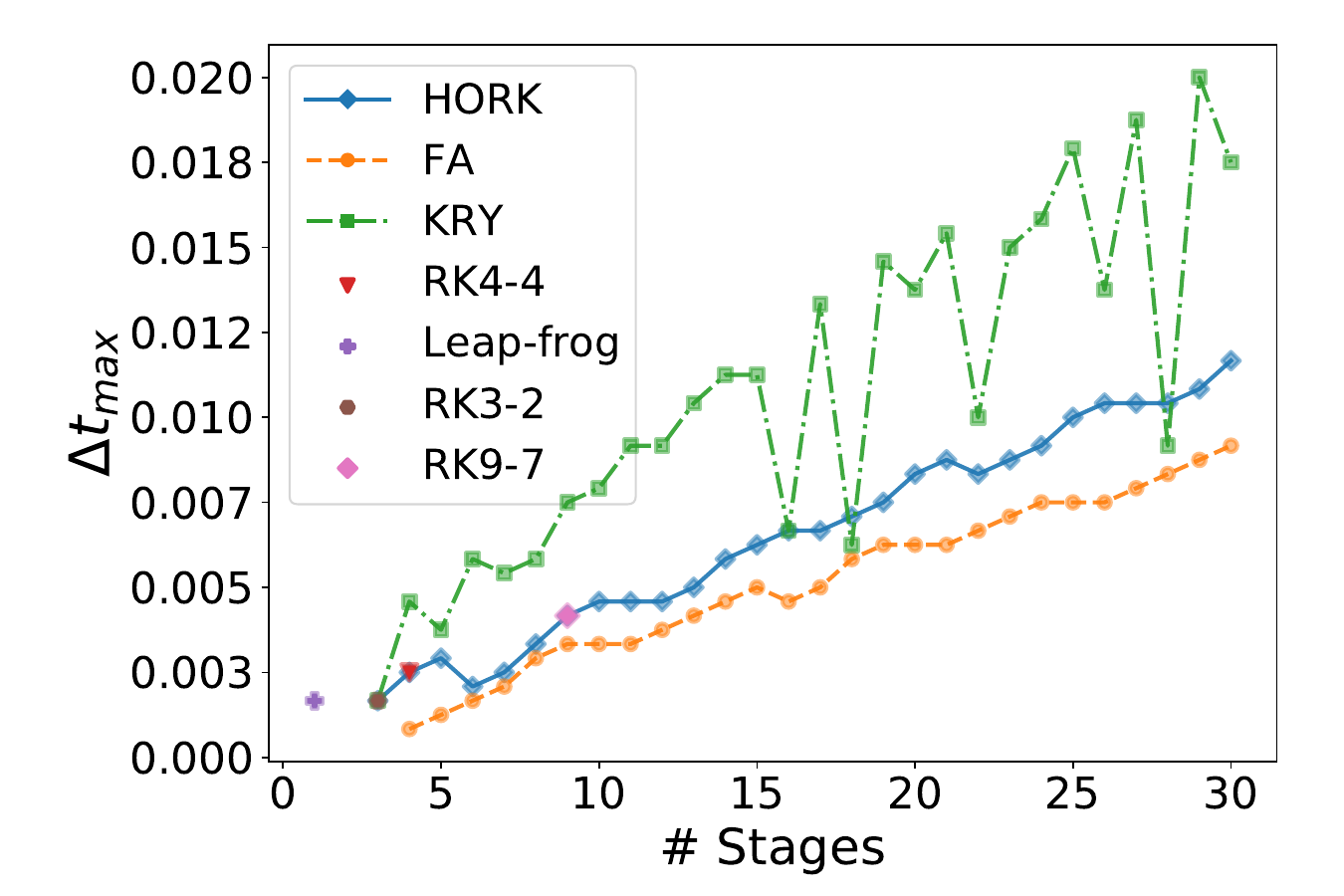}}
    \\
    \subfloat[Santos Basin seismogram until time$\qquad$\linebreak $T=3s$.]{\includegraphics[scale=0.35,valign=t]{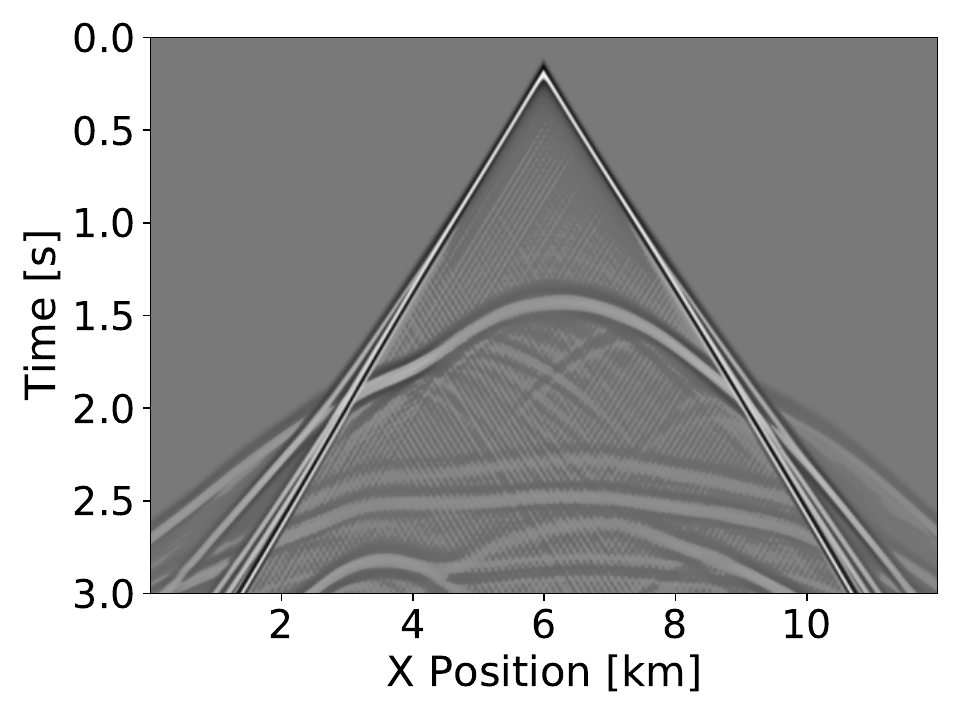}}
    \hfill
    \subfloat[Computing $\Delta t_{\text max}$ for Santos Basin using an error tolerance of $2.65\cdot 10^{-6}$.]{\includegraphics[scale=0.27,valign=t]{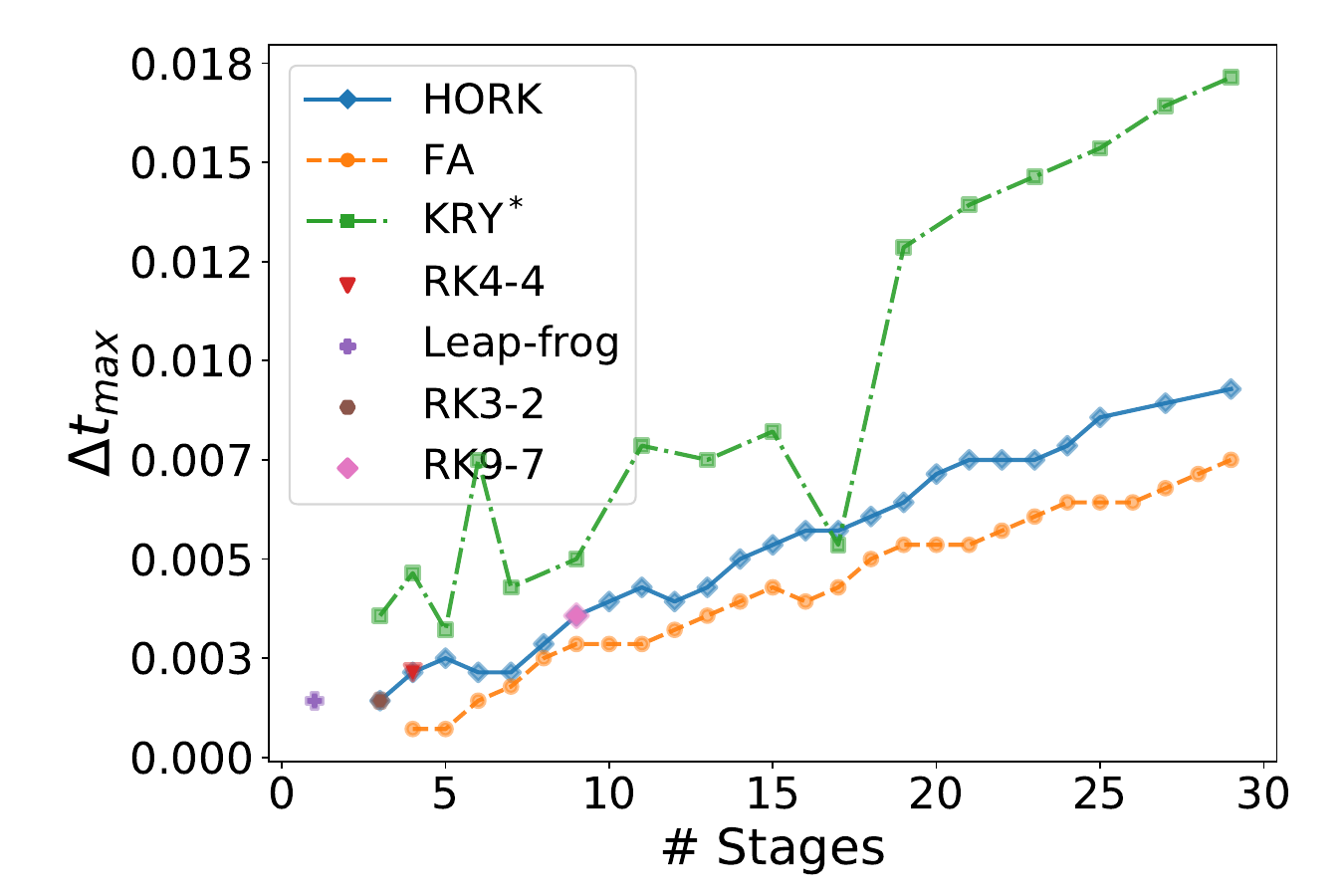}}\\
\subfloat[Marmousi seismogram until time$\qquad$\linebreak $T=3s$.]{\includegraphics[scale=0.35,valign=t]{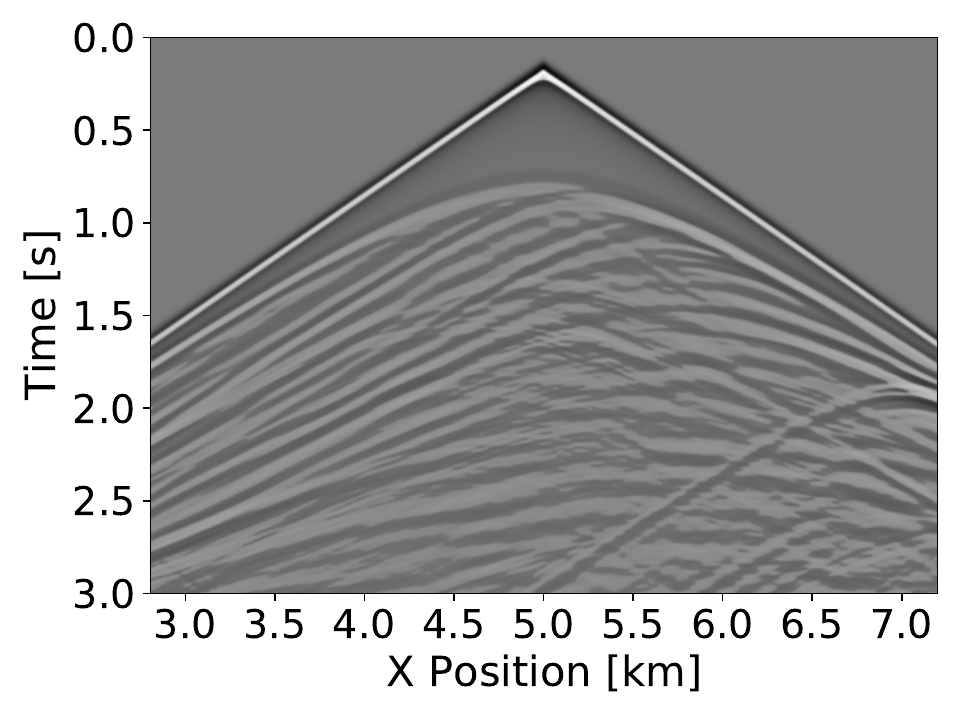}}
    \hfill
    \subfloat[Computing $\Delta t_{\text max}$ for Marmousi using an error tolerance of $1.3\cdot10^{-6}$.]{\includegraphics[scale=0.27,valign=t]{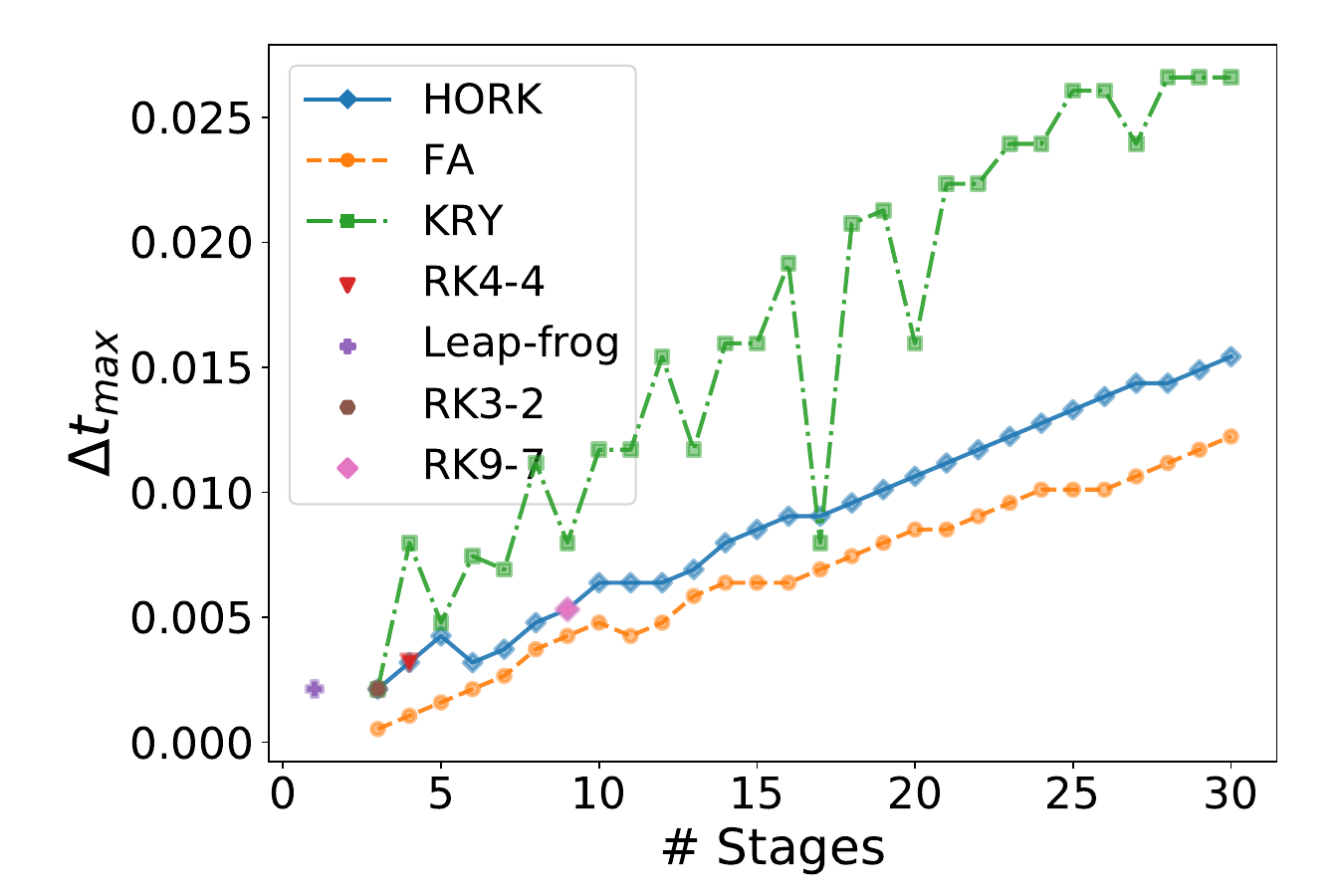}}\\
    \subfloat[SEG/EAGE  seismogram until time$\qquad$\linebreak $T=4s$.]{\includegraphics[scale=0.355,valign=t]{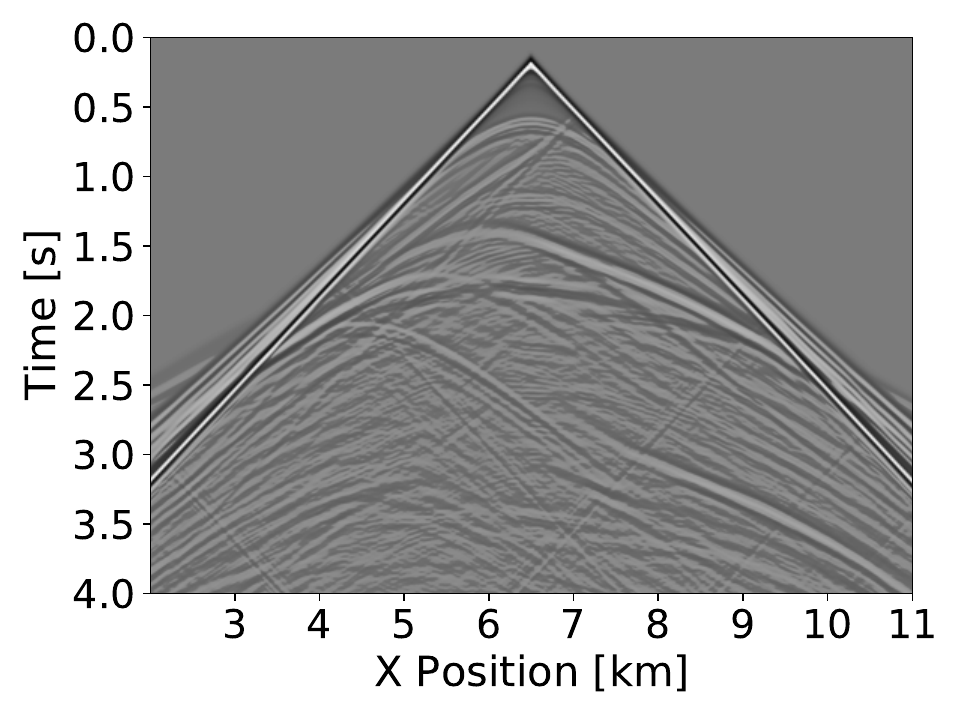}}
    \hfill
    \subfloat[Computing $\Delta t_{\text max}$ for SEG/EAGE using an error tolerance of $4.2\cdot 10^{-6}$.]{\includegraphics[scale=0.27,valign=t]{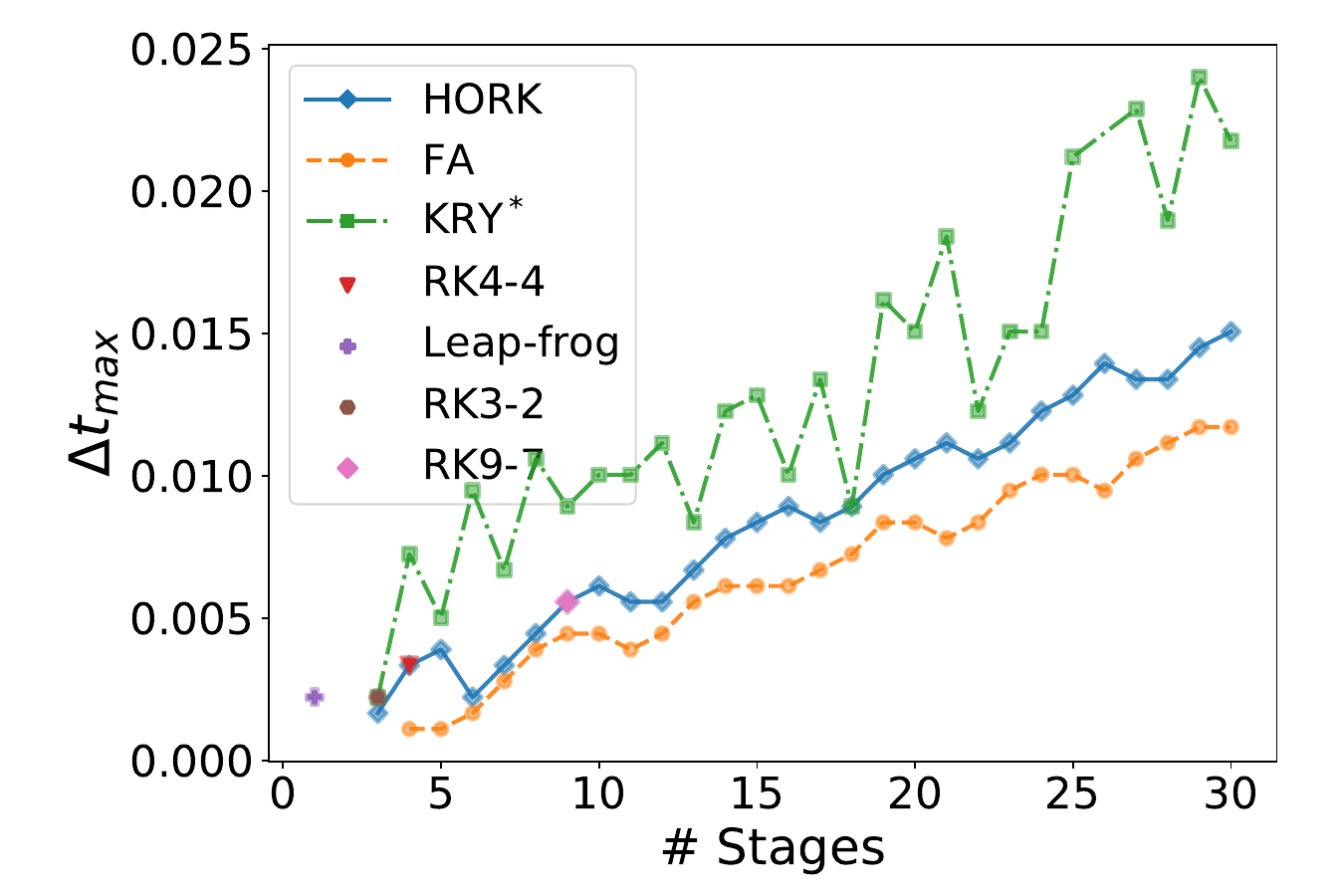}}
    \hfill
    \caption{Seismogram of the reference solution for Corner Model, Santos Basin, Marmousi, and SEG/EAGE numerical tests (left column), and the corresponding $\Delta t_{\text max}$ of each method, ensuring the error remains below a fixed threshold (right column). An increase in the number of stages of the method leads to a larger $\Delta t_{\text max}$. }\label{fig_convergence_time}
\end{figure}
\clearpage
}

Based on the insights gained from Figures \ref{fig_convergence_space} and \ref{fig_convergence_time}, we can conclude that the choice between using the error of the solution at a particular time instant or the seismogram data leads to similar values of $\Delta t_{\text{max}}$ for the methods. Therefore, for the sake of simplicity, we estimate $\Delta t_{\text{max}}$ with the error of the approximation in the physical domain at a specific time instant (as illustrated in Figure \ref{fig_convergence_space}). Next,  we estimate the computational efficiency and memory consumption of each method using the concepts of Section \ref{sec_cost_memory}.

\subsection{Computational efficiency and memory consumption}

From the previous section, we concluded that using a method with a large number of stages allows an increase in the maximum time step such that we have a solution with good accuracy. However, it is unclear if increasing the number of stages to use a larger $\Delta t$ reduces the number of operations or how it helps in utilizing the memory. To answer this question, we apply the ideas discussed in Section \ref{sec_cost_memory} and define the measure of computational efficiency
\begin{equation*}
    \text{N}^{\Delta t}_{\text{op}}=\frac{\text{\# MVOs}}{\Delta t_{\text{max}}},
\end{equation*}
and the indicator of memory consumption to store results for a backward propagation
\begin{equation*}
    \text{N}^{\Delta t}_{\text{mem}}=\frac{T}{\Delta t_{\text{max}}},
\end{equation*}
 where $T$ is the simulation time  defined by Table \ref{tab_test_cases}, for each numerical experiment.

\begin{figure}[H]
\centering
\includegraphics[trim=50 480 0 0,clip,scale=0.33]{figures/legend_b.pdf}\\[-3ex]
\subfloat[Computational cost.]{\includegraphics[scale=0.3]{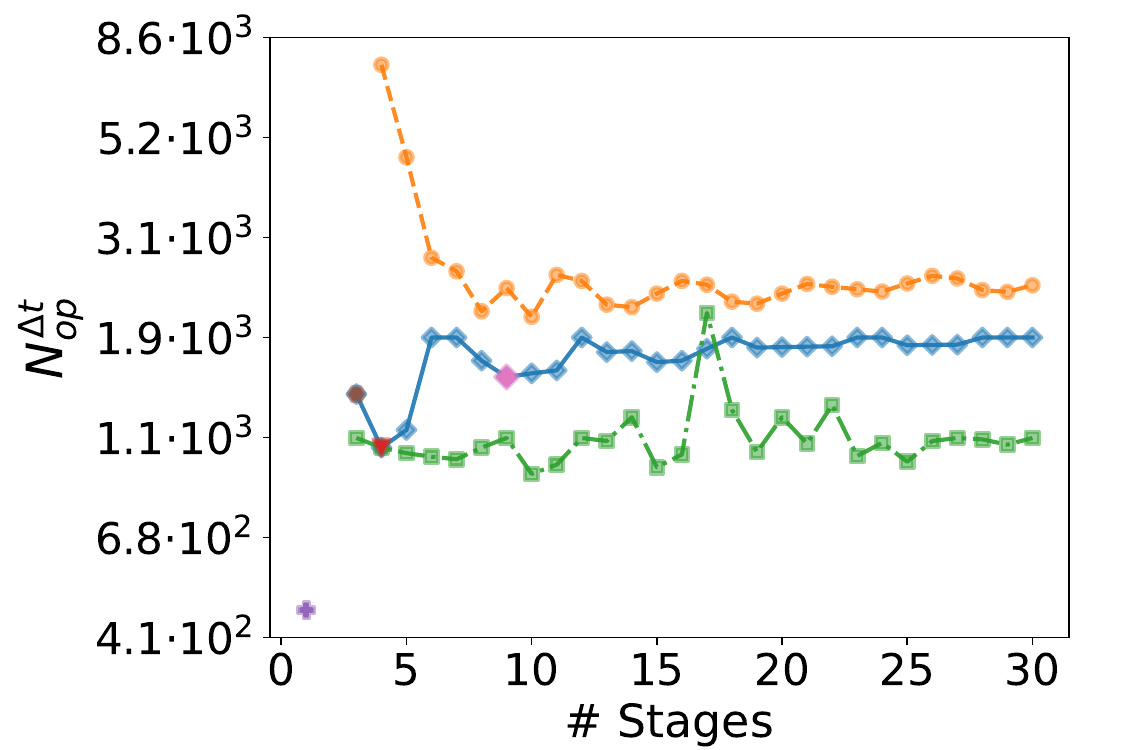}}
    \hfill
    \subfloat[Memory utilization.]{\includegraphics[scale=0.3]{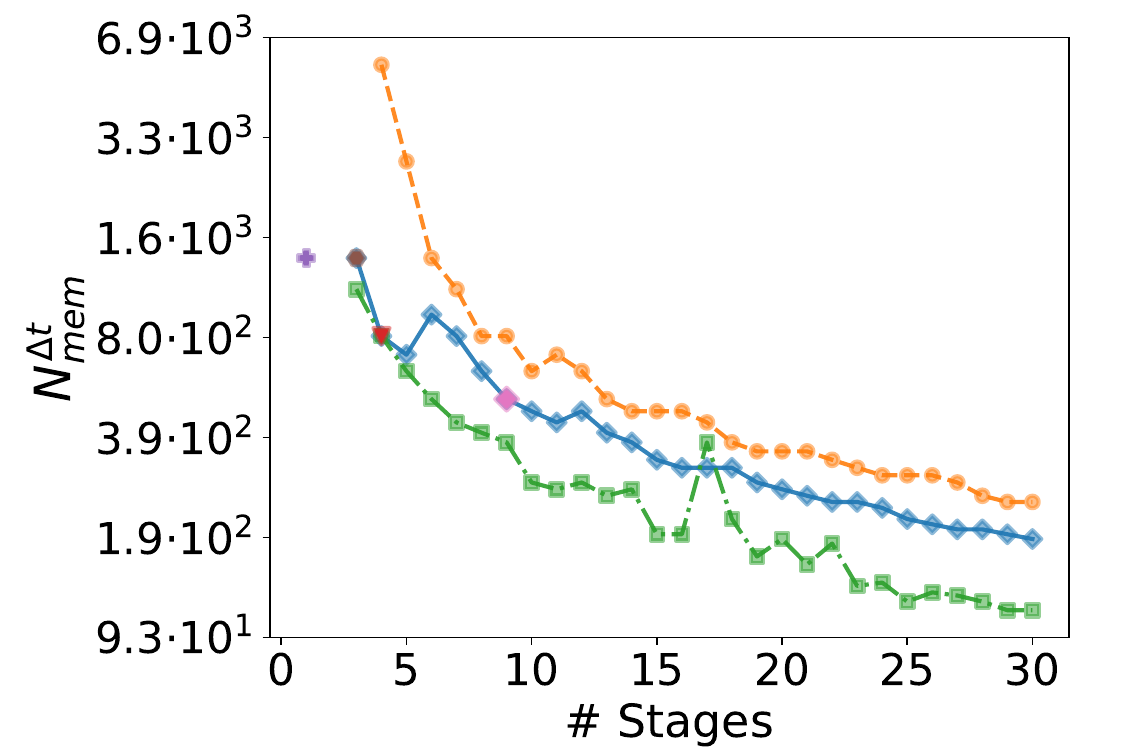}}
    \caption{Dependence of the number of MVOs and amount of stored solution vectors on the polynomial degree for the Marmousi numerical test. As the number of stages increases, the number of computations stabilizes, and memory usage decreases. * Here we neglect the computational complexity of creating the Krylov subspaces. 
    }\label{fig_efficiency_methods}
\end{figure}

Figure \ref{fig_efficiency_methods} illustrates the number of MVOs and the memory usage for all the methods when solving the Marmousi numerical example. The Leap-frog algorithm proves the most efficient among the tested methods. However, in terms of memory utilization, this method requires a substantial amount of memory. On the other hand, among the high-order methods, the Krylov subspace approximations demonstrate the best performance, even comparable to the Leap-Frog scheme.
However, we would like to point out that we are using a simplified model that doesn't consider the creation of the Krylov subspaces.
Nevertheless, we observe a significant decrease in the number of stored vectors for high-order methods in general. Additionally, we notice that further increases in the approximation degree have an attenuated effect in reducing memory utilization, which is negligible for degrees larger than 20.

We observed similar behavior in the other numerical tests, and their corresponding graphs can be found in Appendix \ref{sec_appendix_eff}.

\section{Discussion}\label{sec_discussion}

In this paper, we have implemented seven time-integration schemes, consisting of three arbitrary-order schemes based on exponential integrators and four classical low-order schemes. These algorithms have been compared through various numerical accuracy metrics, including stability, dispersion, and convergence. We have also studied the computational cost and memory requirements for each method across different approximation degrees.

The stability and dispersion analyses were conducted within a homogeneous domain by analyzing the Fourier transform of a single wave generated by a Ricker wavelet. We observed that the high-order methods were capable of using larger time steps as the polynomial degree increased. In general, we found that the Leap-frog method, although requiring smaller time steps, outperformed the high-order methods. Yet, when considering the dissipation error, the high-order methods displayed competitiveness and even surpassed all the low-order methods. 

We conducted extensive tests to evaluate convergence using four distinct velocity fields: three realistic fields and one with sharp interfaces. We assessed the approximation error both in the physical space at a specific time instant and using seismogram data. Remarkably, our results proved consistent and robust across both types of errors and all four numerical experiments.
Moreover, the Krylov method presented the largest time step size in all the tests, resulting in the least amount of solution vectors required to save for the inverse problem.
As a drawback, the Krylov method requires at each time instant to access as many vectors (with the dimensions of the solution of the wave equations) as stages of the method are used. This greatly hinders using the method to solve the direct problem. In general, high approximation degrees allowed for larger time steps, a finding that significantly impacts the number of saved vectors needed for solving the inverse problem. These results provide a different strategy to approach the memory challenges associated with the inverse problem.

This research addresses a gap in the existing literature, as most previous studies on high-order methods have predominantly focused on the spatial dimension \cite{liu2009new, weber2022stability, burman2022hybrid,wilcox2010high}. Additionally, no prior work has comprehensively examined high-order exponential integrators in the context of the wave equation applied to seismic imaging, scrutinizing the performance across a wide range of approximation degrees. Nonetheless, we acknowledge that our implementation of high-order approximations using exponential integrators is naive. Substantial enhancements are possible, particularly in terms of implementing adaptive time-stepping strategies to mitigate the hump phenomena associated with the matrix exponential \cite{moler2003nineteen}. Indeed, adaptive algorithms have been proposed, such as the KIOPS algorithm for the Krylov subspace projections, which significantly outperforms the classical Krylov method used in our study.

\backmatter

\section{Declarations}

\bmhead{Conflict of interest}

The authors declare no conflict of interest regarding this work.

\bmhead{Data avaliability}

The data and codes of the different schemes, with the exception of the Santos Basin velocity field, are available in the git-hub link: \url{https://github.com/fernanvr/Explicit-exponentials}

\bmhead{Funding}

This research was carried out in association with the ongoing R\&D project registered as ANP20714-2 STMI - Software Technologies for Modelling and Inversion, with applications in seismic imaging (USP/Shell Brasil/ANP). It was funded in part by the Coordenação de Aperfeiçoamento de Pessoal de Nível Superior - Brasil (CAPES) - Finance Code 001, and in part by Conselho Nacional de Desenvolvimento Científico e Tecnológico (CNPq) - Brasil. Fundação de Amparo à Pesquisa do Estado de São Paulo (FAPESP) grant 2021/06176-0 is also acknowledged. It has also partially received funding from the Federal Ministry of Education and Research and the European High-Performance Computing Joint Undertaking (JU) under grant agreement No 955701, Time-X. The JU receives support from the European Union’s Horizon 2020 research and innovation programme and Belgium, France, Germany, and Switzerland.

\begin{appendices}
\section{Approximations at the free-surface}\label{sec_appendix_approximation_free_surface}

We present the finite difference approximations of 8th order for the required derivatives of the functions at the points near the free surface. To simplify the notation, we define $u_i=u(x,-i\Delta x)$, and $w_i=w_y\left(x,-(i+\frac{1}{2})\Delta x\right)$. Since we are considering a uniform grid, we have that $\Delta y=\Delta x$, and so, only $\Delta x$ will be used.

\begin{align*}
    \frac{\partial^2 u}{\partial y^2}(x,0)&\approx\left(-\frac{3144919}{352800}u_0+16u_1-14u_2+\frac{112}{9}u_3-\frac{35}{4}u_4+\frac{112}{25}u_5-\frac{14}{9}u_6\right.\\
    &\left.+\frac{16}{49}u_7-\frac{1}{32}u_8\right)\frac{1}{\Delta x^2},\\
    \frac{\partial^2 u}{\partial y^2}(x,-\Delta x)&\approx\left(\frac{271343}{156800}u_0-\frac{1991}{630}u_1+\frac{57}{40}u_2+\frac{13}{60}u_3-\frac{109}{288}u_4+\frac{6}{25}u_5-\frac{11}{120}u_6\right.\\
    &\left.+\frac{179}{8820}u_7-\frac{9}{4480}u_8\right)\frac{1}{\Delta x^2},\\
    \frac{\partial^2 u}{\partial y^2}(x,-2\Delta x)&\approx\left(-\frac{18519}{78400}u_0+\frac{58}{35}u_1-\frac{251}{90}u_2+\frac{22}{15}u_3-\frac{1}{16}u_4-\frac{14}{225}u_5+\frac{1}{30}u_6\right.\\
    &\left.-\frac{2}{245}u_7+\frac{17}{20160}u_8\right)\frac{1}{\Delta x^2},\\
    \frac{\partial^2 u}{\partial y^2}(x,-3\Delta x)&\approx\left(\frac{74801}{1411200}u_0-\frac{37}{140}u_1+\frac{67}{40}u_2-\frac{263}{90}u_3+\frac{53}{32}u_4-\frac{23}{100}u_5+\frac{13}{360}u_6\right.\\
    &\left.-\frac{1}{245}u_7+\frac{1}{4480}u_8\right)\frac{1}{\Delta x^2},\\
    \frac{\partial u}{\partial y}\left(x,-\frac{1}{2}\Delta x\right)&\approx\left(\frac{5034629}{3763200}u_0-\frac{23533}{15360}u_1+\frac{4259}{15360}u_2-\frac{1103}{9216}u_3+\frac{151}{3072}u_4\right.\\
    &\left.-\frac{1171}{76800}u_5+\frac{139}{46080}u_6-\frac{211}{752640}u_7\right)\frac{1}{\Delta x},\\    
    \frac{\partial u}{\partial y}\left(x,-\frac{3}{2}\Delta x\right)&\approx\left(-\frac{363509}{3763200}u_0+\frac{6297}{5120}u_1-\frac{6147}{5120}u_2+\frac{211}{3072}u_3+\frac{3}{1024}u_4\right.\\
    &\left.-\frac{153}{25600}u_5+\frac{29}{15360}u_6-\frac{57}{250880}u_7\right)\frac{1}{\Delta x},\\ 
    \frac{\partial u}{\partial y}\left(x,-\frac{5}{2}\Delta x\right)&\approx\left(\frac{4631}{250880}u_0-\frac{305}{3072}u_1+\frac{1245}{1024}u_2-\frac{3725}{3072}u_3+\frac{275}{3072}u_4-\frac{69}{5120}u_5\right.\\
    &\left.+\frac{5}{3072}u_6-\frac{5}{50176}u_7\right)\frac{1}{\Delta x}.
\end{align*}

\section{Homogeneous medium}\label{sec_appendix_dips_diss}

\begin{figure}[H]
\subfloat[Convergence error.]{\includegraphics[scale=0.33]{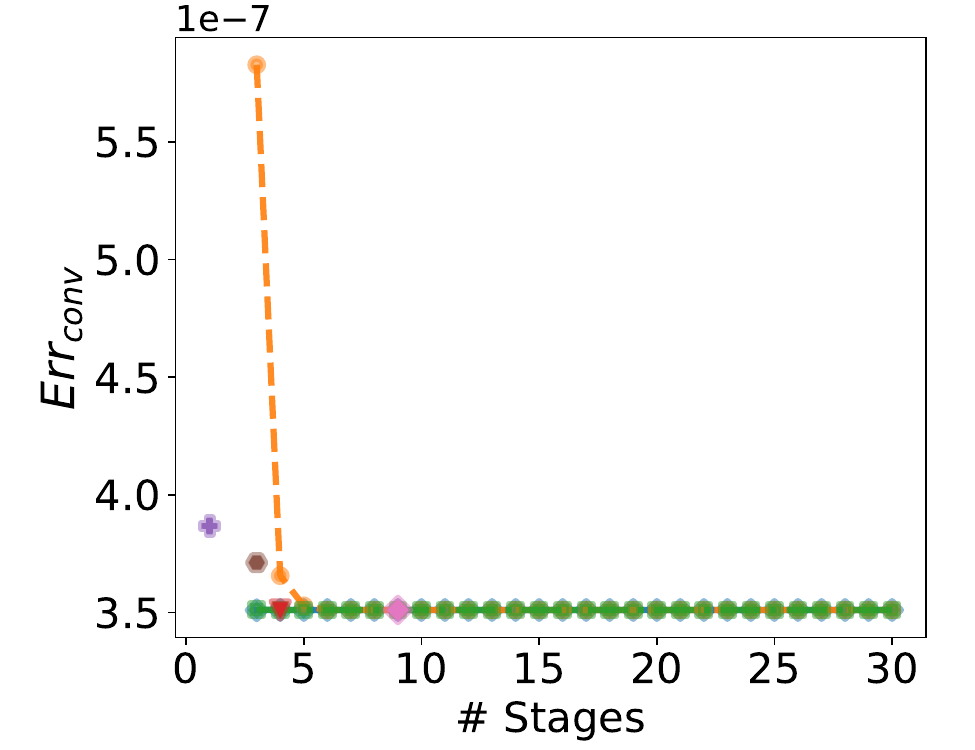}}\hfill
\subfloat[Dispersion error.]{\includegraphics[scale=0.33]{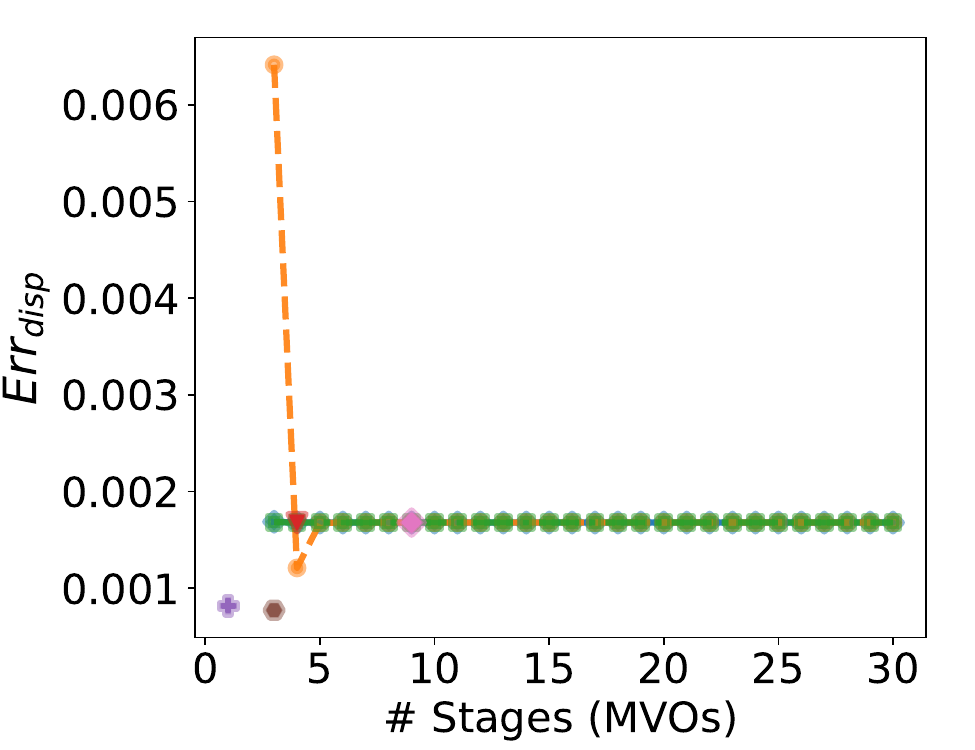}}\\
\subfloat[Dissipation error.]{\includegraphics[scale=0.33]{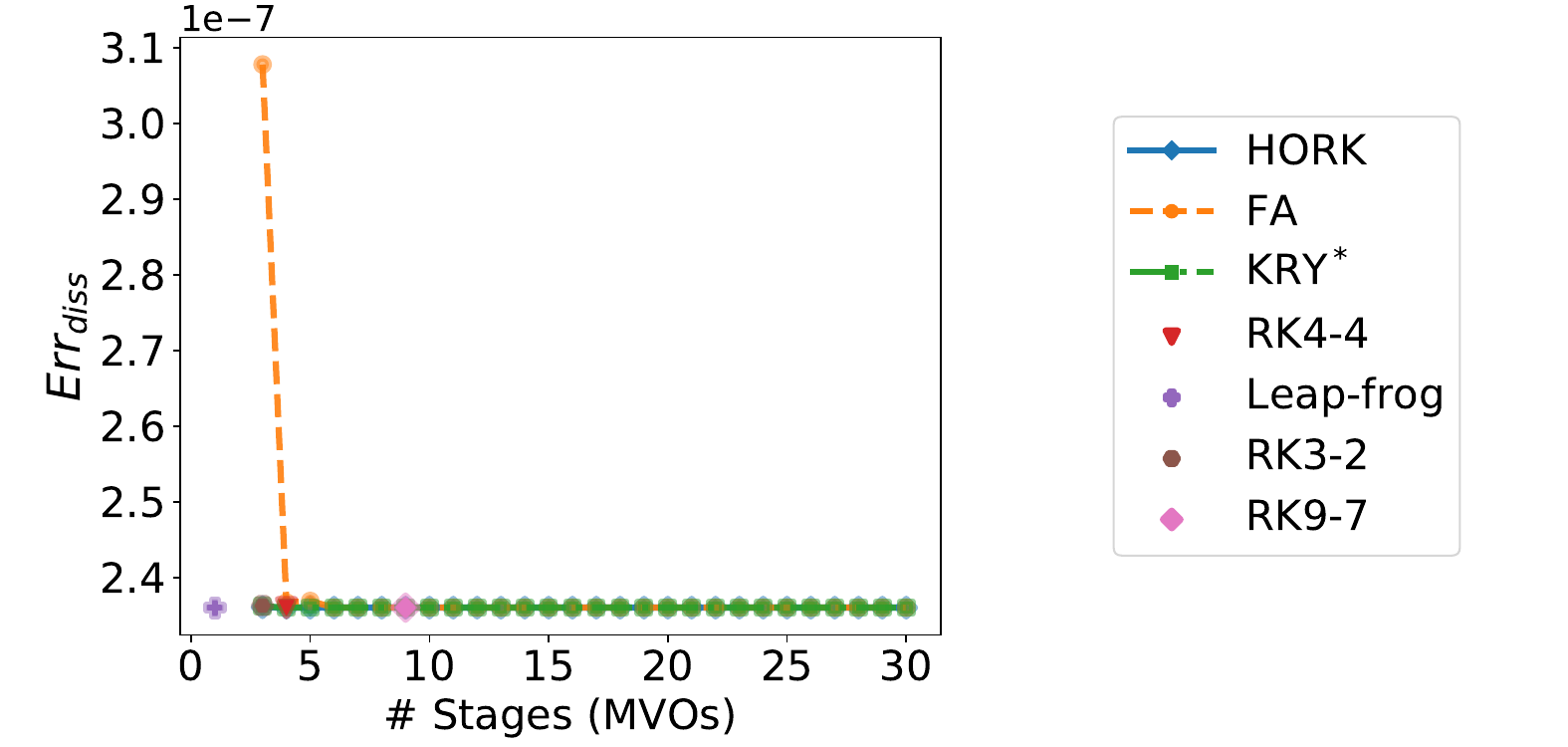}}
\caption{Convergence, dispersion, and dissipation errors using the time-step $\Delta t=\frac{\Delta x}{8c}$ for different numerical methods, with a peak frequency of $f_M=15$Hz. The approximation order does not matter, since there is an error associated to the spatial discretization. }\label{fig_disp_diss_min_delta}
\end{figure}

This section complements the results in Section \ref{sec_disp_diss}. First, we show the convergence, dispersion, and dissipation errors associated with the eighth-order spatial discretization scheme using $\Delta x=10m$ (Figure \ref{fig_disp_diss_min_delta}). Additionally, we present how varying the peak frequencies as $f_M={10,\;15,\;20,\;25}$, impact the maximum allowable time-step $\Delta t_{\text{max}}$ and the number of matrix-vector operations (MVOs) for different schemes and approximation degrees.

\subsection{Dispersion results}\label{sec_appendix_dips_diss_disp}

\begin{figure}[H]
\centering
\includegraphics[trim=50 480 0 0,clip,scale=0.33]{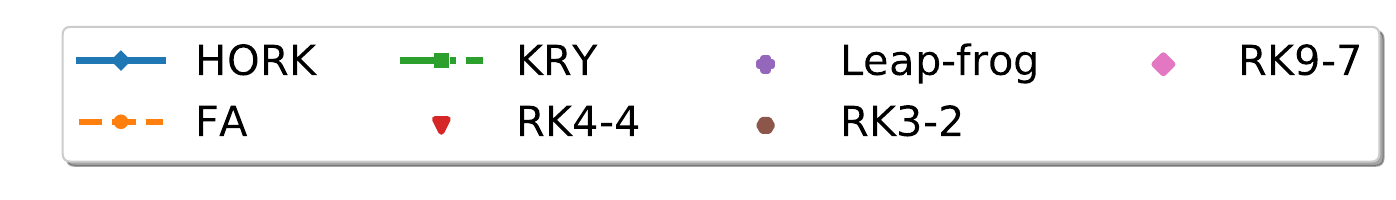}\\[-3ex]
\subfloat[Peak frequency $f_M=10$Hz.]{\includegraphics[scale=0.33]{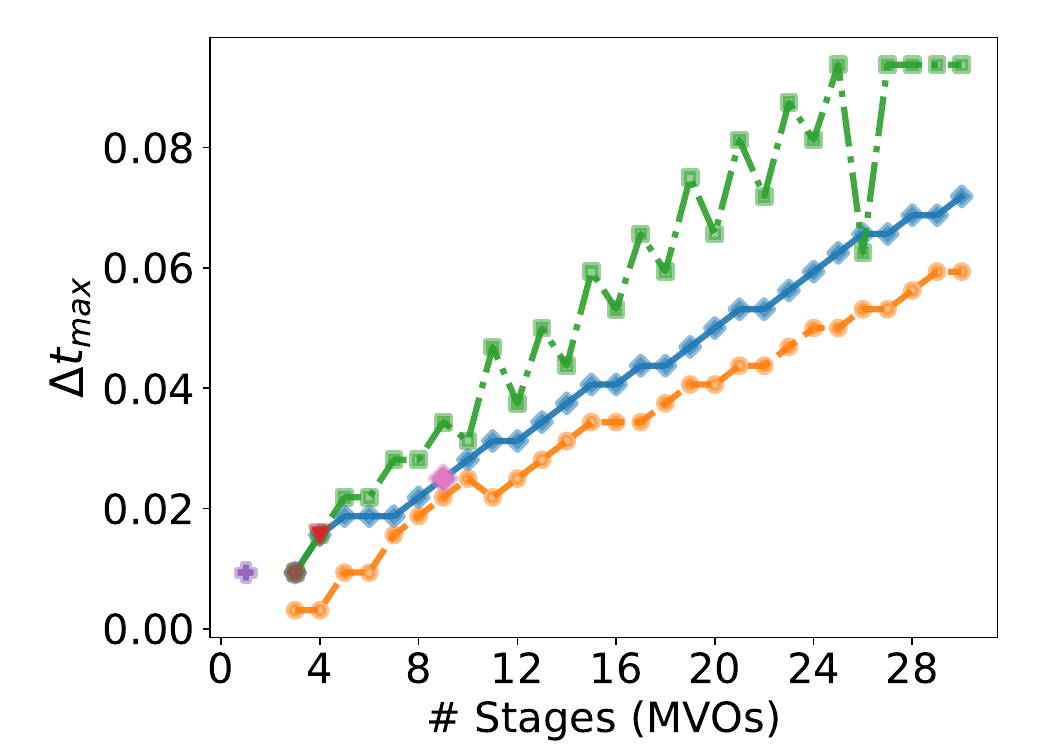}}
    \hfill
    \subfloat[Peak frequency $f_M=15$Hz.]{\includegraphics[scale=0.33]{figures/methods_ord_8_dx_0.01_f0_15_dispersion_max_dt.pdf}}
    \hfill
    \\
    \subfloat[Peak frequency $f_M=20$Hz.]{\includegraphics[scale=0.33]{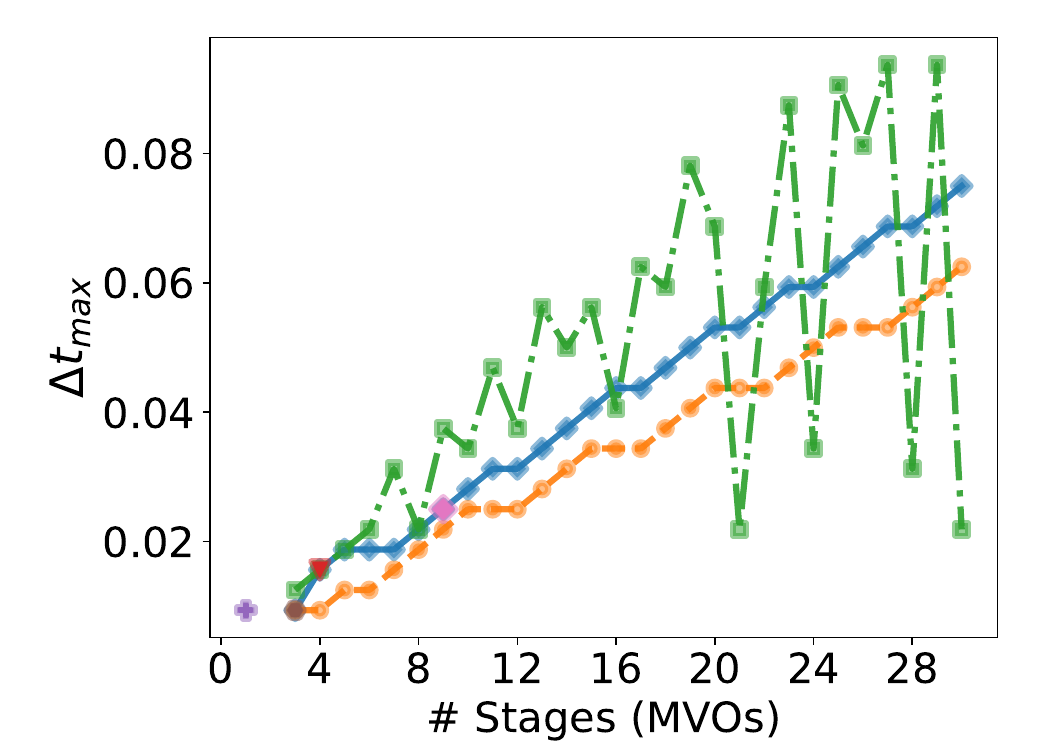}}
    \hfill
    \subfloat[Peak frequency $f_M=25$Hz.]{\includegraphics[scale=0.33]{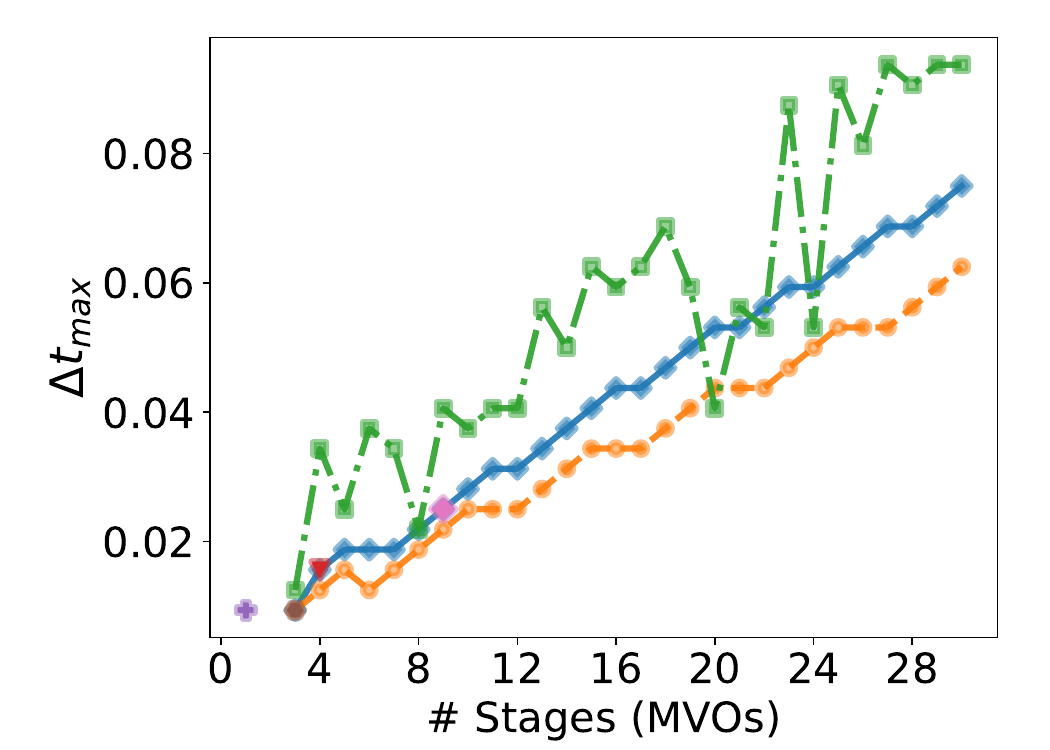}}
    \caption{ Maximum time step ($\Delta t_{\text{max}}$) while controlling the time dispersion error of each method to be below 50\% of the spatial dispersion error concerning different peak frequencies of the Ricker wavelet. A grater number of stages generally allows larger time steps.
    }\label{fig_delta_t_disp}
\end{figure}

From Fig. \ref{fig_delta_t_disp}, we perceive that the general behavior is maintained independent of the peak frequencies. With the difference that when the peak frequency increases, the results for the Krylov method are more oscillatory, and the high-degree approximations using Faber polynomials suffer from more round-off errors.

\begin{figure}[H]
\centering
\includegraphics[trim=50 480 0 0,clip,scale=0.33]{figures/legend_b.pdf}\\[-3ex]
\subfloat[Peak frequency $f_M=10$Hz.]{\includegraphics[scale=0.33]{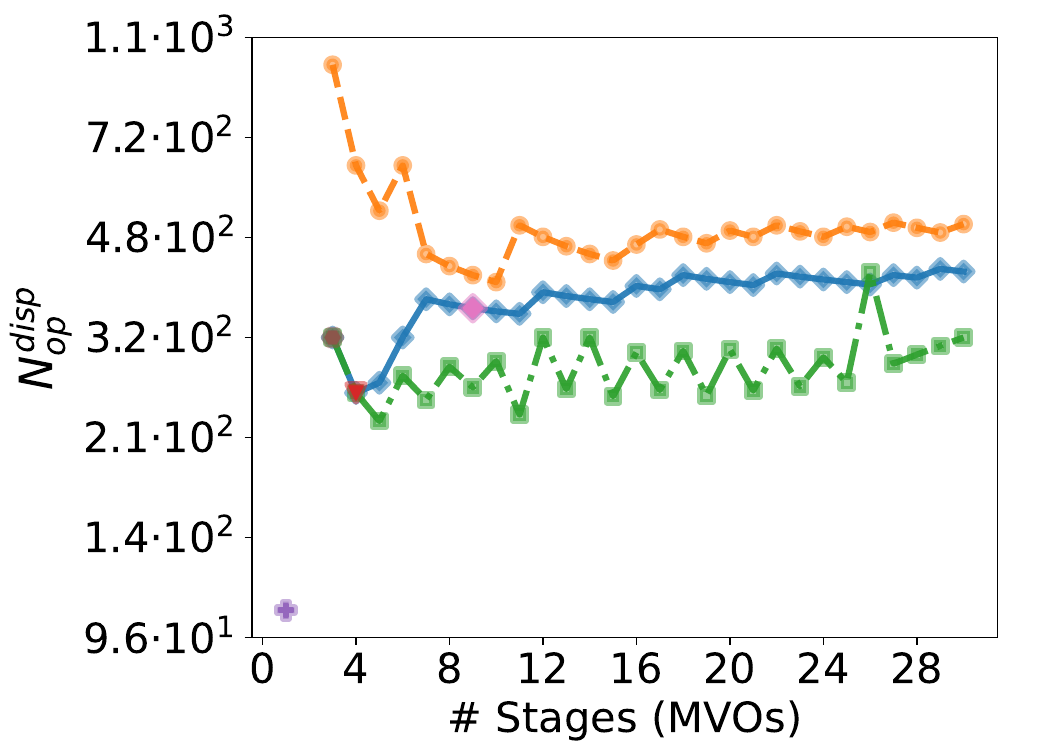}}
    \hfill
    \subfloat[Peak frequency $f_M=15$Hz.]{\includegraphics[scale=0.33]{figures/methods_ord_8_dx_0.01_f0_15_dispersion_eff.pdf}}
    \hfill
    \\
    \subfloat[Peak frequency $f_M=20$Hz.]{\includegraphics[scale=0.33]{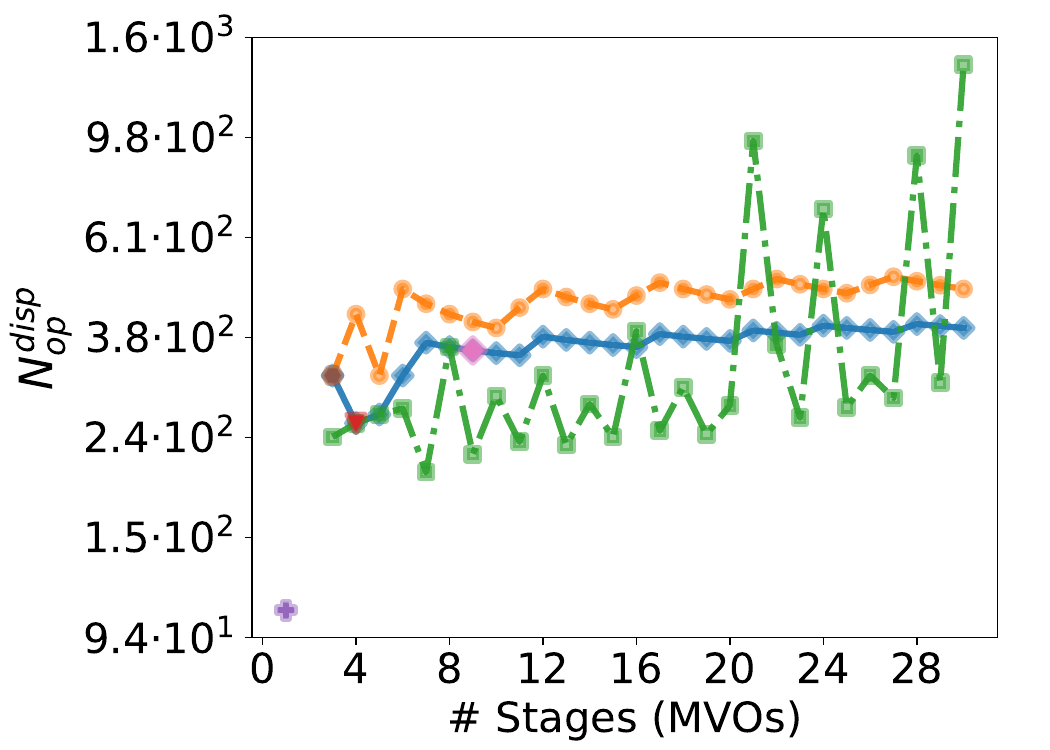}}
    \hfill
    \subfloat[Peak frequency $f_M=25$Hz.]{\includegraphics[scale=0.33]{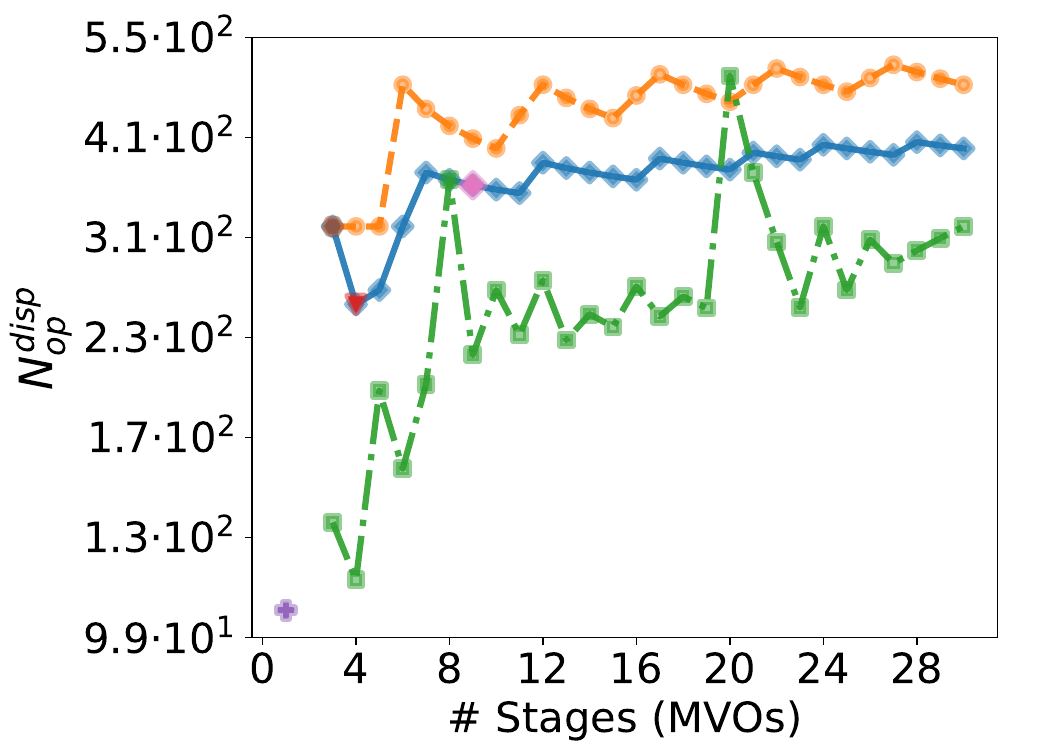}}
    \caption{Dependence of the number of matrix-vector operations and the maximum time-step required to compute the solution on the polynomial degree, considering different peak frequencies. While increasing the number of stages generally leads to a slight increment in computations. * Here we neglect the computational complexity of creating the Krylov subspaces.
    }\label{fig_dispersion_eff}
\end{figure}

In Figure \ref{fig_dispersion_eff}, we still observe that the Leap-frog algorithm requires the least amount of MVOs. The FA and HORK methods share a similar number of computations independent of the peak frequency.

\subsection{Dissipation results}\label{sec_appendix_dips_diss_diss}

\begin{figure}[H]
\centering
\includegraphics[trim=50 480 0 0,clip,scale=0.33]{figures/legend.pdf}\\[-3ex]
\subfloat[Peak frequency $f_M=10$Hz.]{\includegraphics[scale=0.33]{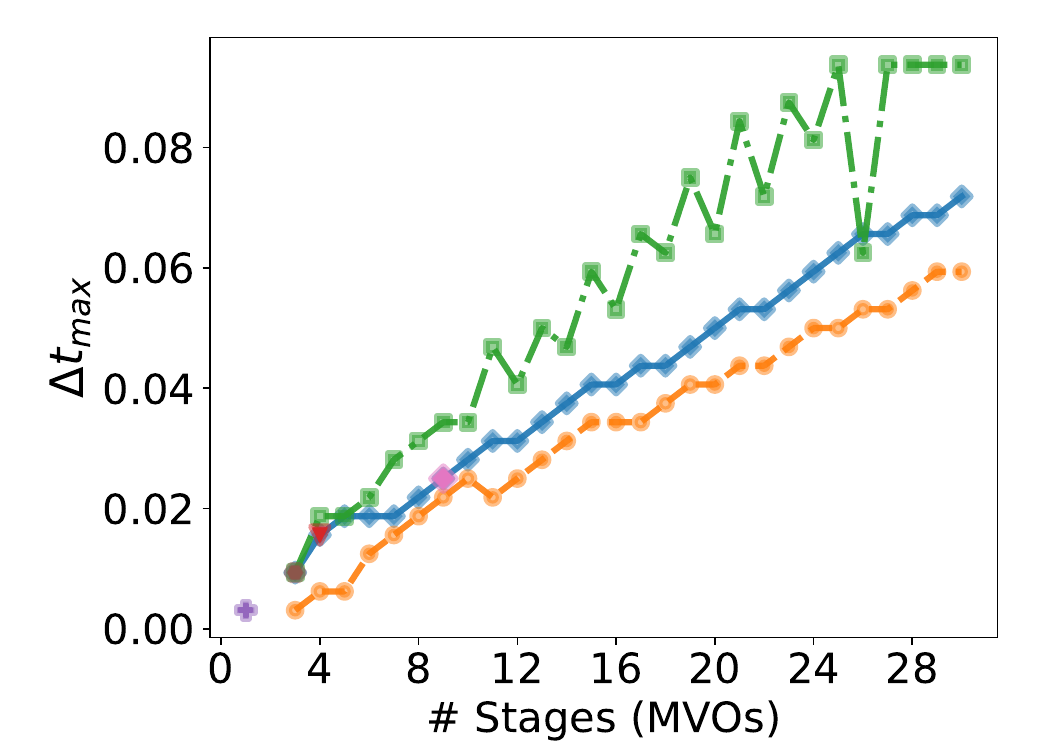}}
    \hfill
    \subfloat[Peak frequency $f_M=15$Hz.]{\includegraphics[scale=0.33]{figures/methods_ord_8_dx_0.01_f0_15_dissipation1_max_dt.pdf}}
    \hfill
    \\
    \subfloat[Peak frequency $f_M=20$Hz.]{\includegraphics[scale=0.33]{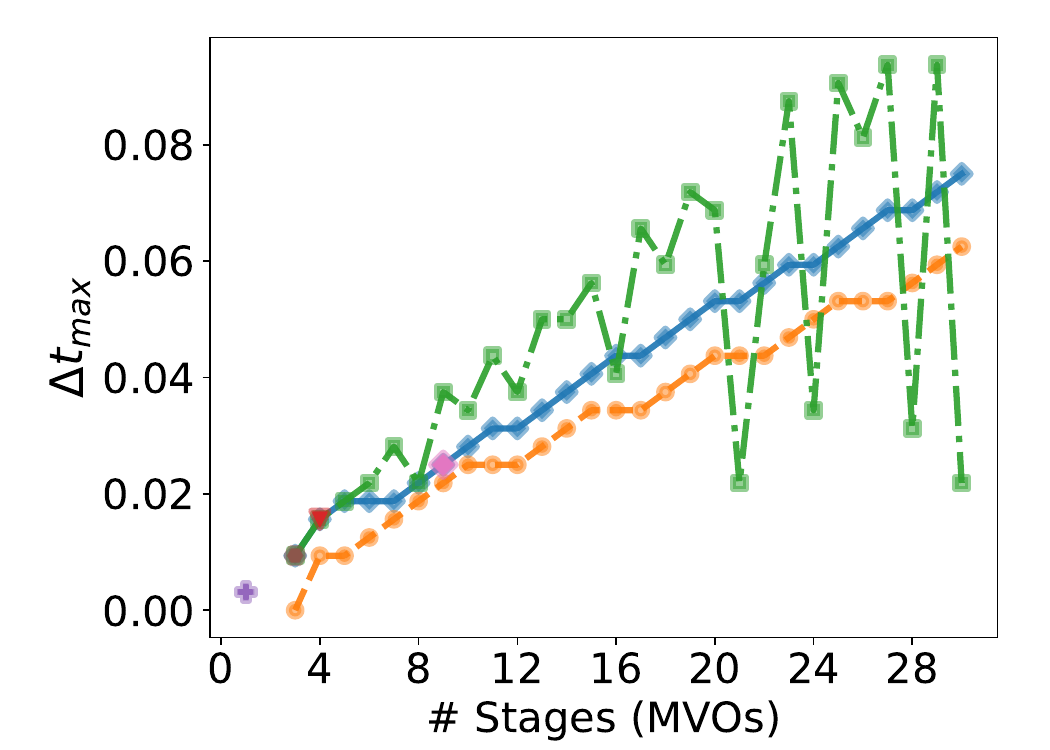}}
    \hfill
    \subfloat[Peak frequency $f_M=25$Hz.]{\includegraphics[scale=0.33]{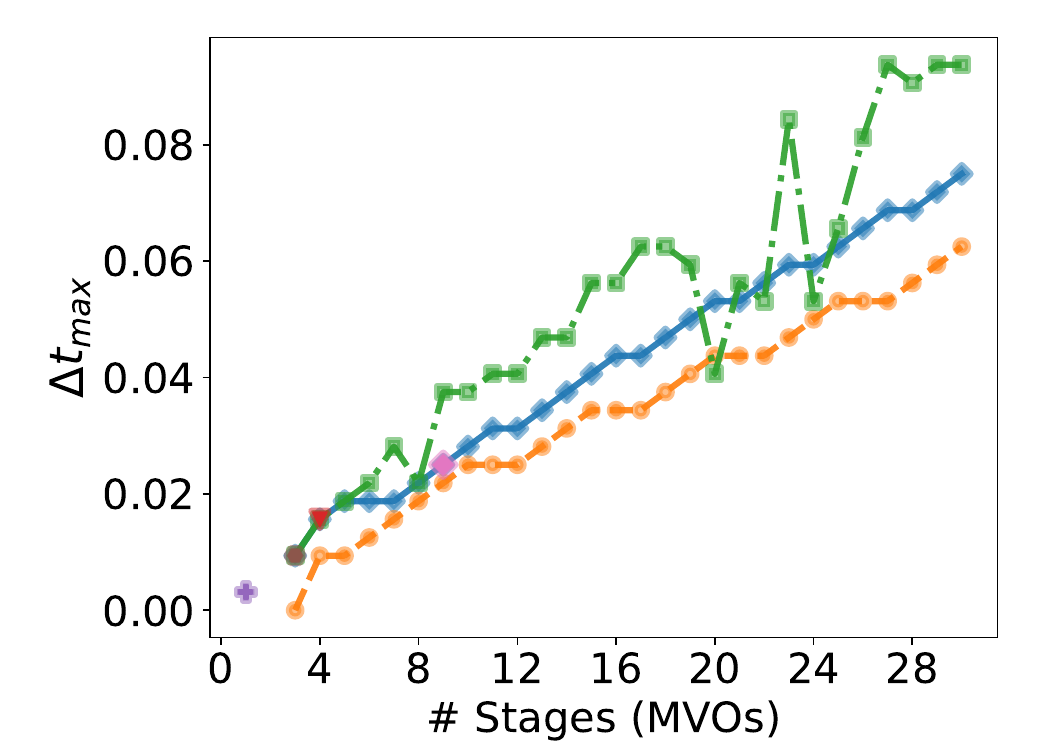}}
    \caption{Maximum time step such that the time dissipation error of each method is less than 50\% of the spatial dispersion error for different peak frequencies of the Ricker wavelet.
    In general, more stages allow larger time steps, except for the Krylov method, where $\Delta t_{\text{max}}$ reach a limit.
    }\label{fig_delta_t_diss}
\end{figure}

A similar trend of Fig. \ref{fig_delta_t_disp} is observed in Figure \ref{fig_delta_t_diss}, as with the dispersion error. The Krylov method still has the worst performance for the different peak frequencies. However, it is noteworthy that the RK9-7 method (red triangle) displays an even better performance concerning the dissipation error.

Regarding computational efficiency in the analysis of the dispersion error, the RK9-7 scheme still maintains an efficient computational performance. The FA and HORK exhibit similar behavior among the high-order methods, with a decline in efficiency for high-order Faber polynomials as the peak frequency increases. Nonetheless, the Krylov method exhibits the best performance in general, but with a very marked oscillatory behavior.

\begin{figure}[H]
\centering
\includegraphics[trim=50 480 0 0,clip,scale=0.33]{figures/legend_b.pdf}\\[-3ex]
\subfloat[Peak frequency $f_M=10$Hz.]{\includegraphics[scale=0.33]{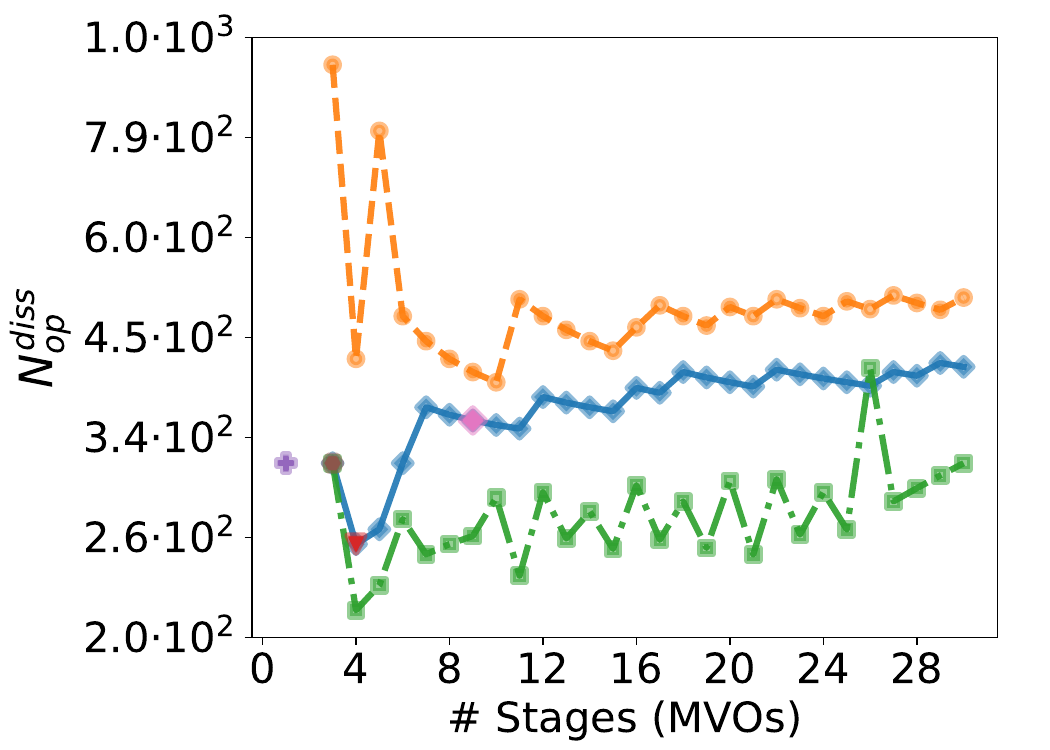}}
    \hfill
    \subfloat[Peak frequency $f_M=15$Hz.]{\includegraphics[scale=0.33]{figures/methods_ord_8_dx_0.01_f0_15_dissipation1_eff.pdf}}
    \hfill
    \\
    \subfloat[Peak frequency $f_M=20$Hz.]{\includegraphics[scale=0.33]{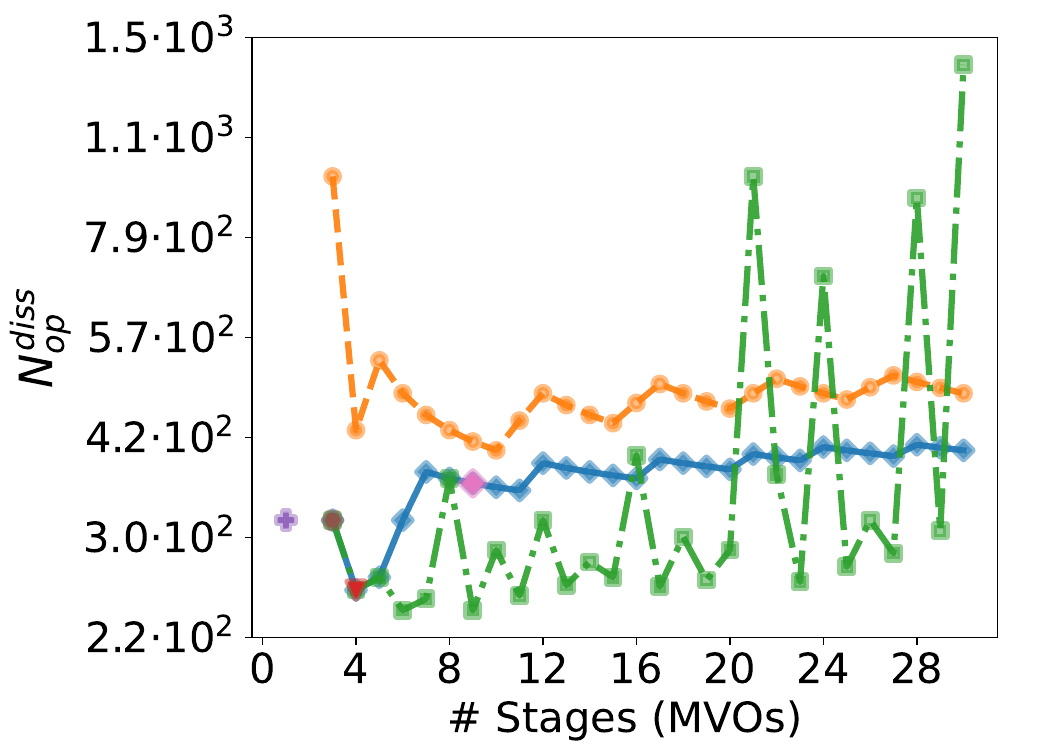}}
    \hfill
    \subfloat[Peak frequency $f_M=25$Hz.]{\includegraphics[scale=0.33]{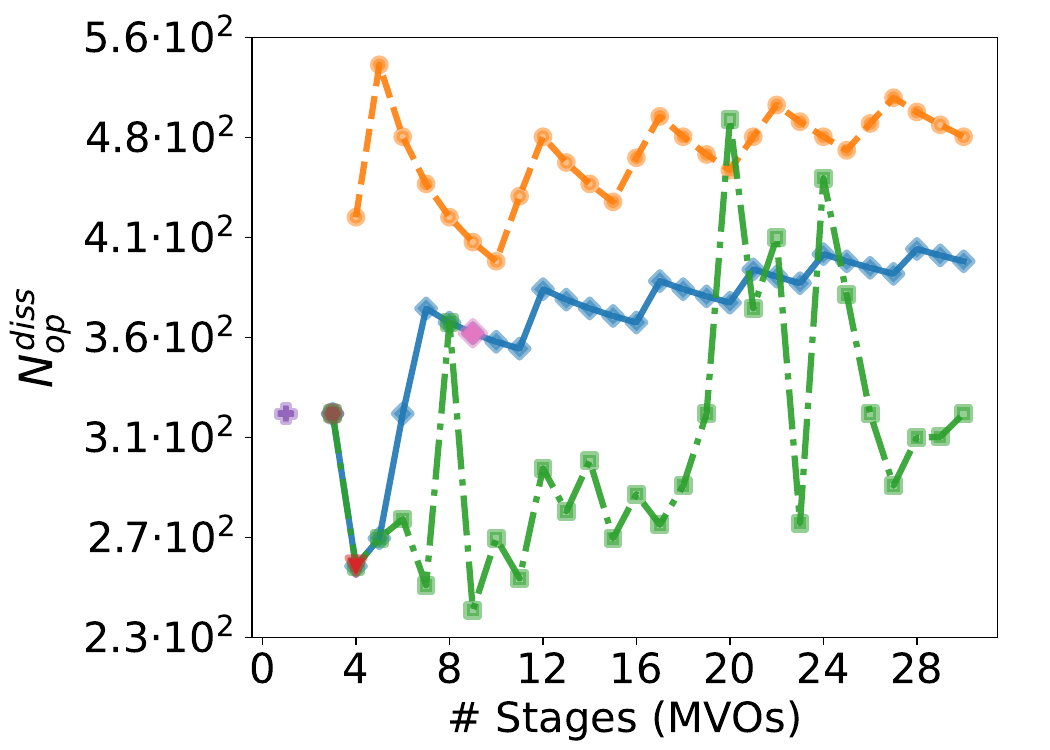}}
    \caption{Dependence on the polynomial degree of the number of matrix-vector operations by maximum time-step required to compute the solution for different peak frequencies. When the number of stages increases, the number of computations increases slightly. * Here we neglect the computational complexity of creating the Krylov subspaces.
    }\label{fig_dissipation_eff}
\end{figure}

\section{Convergence and computational efficiency}\label{sec_appendix_convergence}

In this section, we complement the results of the numerical experiments of Section \ref{sec_convergence}. First, we show the error graphics using the minimum time-step of $\Delta t=\frac{\Delta x}{8c_{\text{max}}}$, where $c_\text{max}$ is the medium maximum velocity. These graphs account for all the methods discussed in Section \ref{sec_methods} and several approximation degrees for the high-order schemes. Following that, we present the graphics of the estimation of $\Delta t_\text{max}$, the computational efficiency, and the memory utilization.

\begin{figure}[H]
\centering
\includegraphics[trim=50 480 0 0,clip,scale=0.33]{figures/legend.pdf}\\[-2ex]
\subfloat[Corner Model.]{\includegraphics[scale=0.32]{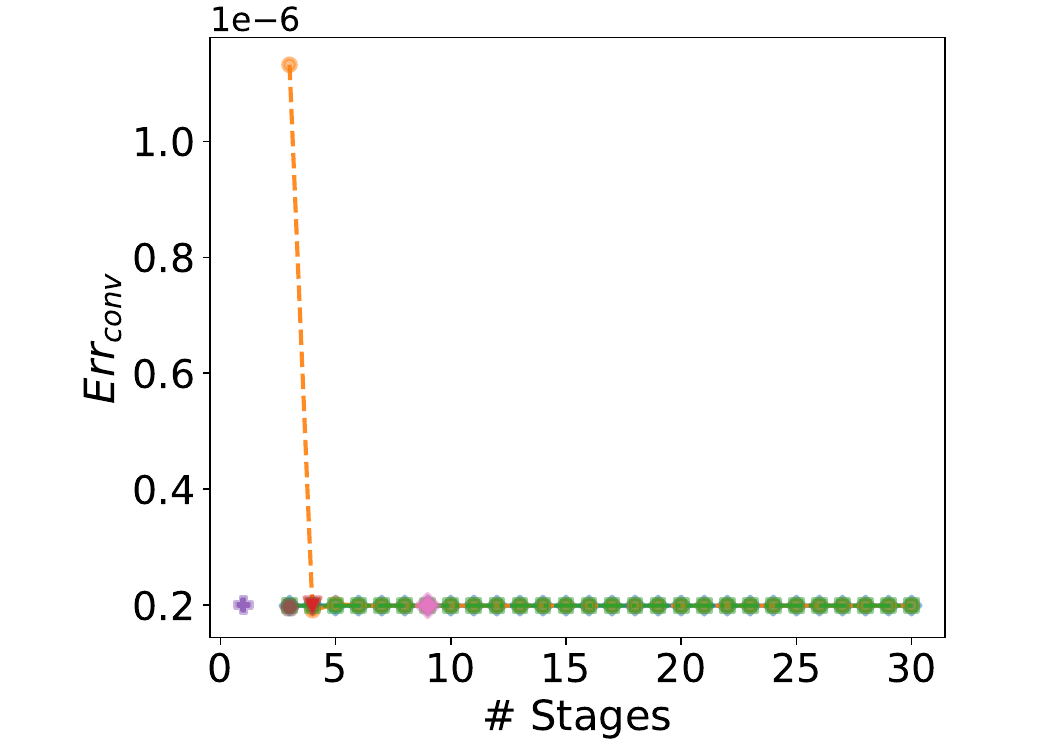}}
    \hfill
    \subfloat[Santos Basin.]{\includegraphics[scale=0.32]{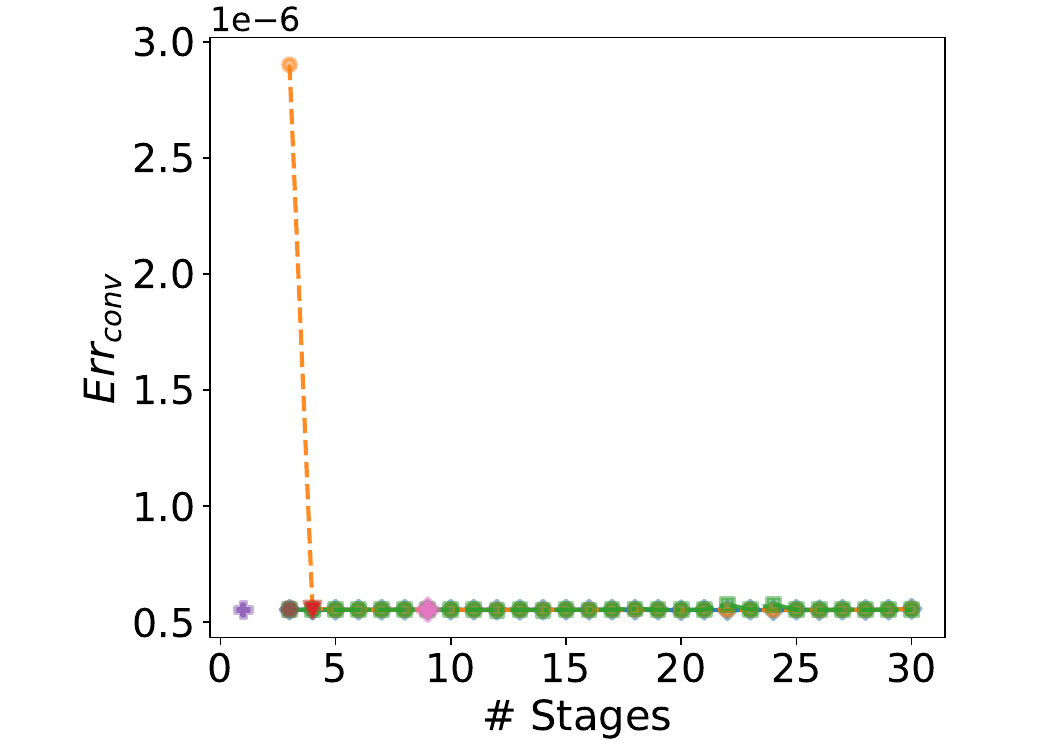}}
    \hfill
    \\
    \subfloat[Marmousi.]{\includegraphics[scale=0.32]{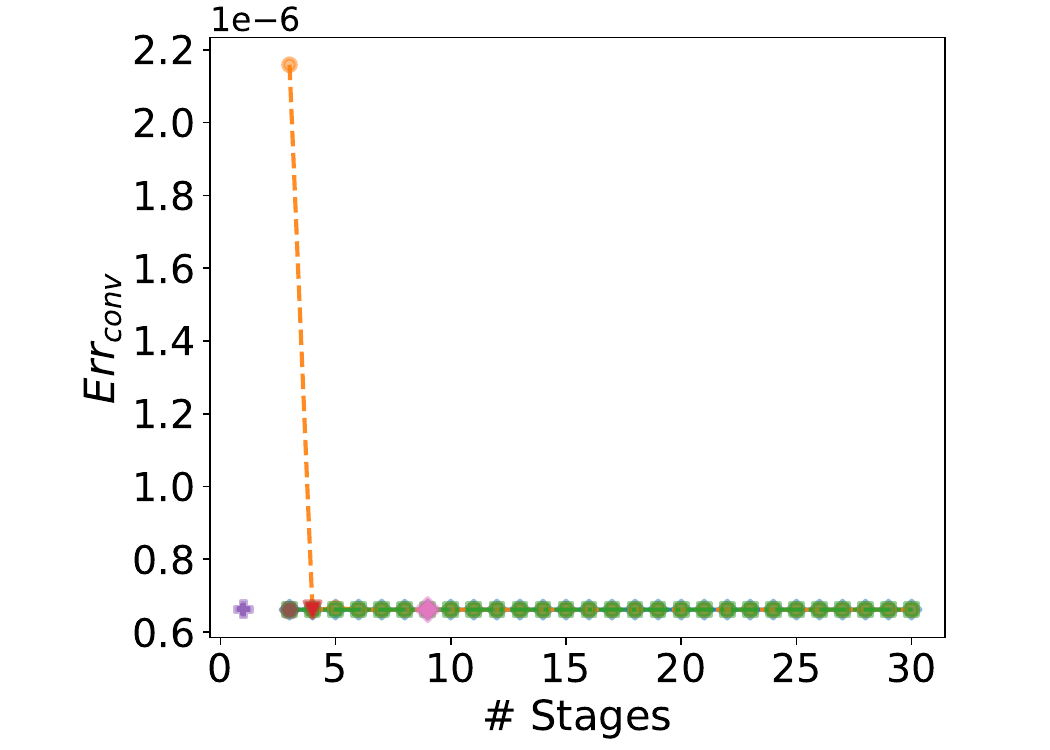}}
    \hfill
    \subfloat[SEG/EAGE.]{\includegraphics[scale=0.32]{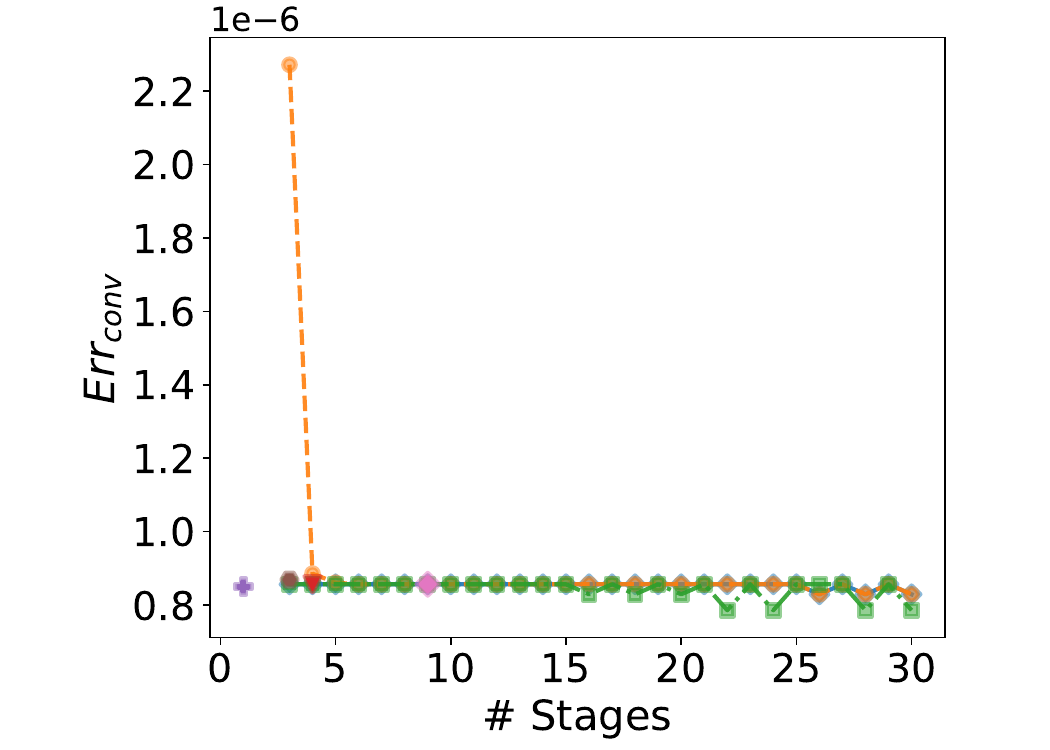}}
    \caption{Error at a time instant in the physical space achieved by each time integrator and several approximation degrees, for all the numerical experiments described in Section \ref{sec_test_cases}, using a time step size of $\Delta t=\frac{\Delta x}{8c_\text{max}}$. Regardless of the order of the method, there is an inferior limit for the error due to spatial discretization step-size size and scheme. }\label{fig_min_delta_t}
\end{figure}

Based on Figure \ref{fig_min_delta_t}, we observe an approximation error in all the numerical examples that do not decrease with the order of the method or with the selected method. This error is independent of the time integration strategy and is produced by the spatial discretization operator. While the dependence of the spatial error on the numerical experiment is weak, it is important to estimate it accurately for a reliable computation of $\Delta_\text{max}$, as quantified in Table \ref{tab_min_delta_t}.

\begin{table}[!htb]
    \centering
    \begin{tabular}{|c|c|c|}
    \hline
    Numerical experiment&Spatial error&Error tolerance\\\hline
         Corner Model& $2.02\cdot 10^{-7}$&$3.03\cdot 10^{-7}$ \\\hline
         Santos Basin& $5.55\cdot 10^{-7}$& $8.33\cdot 10^{-7}$\\\hline
         Marmousi &$6.62\cdot 10^{-7}$ &$9.93\cdot 10^{-7}$ \\\hline
         SEG/EAGE &$8.65\cdot 10^{-7}$ & $1.3\cdot 10^{-6}$\\\hline
    \end{tabular}
    \caption{Numerical error at a time instant in the physical domain produced by the spatial discretization.}
    \label{tab_min_delta_t}
\end{table}

Table \ref{tab_min_delta_t} contains two key columns of information. The first column, labeled ``Spatial error", represents the error stemming from the spatial discretization. Meanwhile, the second column, labeled ``Error tolerance", accounts for the error tolerance of 150\% of the spatial error we defined for the numerical experiment.

For the minimum error using the seismogram data, we have the respective error graphics and tolerance for each numerical test.

\begin{figure}[H]
\centering
\includegraphics[trim=50 480 0 0,clip,scale=0.33]{figures/legend.pdf}\\[-2ex]
\subfloat[Corner Model.]{\includegraphics[scale=0.32]{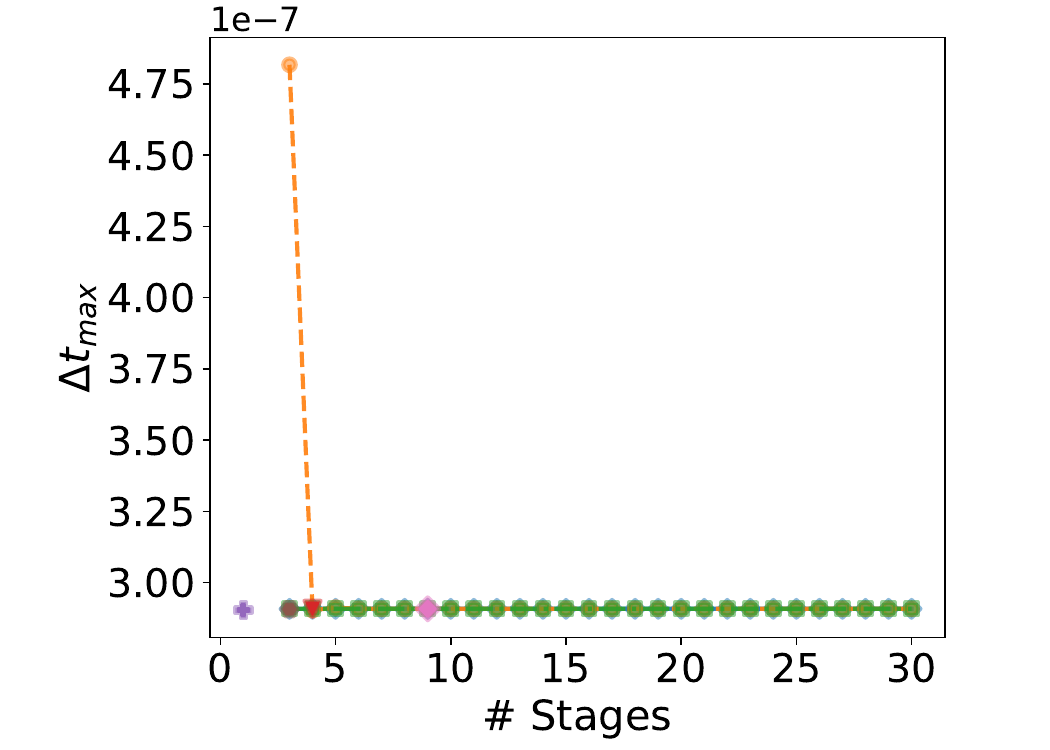}}
    \hfill
    \subfloat[Santos Basin.]{\includegraphics[scale=0.32]{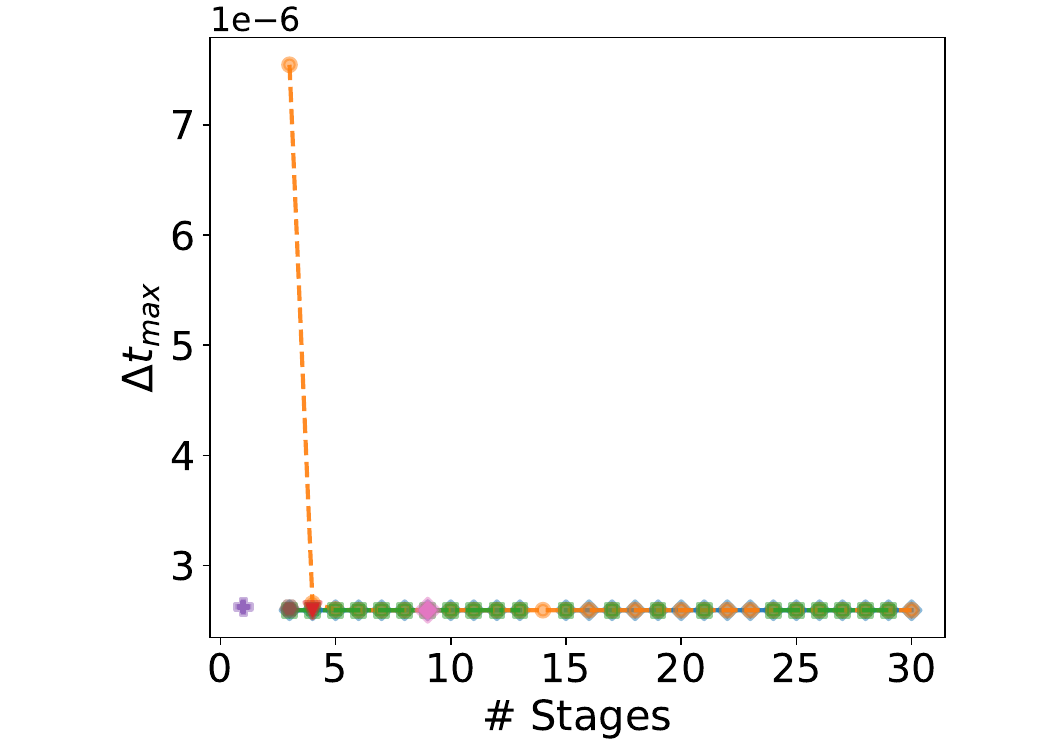}}
    \hfill
    \\
    \subfloat[Marmousi.]{\includegraphics[scale=0.32]{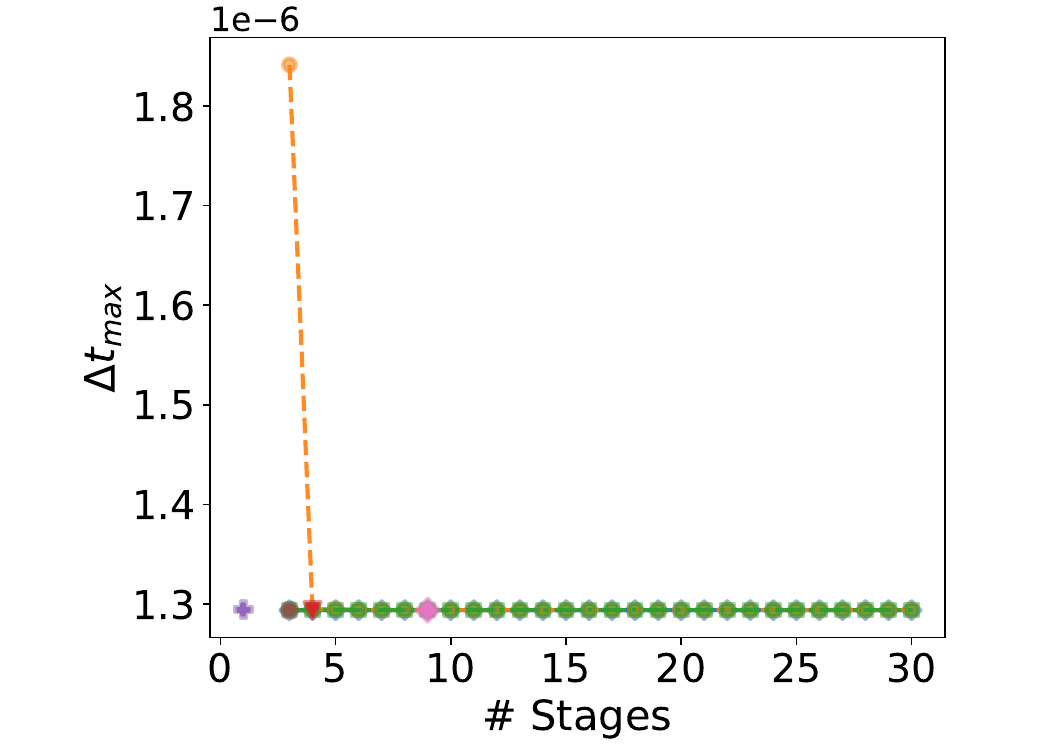}}
    \hfill
    \subfloat[SEG/EAGE.]{\includegraphics[scale=0.32]{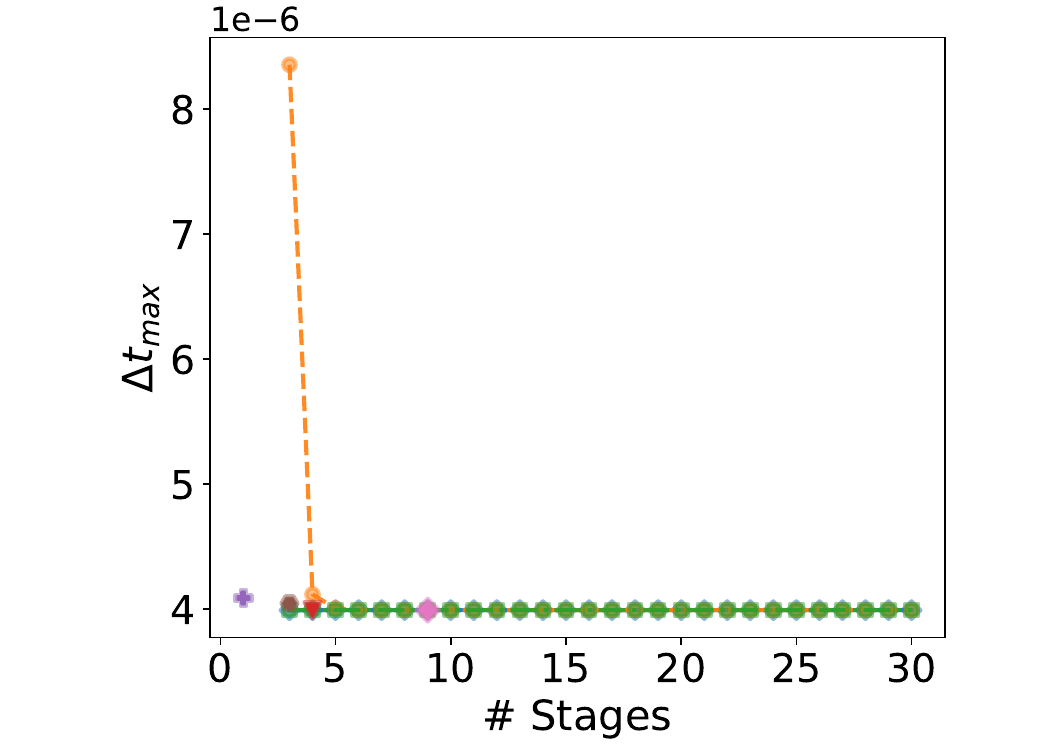}}
    \caption{Error using the seismogram data achieved by each time integrator and several approximation degrees, for all the numerical experiments described in Section \ref{sec_test_cases}, using a time step size of $\Delta t=\frac{\Delta x}{8c_\text{max}}$. Regardless of the order of the method, there is an inferior limit for the error due to spatial discretization step-size size and scheme. }
\end{figure}

\begin{table}[!htb]
    \centering
    \begin{tabular}{|c|c|c|}
    \hline
    Numerical experiment&Spatial error&Error tolerance\\\hline
         Corner Model& $2.92\cdot 10^{-7}$&$4.38\cdot 10^{-7}$ \\\hline
         Santos Basin& $2.65\cdot 10^{-6}$& $3.97\cdot 10^{-6}$\\\hline
         Marmousi &$1.3\cdot 10^{-6}$ &$1.95\cdot 10^{-6}$ \\\hline
         SEG/EAGE &$4.2\cdot 10^{-6}$ & $6.3\cdot 10^{-6}$\\\hline
    \end{tabular}
    \caption{Numerical error utilizing the seismogram data produced by the spatial discretization.}
\end{table}

\subsection{Computational efficiency and memory consumption}\label{sec_appendix_eff}

Figure \ref{fig_eff_memory} displays each time the integrator's computational cost and memory utilization for the numerical tests Corner Model, Santos Basin, and SEG/EAGE. Although there are some variations between the experiments, the general behavior remains consistent.
High-order methods require significantly less memory; in some cases, they are competitive with low-order methods, such as the Leap-Frog scheme.

\afterpage{
\begin{figure}[p]
\centering
\includegraphics[trim=50 480 0 0,clip,scale=0.33]{figures/legend_b.pdf}\\[-3ex]
\subfloat[Computational cost of Corner Model.]{\includegraphics[scale=0.3]{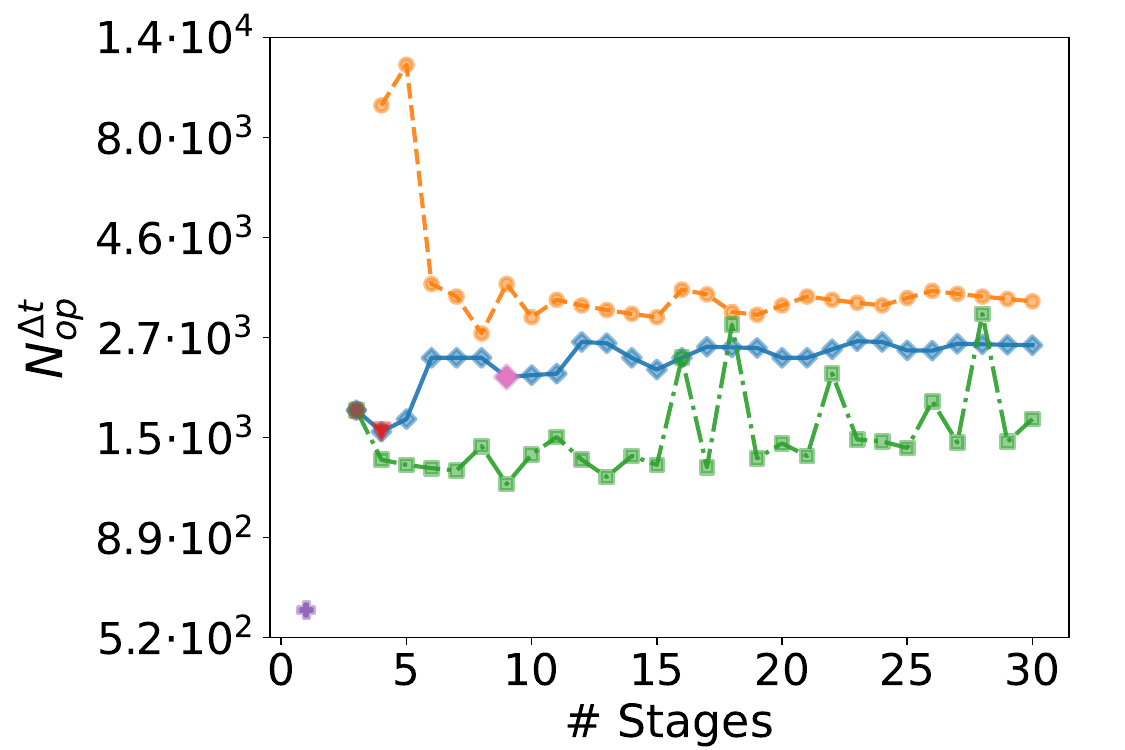}}
    \hfill
    \subfloat[Memory utilization of Corner Model.]{\includegraphics[scale=0.3]{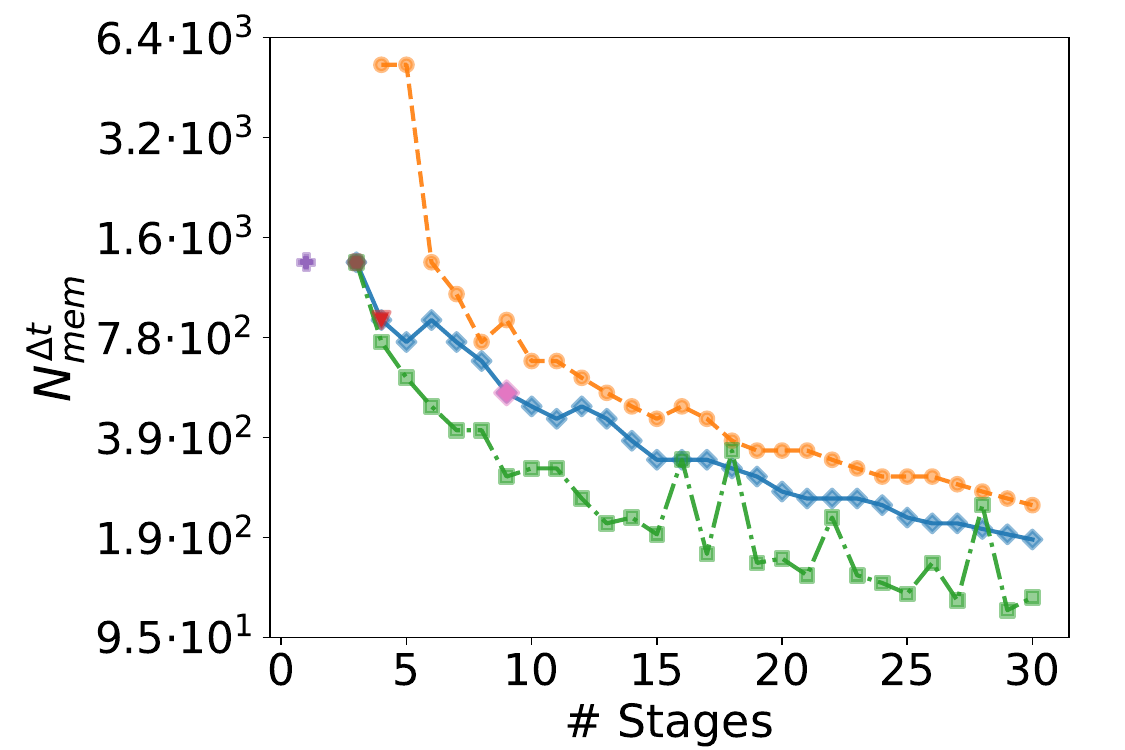}}
    \hfill
    \\
    \subfloat[Computational cost of Santos Basin.]{\includegraphics[scale=0.26]{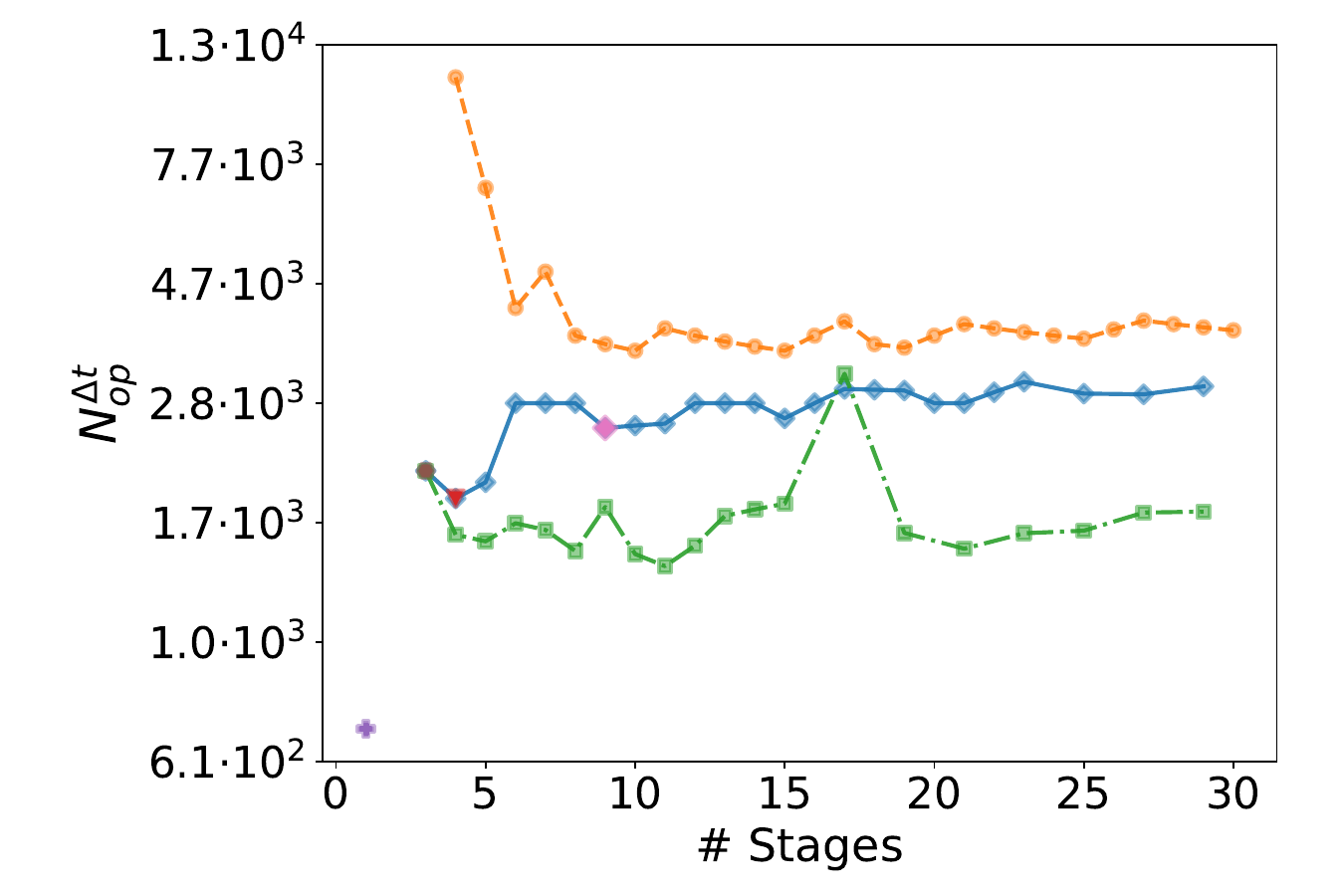}}
    \hfill
    \subfloat[Memory utilization of Santos Basin.]{\includegraphics[scale=0.26]{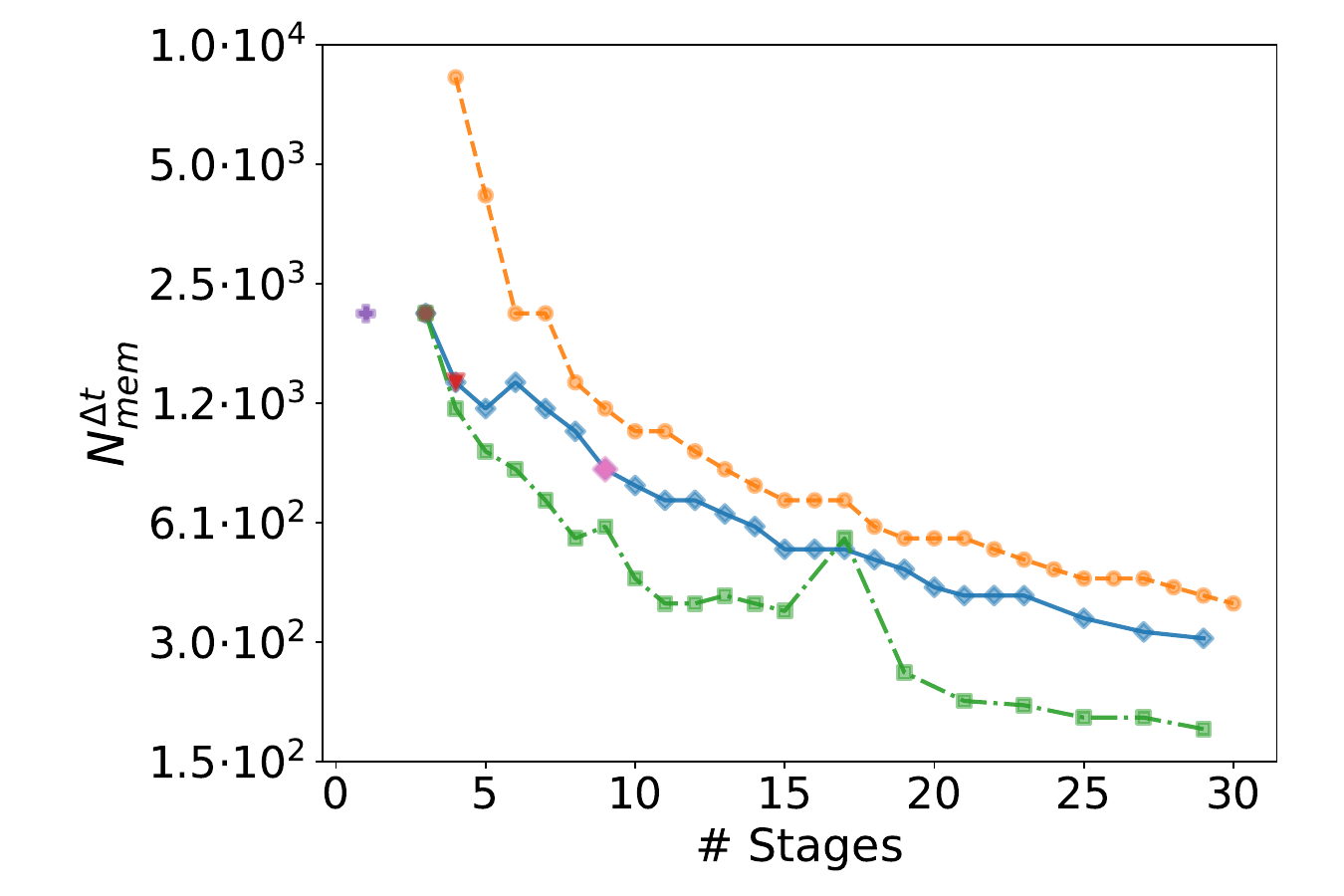}}\\
    \subfloat[Computational cost of SEG/EAGE.]{\includegraphics[scale=0.26]{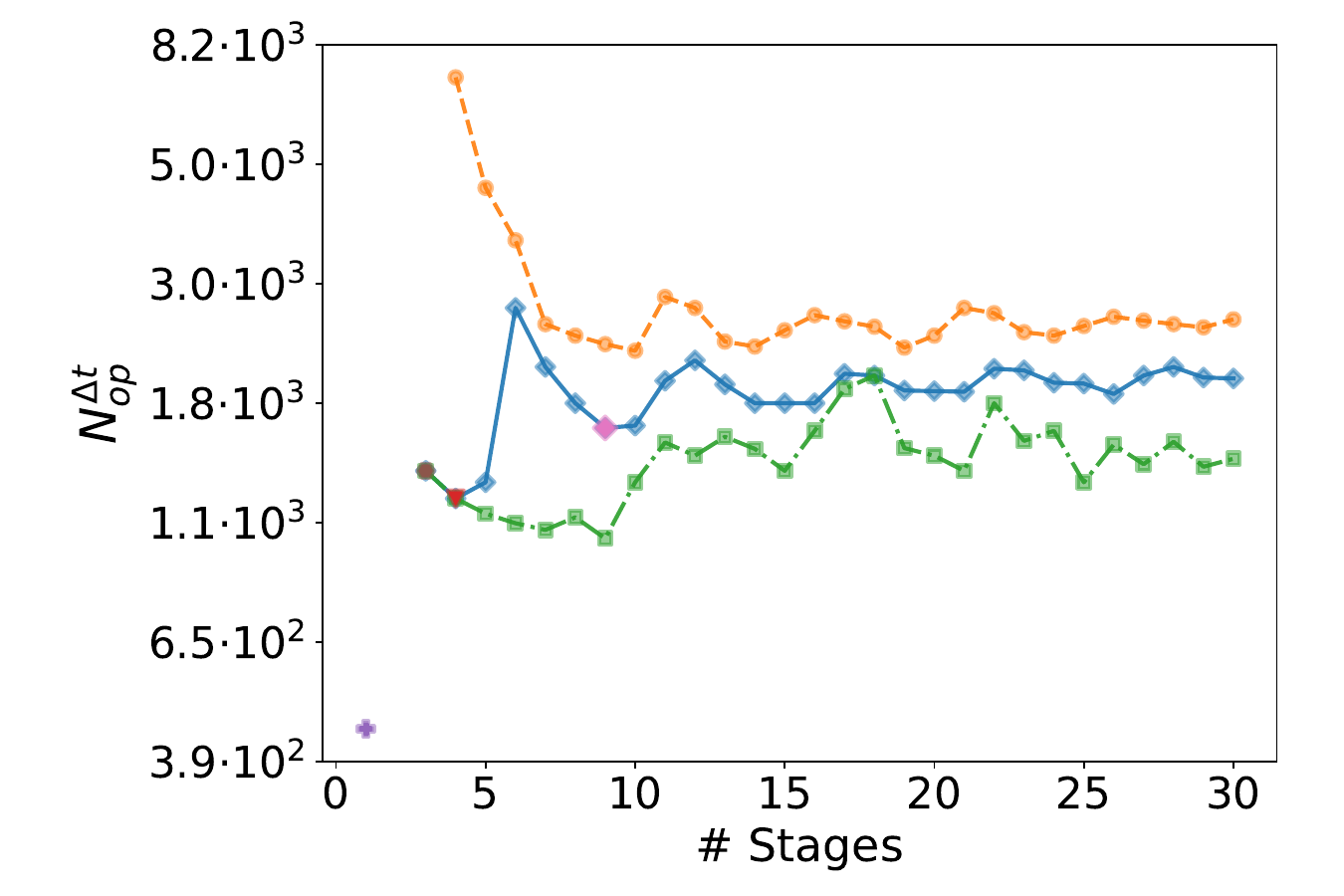}}
    \hfill
    \subfloat[Memory utilization of SEG/EAGE.]{\includegraphics[scale=0.26]{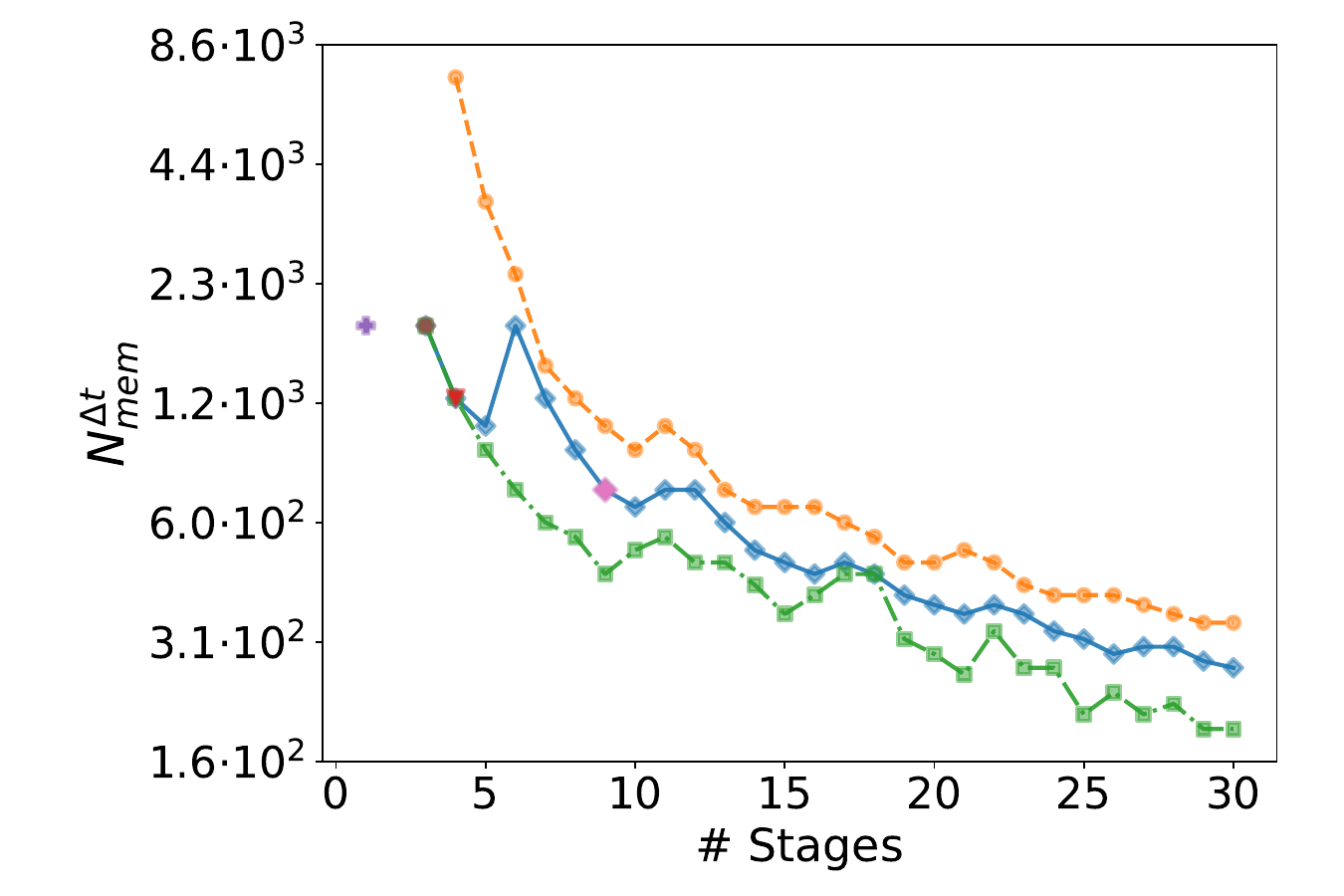}}
    \caption{Dependence of the number of MVOs and amount of stored solution vectors on the polynomial degree, for the Corner Model (first line), Santos Basin (second line), and SEG/EAGE (third line) numerical tests. As the number of stages increases, the number of computations stabilizes, and memory usage decreases. * Here we neglect the computational complexity of creating the Krylov subspaces.}\label{fig_eff_memory}
\end{figure}
\clearpage}

The relationship between the number of MVOs and the quantity of stored solution vectors concerning the polynomial degree is illustrated for the Corner Model (first line), Santos Basin (second line), and SEG/EAGE (third line) numerical tests. As the number of stages increases, there is a stabilization in the number of computations, and memory usage decreases.

\end{appendices}

%%===========================================================================================%%
%% If you are submitting to one of the Nature Portfolio journals, using the eJP submission   %%
%% system, please include the references within the manuscript file itself. You may do this  %%
%% by copying the reference list from your .bbl file, paste it into the main manuscript .tex %%
%% file, and delete the associated \verb+\bibliography+ commands.                            %%
%%===========================================================================================%%
\newpage
\bibliography{sn-bibliography}% common bib file
%% if required, the content of .bbl file can be included here once bbl is generated
% \input sn-article.bbl

\end{document}